\documentstyle[11pt]{article}

\newtheorem{theorem}{Theorem}
\newtheorem{proposition}{Proposition}
\newtheorem{lemma}{Lemma}
\newtheorem{corollary}{Corollary}

\newenvironment{proof}{\bigskip \noindent
         {\bf Proof.}}{ \hfill $\Box$ \medskip}

\newcommand{\at}{A_t}

\newcommand{\cAN}{{\cal A}^{{\bf \small N}}}
\newcommand{\cAZ}{{\cal A}^{{\bf \small Z}}}
\newcommand{\cA}{{\cal A}}
\newcommand{\cAp}{{\cal A}_{+}}
\newcommand{\cAm}{{\cal A}_{-}}

\newcommand{\cB}{{\cal B}}
\newcommand{\overb}{\overline{b}}

\newcommand{\cD}{{\cal D}}

\newcommand{\cE}{{\cal E}}
\newcommand{\cG}{{\cal G}}
\newcommand{\cGd}{{\cal G}_d}
\newcommand{\cF}{{\cal F}}

\newcommand{\cK}{{\cal K}}
\newcommand{\cL}{{\cal L}}

\newcommand{\cO}{{\cal O}}
\newcommand{\cP}{{\cal P}}

\newcommand{\overps}{\overline{\Phi}_{1,s}}  
\newcommand{\cQ}{{\cal Q}}
\newcommand{\cR}{{\cal R}}
\newcommand{\overR}{\overline{{\cal R}}}
\newcommand{\cS}{{\cal S}}
\newcommand{\cT}{{\cal T}}
\newcommand{\cTd}{{\cal T}_{2d}}

\newcommand{\cY}{{\cal Y}}

\newcommand{\bR}{{\bf R}}

\newcommand{\bP}{{\bf P}}

\newcommand{\bfl}{{\bf \lambda}}

\begin{document}

\title{{ \bf ANISOTROPIC CONTACT PROCESS ON HOMOGENEOUS TREES}}
\author{Irene Hueter} 
\date{}
\footnotetext[1]{Date: May 15, 1998. 
 Revised: July 7, 1999; April 19, 2001}
\footnotetext[2]{Department of Mathematics, University
                 of Florida, PO Box 118105, Gainesville, FL 32611-8105,
                 USA, email: {\sf hueter@math.ufl.edu}.} 
 \footnotetext[3]{Mathematics Subject Classification:
{\em Primary:} 60K35, 60J80 {\em Secondary:} 28A80.}
\footnotetext[4]{Key words and phrases: 
         Anisotropic contact process,
         Hausdorff dimension,
         homogeneous tree,
         phase transition,
         weak survival}

\maketitle

\begin{abstract}
The existence of a weak survival region is established for the
{\em anisotropic} symmetric contact process on a homogeneous tree 
$\cT_{2d}$ 
of degree $2d \geq 4:$ For parameter values in a certain connected region
of positive Lebesgue measure, the population survives
forever with positive probability but ultimately vacates every
finite subset of the tree with probability one. In this phase, 
infection trails must converge to the {\em geometric
boundary} $\Omega$ of the tree. The random subset $\Lambda$ of the
boundary consisting of all ends of the tree in which the infection
survives, called the limit set of the process, is shown to have 
Hausdorff dimension no larger than one half the Hausdorff dimension
of the entire geometric boundary. In addition, there is {\em strict
inequality} at the transition between weak and strong survival
except when the contact process is {\em isotropic}.
It is further shown that in all cases there is a distinguished
probability measure $\mu,$ supported by $\Omega,$ such that
the Hausdorff dimension of $\Lambda \cap \Omega_{\mu},$ where
$\Omega_{\mu}$ is the set of $\mu$-generic points of $\Omega,$
converges to one half the Hausdorff dimension of $\Omega_{\mu}$
at the phase separation points. Exact formulae for the Hausdorff
dimensions of $\Lambda$ and $\Lambda \cap \Omega_{\mu}$ are obtained.
We also prove that the contact process at the transition
between extinction and weak survival does not survive.
The method developed shows that the contact process
at the phase transition to strong
survival survives {\em weakly} for $d \geq 2.$
\end{abstract}


\section{Introduction}

\subsection{Background: Transition from Weak to Strong Survival}
This paper considers the anisotropic contact process on an infinite
homogeneous tree when in weak survival or at the transition to 
strong survival. The process was introduced in \cite{pema} and
pursued intensely in the isotropic case, as briefly surveyed
in \cite{lig2} (see \cite{harr,lig1} for the contact process in general).
In fact, the contact process is a stochastic growth process which,
along with branching random walks and percolation processes in spaces
with hyperbolic geometries, exhibits an intermediate phase not
present in the corresponding processes living in spaces with
Euclidean geometry. This is the {\em weak survival region},
in which the ``population" survives forever with positive probability
but, with probability one, eventually vacates every compact subset of the 
ambient space. Weak survival is known to arise for the isotropic 
contact process on a homogeneous tree \cite{pema,lig4,stac}, 
for the anisotropic symmetric branching random  walk on a 
homogeneous tree \cite{lahu}, branching Brownian motion in the 
Poincar\'{e} plane \cite{lase1}, and site percolation on 
a co-compact Fuchsian group \cite{lal1}.
The transition from weak to strong
survival still eludes complete understanding for many particle
systems.

This paper establishes the {\em existence} of a weak survival region
and analyzes its features for the {\em anisotropic} symmetric 
contact process on an infinite homogeneous tree
of even degree \mbox{$\geq 4.$} For ease of exposition, we shall
restrict our attention to the contact process on the infinite
homogeneous tree $ \cT_{2d}$ of degree $2d$ for $d \geq 1.$ 
Applying our technique to the anisotropic nearest 
neighbour contact process with some symmetry assumption
on homogeneous trees of odd degree is more subtle 
(see Remark (4) at the end of the Introduction) 
and leads to slightly different algebra.

The anisotropic case poses some genuine difficulties
not present in the isotropic case. 
The approach pursued here is completely different from 
the existence proofs for weak survival as previously 
given for the isotropic process.
For instance, the proof in this paper does not use the weight function
and the fact that the contact process becomes extinct 
at the transition to survival (even though the latter result will
be proved subsequently).
Our main interest is to study the {\em limit set} $\Lambda,$ 
defined to be the subset of the geometric boundary 
$\Omega$ (the set of ends) of the tree
in which the infection survives, when the contact process is
in weak survival and to explore the behaviour of the 
transition between weak and strong survival. 
As our explicit formulae show, the Hausdorff
dimension $\delta_H(\Lambda)$ of the limit set never exceeds
$\frac{1}{2}$ the Hausdorff dimension $\delta_H(\Omega)$ of the
geometric boundary, and importantly, the equality $\delta_H(\Lambda)
= \frac{1}{2} \delta_H(\Omega)$ is valid exactly at the transition
between both survival regions if and only if the contact process is
{\em isotropic.} This confirms a conjecture, raised in \cite{lase2},
that there is equality in the isotropic case. Moreover, we
will prove that, at the transition between extinction
and weak survival, the contact process on  $ \cT_{2d}$ 
does not survive for $d>1.$ 
Our method shows that at the phase transition to strong
survival, the contact process on $ \cT_{2d}$ 
survives {\em weakly} for $d>1.$

Additionally, we will investigate 
(a) the exponential rate $\eta$ of decay in time $t$ of the
probability that the initial infected site is infected
at time $t,$ \, 
(b) the exponential rate of growth
in space-time on the event of survival, \, 
(c) the distances of the nearest and furthest infected vertices from
the root vertex at time $t$ on the event of survival, and \, 
(d) the distribution of the limiting points in weak survival.
Key ingredients to our analysis 
are {\em shift-invariant probability measures} supported by $\Omega$ 
and their associated transition matrices related to a stationary
one-step Markov chain on some finite set of generators.
These probability measures arise from ``normalizing" the infection
probabilities at each large distance from the root vertex as the
infection is moving off to the boundary of the tree.

This paper is an expanded version of {\bf talks} that I presented
at an AMS Special Session in Gainesville, Florida (March 12-13, 1999),
at a Probability Meeting at Colorado Springs (May 28-30, 1999),
at the World Congress of the Bernoulli Society and IMS Meeting
in Guanajuato, Mexico (May 15-20, 2000),
at an AMS Special Session in New Orleans (January 13, 2001),
at a DIMACS/DIMATIA Workshop at Rutgers University, New Jersey 
(March 19-21, 2001), and
in Colloquia or Seminars at
the New College of the University of South Florida, Sarasota
(October 1, 1998),
at the University of Wisconsin at Milwaukee (November 13, 1998),
at Duke University, Durham, North Carolina (September 15, 2000), and
at the University of Berne, Switzerland (April 23, 2001).


\subsection{Anisotropic Contact Process}

We shall restrict our attention to an anisotropic, symmetric contact
process on a homogeneous tree $\cTd $ of even degree $2d$ ($d \geq 1$).
An {\em anisotropic contact process} on the tree
$ \cT = \cTd $ is a continuous time Markov process $A_t$
on the set of finite subsets of (the vertex set of) $ \cT $ that evolves
as follows. 
Infected sites (members of $A_t$) recover at rate $1$ and
upon recovery are removed from $A_t.$ Healthy sites (members of $A_t^c$)
become infected at a rate that equals the sum of the infection rates 
attached to the edges leading to infected nearest neighbours
and upon infection are added to $A_t.$
Under the default probability measure $ P,$ the initial state $A_0$
is the singleton set $\{1 \} $ (where $1$ is a distinguished element of 
$\cT,$ called the root). Each vertex $x$ of $\cT$
has exactly $2d$ neighbours. The tree $\cT$  is {\em homogeneous}
in that for any two vertices $x$ and $y$ there is an isometry that
maps $x$ to $y.$ Associated with each of the $2d$ emanating edges is an
infection rate. The symmetry assumption guarantees that each infection
rate be used twice for the set of emanating edges of each vertex.
In particular, the same infection rate is attached to an edge when
the ``infection crosses the edge" forwards and backwards.
Hence, in the notion explained in Section 1.3 below, the
set of infection rates is $\lambda_{a_1}, \ldots, \lambda_{a_d},
\lambda_{a_1^{-1}}, \ldots, \lambda_{a_d^{-1}},$ and for each of
the $d$ letters $j,$ we assume $\lambda_{j} = \lambda_{j^{-1}}.$
If the infection rates are all equal, then the contact process is called
{\em isotropic.}
Note that, alternatively, in considering the question of existence
of weak survival for the non-isotropic contact process
one may vary the recovery rates instead of the infection rates
to define an anisotropic process. 


\subsection{The Tree as a Cayley Graph}

Representing the tree ${\cal T} = {\cal T}_{2d} $ as the Cayley graph 
$\cG = \cGd$ on $d$ generators proves its worth to understand
the anisotropic contact process on the tree ${\cal T}.$ 
Let $\cAp= \{a_1, a_2, \ldots, a_d \} $
be a set of $d$ letters, let $\cAm = \{ a_1^{-1}, a_2^{-1}, \ldots,
a_d^{-1} \} $ be the set of formal inverses of the letters in $\cAp,$ and
set $\cA = \cAp \cup \cAm.$ The free group $\cG$ with generators $\cAp$
is the set of finite reduced words from the alphabet $\cA$ (a word is 
{\em reduced} if no letter $a \in \cA$ is adjacent to its inverse),
where multiplication is concatenation followed by reduction and the
group identity $1$ is the empty word. There is a natural bijection
between $\cG$ and the set of vertices of $\cT,$ in which $g,h \in \cG$ are
mapped to adjacent vertices of $\cT$ if and only if $gh^{-1} \in \cA.$
In other words, vertices are uniquely represented by finite reduced
words from $\cA.$ In the subsequent discussion,
we shall not be careful to distinguish
between vertices of $\cT$ and the words (or group elements) representing
them, and we shall refer to $\cG$ as the vertex set of $\cT.$
For any vertex $z,$ denote by $\vert z \vert $ the length of its 
representative word. Note that $\vert z \vert $ as well is the distance
from vertex $z$ to vertex $1$ in the graph $\cT.$
For every integer $n \geq 0,$ let $\cG_n$ denote the set of all
vertices $x \in \cG$ at distance $n$ from the root vertex (i.e.\ 
$\vert x \vert = n$).

In a canonical way, the bijection between
$\cT$ and $\cG$ induces a bijection between the
natural boundary of $\cT$ and the {\em geometric boundary}
$\Omega$, the set
of semi-infinite reduced words from the alphabet $\cA.$ 
A {\em geodesic} in $\cT$ is a finite or semi-infinite sequence
of distinct vertices $v_1, v_2, \ldots $ such that for every $i \geq 1,$
the vertices $v_i$ and $ v_{i+1} $ are nearest neighbours. 
An {\em end} of $\cT$ is an equivalence class of semi-infinite geodesics,
two geodesics being equivalent if and only if the sets of 
vertices through which they pass differ in at most finitely many 
vertices. 
If $\omega = x_1 x_2 \ldots \in \Omega$ then $\omega$ corresponds to
the end of $\cT$ represented by the semi-infinite geodesic that passes
through the vertices $1, x_1, x_1 x_2, \ldots $ in succession. 
For each real number $\alpha \in (0,1),$ there is a natural metric
$d_{\alpha}$ on $\Omega,$ defined by
\begin{equation}
  \label{alphametric}
     d_{\alpha}(\omega, \omega') = \alpha^{N(\omega, \omega')},
\end{equation}
where $N(\omega, \omega')$ is the largest integer $n$ such that
the sequences $\omega$ and $\omega'$ agree in entries $1,2, \ldots,n.$
For any choice of $\alpha,$ the corresponding topology on $\Omega$ is 
the topology of coordinatewise convergence.
For any vertex $z $ of $\cT,$ define $\cT(z)$ to be the set
of vertices $v$ such that the geodesic segment from $1$ to $v$ passes
through $z, $ equivalently, 
such that the unique word representing $z$ 
is a prefix of the word representing $v.$ Similarly, define $\Omega(z)$
to be the set of infinite reduced words $ \omega = x_1 x_2 \ldots $
such that, for some finite $n, $ the word $z$ is represented by the 
word $x_1 x_2 \ldots x_n.$ Observe that, for every integer $n \geq 1,$
the set $ \{ \Omega(z): \vert z \vert = n \} $ is a finite open cover
of the geometric boundary $\Omega.$
Finally, define $\Sigma$ to be the set of all doubly 
infinite reduced words 
 $\xi = (x_n)_{n=- \infty}^{\infty}$ from $\cA.$


\subsection{Anisotropic Contact Process on $\cG$}

The symmetric contact process on a homogeneous tree of degree $3$ or
larger distinguishes itself from the symmetric process on the integer
lattice ${\bf Z}^d$ in that there are two different survival regions
\cite{pema,lig4,stac}. More precisely, on $\cT_{2d}$ for $d > 1$ 
(but not for $d=1$), 
there is a partition of the
parameter space $[0, \infty)^d = {\bf R}_+^d \ni \bfl =
( \{ \lambda_a \}_{a \in \cAp} )$ into three regions 
${\cal R}_1, {\cal R}_2$ and $ {\cal R}_3$ such that
\begin{enumerate}
\item[(a)]
if $ \bfl \in {\cal R}_1,$
then $A_t = \emptyset$ eventually, with probability $1,$
\item[(b)]
if $ \bfl \in {\cal R}_2,$
then $ P\{ \vert A_t \vert \rightarrow \infty\} > 0,$ but 
$\forall x \in {\cal T},$ \mbox{$ P \{ x \in A_t \; \mbox{for arbitrarily
large t} \} =0,$}
\item[(c)]
if $\bfl \in {\cal R}_3,$
then with positive probability $\vert A_t \vert \rightarrow \infty$
and, for all $x \in {\cal T}, $ for arbitrarily large values of $t,$
$ x \in A_t. $ 
\end{enumerate}
On ${\cal R}_1,$ the contact process is called {\em subcritical}, on
${\cal R}_2,$ 
{\em weakly supercritical}, and on ${\cal R}_3,$  
{\em strongly supercritical.} The main results of this paper concern
the weak survival regime 
\begin{equation}
      \bfl \in \cR_2,
\end{equation} 
whose existence needs to be established first. Additionally,
our findings will shed some light onto the nature of the
boundary of $\cR_2,$ that is, the phase transitions.

We begin with introducing some terminology for the anisotropic
contact process with infection rates $  \{ \lambda_a \}_{a \in \cAp}$ 
and recovery rate $1,$ in order to state the
main results. Thus, for every $x \in \cG,$ 
vertex $x$ attempts to infect vertex $x i,$
$i \in \cA,$ at rate $\lambda_i$ (note that vertex $xi$
may already be infected).
Recall that $1$ denotes the root vertex.                        
Define 
\begin{equation}
 \label{rootinf}
   \eta = \eta(\bfl) = \lim_{t \rightarrow \infty}
                      (P\{ 1 \in A_t \} )^{1/t}
   = \sup_{t > 0} (P\{ 1 \in A_t \} )^{1/t} \leq 1.
\end{equation}
The limit exists because by the Markov property (see Section 2.1) and 
monotonicity properties of the contact process, $P \{ 1 \in A_t \}
P \{ 1 \in A_s \} \leq P \{ 1 \in A_{t+s} \}$ for all real $ s,t > 0.$
Thus, an easy subadditivity argument applies. It follows as well
that $P \{ 1 \in A_t \} \leq \eta(\bfl)^t$ for all $ t> 0.$
It is obvious that the function $\eta(\cdot)$ is nondecreasing 
in each $\lambda_j.$ 
For any vertex $x \in \cT,$ define
\begin{equation}
  \label{functionux}
   u_x = u_x(\bfl) = P  \{ x \in \at \, \, \mbox{for some } \, t > 0 \}.
\end{equation}
The strong Markov property together with
the monotonicity properties of the contact process
and the homogeneity of the process at each vertex
implies that $ u_{xy} \geq u_x u_y$
for each $x,y \in \cT$ such that 
$\vert x  y \vert = \vert x \vert + \vert y \vert$
(no reduction occurs when $x$ and $y$ are concatenated).
If the contact process is weakly supercritical, then $u_x < 1 $
for every $x \not =1$ in $\cT$ because 
$u_x = 1$ would imply that, with probability
one, the root be reinfected at indefinitely large times, thus, the
process would be strongly supercritical.
A subadditivity argument shows that, for every $a \in \cA$
and each vertex $x = a a \ldots a \in \cG_n,$ the limit
\begin{equation}
  \label{betalimit}
   \lim_{\vert x \vert =n \rightarrow \infty}
   u_x(\bfl)^{1/n} = \beta_a = \beta_a(\bfl)
\end{equation}
exists and that $u_x(\bfl) \leq \beta_a(\bfl)^n$ for all $n \geq 0.$ 
Moreover, for every integer $k$ and $x \in \cG_k,$ by subadditivity
(more precisely, by supermultiplicativity),
for each {\em periodic} sequence   \linebreak
$y_n = x x \ldots x \in \cG_{nk},$ 
the limit 
\begin{equation}
  \label{uperiodicsubadd}
   \lim_{n \rightarrow \infty} u_{y_n}(\bfl)^{1/n} = \beta_x(\bfl)
\end{equation}
exists for every $\bfl$ and $ u_{y_n}(\bfl) \leq \beta_x(\bfl)^n$ 
for every $n \geq 0.$
Clearly, the functions $\beta_x(\cdot)$ are nondecreasing in each 
infection parameter $\lambda_j.$ We will discuss the strict monotonicity
properties of the $\beta_x(\bfl),$ $\eta(\bfl),$ and other functions
(Sections 2.5 and 4.1) and their continuity properties (Section 5).


\subsection{Main Results: Weak Survival Region and Limit Set} 

There has been a wealth of results on {\em isotropic} contact processes
on trees,
among them \cite{pema, dusc, msz, stac, lig3, lig2, lase2,lal2, sch1}.
In contrast to the contact process on Euclidean lattices, the isotropic
contact process on homogeneous trees was shown to exhibit two essentially 
different survival phases, the weak one being a novelty \cite{pema,lig4}.
One might suspect the {\em anisotropic} contact process
to as well have an intermediate region between extinction and strong 
survival. A result of this paper confirms this previously 
conjectured behaviour for the symmetric anisotropic 
contact process on $\cT_{2d}.$

Further detailing the boundary of the weak survival phase involves
functions \mbox{$ 0 \leq \overb_i(\bfl) \leq 1$} $(i \in \cA),$
nondecreasing in each argument $\lambda_j,$ 
which unfortunately do not come with a short description but are
defined in a number of steps and merely are intermediate tools in 
our study (as opposed to being of interest on their own).
However, to characterize the limit set at the weak/strong survival
transition, the $ \overb_i(\bfl) $ provide a crucial link 
to an analytic tool available for the symmetric anisotropic 
random walk \cite{lahu}.
The $ \overb_i(\bfl)^2$ arise as entries of a certain $2d \times 2d$
Perron-Frobenius matrix with a lead eigenvalue that coincides with
the exponential of the expectation of a potential function 
relative to a Gibbs measure. The full description is 
deferred to Sections 3 and 4. Importantly for the next result, 
each $0 < \overb_i < 1.$  

\begin{theorem}
   \label{weaksurvival}
The weak survival region $\cR_2$ is nonempty unless $d=1$ and
enjoys the following properties:
\begin{enumerate}
\item[(a)]
\vspace*{-0.15cm}
There are functions $\overb_i(\bfl),$ defined in (\ref{rhomatrix}),
such that the boundary $ \overline{\cR}_1 \cap \overline{\cR}_2$ consists
of all $\bfl$ for which the $\overb_i(\bfl)$ satisfy
$$
       \sum_{i \in \cA} \, \frac{ \overb_i(\bfl)}{
     1 + \overb_i(\bfl)} = 1
$$
and such that the boundary $ \overline{\cR}_2 \cap \overline{\cR}_3$
consists of all $\bfl$ so that 
$$
       \sum_{i \in \cA} \, \frac{ \overb_i(\bfl)^2}{
     1 + \overb_i(\bfl)^2} = 1.
$$ 
\item[(b)]
\vspace*{-0.2cm}
Every line in the interior of the first quadrant in 
${\bf R}^d$ that passes 
through the origin has an intersection with $\cR_2$ that is 
a line segment. 
\item[(c)]
\vspace*{-0.2cm}
The region $\cR_2$ has positive $d$-dimensional Lebesgue measure.
\item[(d)]
\vspace*{-0.2cm}
The region $\cR_2$ is connected and is a symmetric region 
in the $d$ parameters $\lambda_i.$
\end{enumerate} 
\end{theorem}

The {\em critical} contact process behaves as follows.

\begin{theorem}
  \label{survivalattransition}
For $\bfl \in \overR_2 \cap \overR_3,$ the contact process 
survives {\em weakly}.
\end{theorem} 

For a proof of this result, see Corollary \ref{cpatsecondphasetrans}.

\begin{theorem}
  \label{extinctionattransition}
$ \overR_1 \cap \overR_2 \subset \cR_1,$ that is,
for $\bfl \in \overR_1 \cap \overR_2,$ the contact process 
on $ \cT_{2d}$ for $d >1$ almost surely becomes extinct.
\end{theorem} 

In the isotropic case, the result is in \cite{pema} for degree $2d \geq 4$
and in \cite{msz} for degree $3.$ Moreover, 
the critical contact process on ${\bf Z}$ and in all other
Euclidean lattices dies out \cite{begr}.
Theorems \ref{weaksurvival} through \ref{extinctionattransition}
will be proven in Section 5.6.
For any set $B \subset {\bf R}^d,$
let $B^c$ denote its complement in ${\bf R}^d.$

\begin{theorem}
  \label{etaatboundary}
For all $d \geq 1,$ we have
 $\eta(\bfl) < 1$ for $\bfl \in \mbox{int}(\cR_1 \cup \cR_2),$
and,
$\eta(\bfl)=1$ for $\bfl \in \overline{\cR_3^c} \cap \overR_3.$
\end{theorem}

In the weak survival region
$\cR_2,$ the ``population'' eventually vacates 
every finite subset of vertices of the tree with probability one.
Therefore, the population has a well-defined limit set.
Define the {\em limit set} $\Lambda$ of the contact process on $\cT$
to be the (random) set of $\omega= x_1 x_2 \ldots \in \Omega$
such that each vertex $x_1 x_2 \ldots x_k$
of $\omega$ is infected at some time. It is
easily seen that, if the contact process is supercritical, then 
on the event of survival, $\Lambda$ is nonempty and compact
(relative to any of the metrics $d_{\alpha}$). 

\begin{theorem}
  \label{limitset}
For $\bfl \in \cR_2,$ almost surely on the event of survival,
the Hausdorff dimension $\delta(\bfl)$ of $\Lambda$ (relative to the
metric $d_{\alpha},$ defined in (\ref{alphametric}))
is given by 
\begin{equation}
    \label{dimlimit}
         \delta(\bfl) = -\frac{\log \theta(\bfl)}{\log \alpha}, 
\end{equation}
where $\theta(\bfl)$ is the leading eigenvalue of some 
Perron-Frobenius matrix, described in \mbox{Section 4,} and
is the unique positive number such that
\begin{equation}
 \label{betatheta}
     \sum_{i \in \cA} \, \frac{ \overb_i(\bfl)}{ \theta(\bfl) + 
                 \overb_i(\bfl)} = 1,
\end{equation}
where the $\overb_i(\bfl)$ are defined in (\ref{rhomatrix}).
For $\bfl \in \cR_2,$
the functions
$ \theta(\bfl)$ and $\delta(\bfl) $ are 
continuous functions in
each of the variables $\lambda_j$ and are strictly increasing
along ``directions of increase" for $\bfl$ (see Section 2.5).
Furthermore, if $ \delta_H(\Omega)$ denotes the Hausdorff dimension
of $\Omega,$ 
\begin{equation}
\label{halfdim1}
    \delta(\bfl) \leq \frac{1}{2} \, \delta_H(\Omega),
\end{equation}
with equality holding if and only if the underlying contact process
is isotropic and $ \bfl \in \overline{\cR}_2 \cap \overline{\cR}_3$
(at the transition between weak and strong survival).
\end{theorem}

The equality in (\ref{halfdim1}) confirms a conjecture, first raised
in \cite{lase2}, concerning the behaviour of the isotropic contact process 
at the transition between weak and strong survival.
Any function in $\bfl$ that is strictly increasing along ``directions 
of increase" for $\bfl$ necessarily is strictly increasing
if every component of $\bfl$ is increased.  
The following properties will be valuable.

\begin{theorem}
  \label{properties}
For each $i,j \in \cA,$ the functions
\begin{eqnarray*}
    \label{mainfctcont}
      \lambda_j & \rightarrow &  \overb_i(\bfl), \qquad \qquad 
       \lambda_j   \rightarrow  \beta_i(\bfl), 
                           \nonumber \\
         \lambda_j &  \rightarrow & \eta(\bfl),
            \qquad  \qquad 
        \lambda_j  \rightarrow \theta(\bfl)      
\end{eqnarray*} 
are continuous for each $\bfl \in \mbox{int}(\cR_1 \cup \cR_2),$
are left-continuous for $\bfl \not \in \cR_3$ such that 
each $\lambda_k >0.$ Moreover, on $\mbox{int}(\cR_1 \cup \cR_2),$ 
each of these functions
is strictly increasing along ``directions of increase" for
$\bfl$ as defined in Section 2.5.
\end{theorem}

Now we turn back to the limit set of the contact process.
The more intriguing part of our studies revolves around the behaviour
of the limit set $\Lambda$ when $\bfl$ takes a critical value 
in $ \overline{\cR}_2 \cap \overline{\cR}_3.$
In order to understand the fine structure of $\Lambda,$
we partition the geometric boundary $\Omega$ into
{\em measure classes} $ \Omega_{\mu} $ and ask about the size of
the intersection of each equivalence class with the limit set. 
Recall that $\Omega$ is the
set of semi-infinite reduced words from $\cA.$  Let $\sigma: 
\Omega \rightarrow \Omega$ be the one-sided forward shift operator on 
$\Omega,$ that is,
$$
  \sigma(x_1 x_2 \ldots) = x_2 x_3 \ldots. 
$$
For any ergodic, $\sigma$--invariant probability measure $\mu$ on
$\Omega,$ define $\Omega_{\mu}$ to be the subset of $\Omega$ 
consisting of all $\omega \in \Omega$ such that for every continuous
real-valued function $f: \Omega \rightarrow \bR,$
\begin{equation}
   \label{mugeneric}  
    \lim_{n \rightarrow \infty} \frac{1}{n} \sum_{i=1}^n 
     f(\sigma^i \omega) = \int_{\Omega} f d\mu.
\end{equation} 
Birkhoff's ergodic theorem implies that $\mu(\Omega_{\mu}) = 1 $
since the space of continuous functions on $\Omega$ is separable
in the sup norm topology. Moreover, if $\mu $ and $\nu$ are distinct
ergodic probability measures, then $\Omega_{\mu} \cap \Omega_{\nu}
= \emptyset. $

For any ergodic, $\sigma$--invariant probability measure $\mu$ on 
$\Omega,$ let $h(\mu)$ denote the Kolmogorov-Sinai entropy 
of the measure-preserving system $(\Omega, \mu, \sigma)$
(for the definition, see e.g.\ \cite{walt}, Chapter 4).
Define the function $\varphi_{\bfl} : \Omega \rightarrow \bR$
by 
\begin{equation}
 \label{phifunction}
 \varphi_{\lambda}(x_1 x_2 \ldots ) = \log \overb_{x_1}(\bfl), 
\end{equation}
where the functions $\overb_i(\bfl)$ are defined in (\ref{rhomatrix}).

\begin{theorem}
 \label{halfboundary}
Let $\bfl \in \cR_2$ and let $\mu$ be any ergodic, 
$\sigma$--invariant probability measure on $\Omega.$
If $ h(\mu) < - \int 
\varphi_{\bfl} d\mu, $ then almost surely,
 $\Lambda \cap \Omega_{\mu} = \emptyset. $ If 
$h(\mu) \geq 
- \int \varphi_{\bfl} d \mu,$
then almost surely on the event of survival, 
the set $\Lambda \cap \Omega_{\mu}$
has Hausdorff dimension $ \delta(\bfl;\mu)$ (relative to the metric
$d_{\alpha}$)
\begin{equation}
     \label{dimintersect}
   \delta(\bfl;\mu) =  - \frac{h(\mu) + \int_{\Omega} \varphi_{\bfl}
                          d\mu}{ \log \alpha} .
\end{equation}
This Hausdorff dimension satisfies the inequality 
\begin{equation}
   \label{halfdim2}
    \delta(\bfl;\mu) \leq \frac{1}{2} \delta_H (\Omega_{\mu}),
\end{equation}
where equality holds in (\ref{halfdim2}) for one and only one 
ergodic probability measure $\mu_*$ and only when $\bfl \in
\overline{\cR}_2 \cap \overline{\cR}_3.$
\end{theorem}                       

The second display is the more curious part of Theorem 
\ref{halfboundary}, the existence of a shift-invariant
probability measure $\mu_*$ on $\Omega$ for which equality
holds in (\ref{halfdim2}), and provides a more vigorous instance
of equality (\ref{halfdim1}), its analogue in the isotropic case.
In fact, the transition from
weak to strong survival happens precisely when, for some 
$\Omega_{\mu},$ the set $\Lambda$ fills a subset of half the
Hausdorff dimension of $\Omega_{\mu}.$ 
The distinguished probability measure
$ \mu_*$ is defined as follows. 

Define the {\em backscatter matrix} $M_2 = M_2(\bfl)$ to be the $2d 
\times 2d$ matrix, indexed by elements of $\cA,$ whose entries 
are given by 
\begin{eqnarray}
  \label{backscattermatrix}
            (M_2(\bfl))_{ij} & = & \overb_j(\bfl)^2 \quad
                        \mbox{if }  j \not = i^{-1}, \\
                            & = & 0 \qquad \, \,
                        \quad \mbox{if }  j = i^{-1}. \nonumber 
\end{eqnarray}
We will see that, for $\bfl \in \cR_2$ such that each
$\lambda_k >0,$
this is an irreducible nonnegative matrix, 
thus, the Perron-Frobenius theorem applies. 
Proposition \ref{secondtrans}
 below will show that the lead
(Perron-Frobenius) eigenvalue of $M_2$ is $1$
for every critical value $\bfl \in 
\overline{\cR}_2 \cap \overline{\cR}_3.$
It follows that, if $v$ is the (positive) right eigenvector, then
for every $\bfl \in \overline{\cR}_2 \cap \overline{\cR}_3,$
\begin{equation}
  \label{stochmatrix}
     p_2(i,j) = \frac{(M_2(\bfl))_{ij} v_j}{v_i} 
\end{equation}
are the entries of an irreducible stochastic matrix $\bP_2.$ The
probability measure $\mu_*$ is the unique probability measure on
$\Omega$ such that the induced coordinate process is the
stationary Markov chain with transition probability matrix $\bP_2.$
In fact, it is only possible that such
an invariant measure exists for each $\bfl$
at the phase transition $\overline{\cR}_2 \cap \overline{\cR}_3$
due to the fact that the critical contact process survives weakly. 
Observe that in a similar fashion, for each $\bfl \in \cR_2$ and
each positive real $\rho$ above some (critical) value,
if $M_{\rho}$ represents the matrix that results when
each power $2$ is replaced by power $\rho$ in 
(\ref{backscattermatrix}), then
$$
   p_{\rho}(i,j) = \frac{\overb_j(\bfl)^{\rho} v_j}{
          \theta(\rho;\bfl) v_i} \, ( 1 - \delta_i(j^{-1}))
$$
denotes the transition probabilities of a {\em stationary}
one-step Markov chain, where $v$ is the  right eigenvector
associated with the lead eigenvalue 
$ \theta(\rho; \cdot)$ of the 
Perron-Frobenius matrix $M_{\rho}$ and $\delta_{.}(\cdot)$ 
denotes the Kronecker delta function. For each $\bfl \in \cR_2$ 
and each such $\rho,$ there exists a {\em unique shift-invariant}
probability measure supported by $\Omega$ 
(see Section 3 for more details). 

Moreover, it will be demonstrated (Theorem \ref{statprocess}) that, if   
$\omega_1  \omega_2 \ldots $ is a limit point in $\Omega$ and if  
$\mu_n$ denotes the distribution
 under $P$ of the process  $\omega_n , \omega_{n+1}, \ldots,$ 
then for every $n \geq 1,$ the measure $\mu_n$ 
is absolutely continuous with respect to a {\em Gibbs state} 
$\mu_{\varphi}$ and $\mu_n \stackrel{\cD} {\rightarrow} \mu_{ \varphi}$
as $n \rightarrow \infty.$ The limiting distribution will be 
identified. The stationary process induced by 
$\mu_{\varphi}$ is isomorphic to a Bernoulli shift \cite{bowe}.
These measures decay
exponentially in the distance from the root vertex
and are spherically symmetric if and only if the contact
process is isotropic (see also Proposition \ref{pstatprocess}
in Section 3.5 as well as the remarks in the paragraph thereafter).
 
Our final result offers insight into the dispersal behaviour
of the infection in space-time.
The symmetric anisotropic contact process moves at linear distance 
with time as does the isotropic contact process \cite{lal2}.
Following the same notation, we define $r_t$ and $R_t$ to be
the {\em smallest} and {\em largest} distances among the infected
sites $x \in A_t$ and $N_n(ns)$ to be the number of vertices $x \in
A_{ns}$ at distance $n$ from the root that are infected at time
$ns.$
In the time-dependent case, 
let $\exp \{ \overline{\Phi}_{1,s; \bfl} \}$ 
be the analogous function to $\theta (\bfl),$ defined in 
(\ref{betatheta}), at time scale $s$, 
thus, $ \exp \{ \overline{\Phi}_{1,s; \bfl} \} $ 
is the unique positive number such that
\begin{equation}
 \label{timetheta}
     \sum_{i \in \cA} \, \frac{ \overb_{i,s}(\bfl)}{ 
         \exp \{\overline{\Phi}_{1,s; \bfl} \} + 
                 \overb_{i,s}(\bfl)} = 1,
\end{equation}
where the  $\overb_{i,s}$ are the time-dependent 
equivalents of the functions $\overb_i,$ 
encountered in \mbox{Theorem \ref{weaksurvival}} 
(for a fuller description, see 
(\ref{htimematrices}) and (\ref{bigphimeans})).
In \cite{lal2}, the function $V(s) =
\overline{\Phi}_{1,s} - \log (2d-1)$ was called
the {\em growth profile} of the contact process.
  
\begin{theorem}
   \label{tgrowth}
Let $d>1 $ and $\bfl$ in $ \overR_2.$  
Then there exist smallest and largest
solutions $ 0 < s_1 \leq s_2 < \infty$ of 
 $\overline{\Phi}_{1,s; \bfl} = 0.$ 
Almost surely on the event of survival,
\begin{eqnarray}
     \label{tsmalldist}
         \lim_{t \rightarrow \infty} r_t /t = 1/s_2, \\*[0.15cm]
       \label{tlargedist}      
           \lim_{t \rightarrow \infty} R_t /t = 1/s_1.
\end{eqnarray}
Moreover, for each $s > 0$ such that $\overline{\Phi}_{1,s; \bfl} > 0,$
almost surely on the event of survival, 
\begin{equation}
   \label{tnumberinfected}
     \lim_{n \rightarrow \infty} \frac{1}{n} \log N_n(ns) = 
     \overline{\Phi}_{1,s; \bfl} .
\end{equation}      
\end{theorem}

See Section 10 for the proof.

\bigskip
{\bf Remarks.}
\begin{enumerate}
\item[(1)]
\vspace*{-0.2cm}
The phenomenon that the lead eigenvalue of $M_2$  converges to $1$
at the phase transition will be shown to be responsible
for both convergences of the Hausdorff dimension in the lower
bounds to $ \frac{1}{2}$ of some Hausdorff dimension 
in (\ref{halfdim1}) and (\ref{halfdim2}),
 and, offers a useful {\em characterization}
of the phase transition. 
\item[(2)]
\vspace*{-0.2cm}
Analogous results can be proved for a contact process
with infection rates between sites within {\em finite distance} 
(but not exclusively nearest neighbours).
\item[(3)]
\vspace*{-0.2cm}
It is worthwhile noticing that the formulae for the Hausdorff dimensions
resemble the ones for the branching random walks on $\cT_{2d},$ with
the average asymptotic infection probabilities $\overb_i$ corresponding
to the generating functions $F_i.$
Surprisingly, for the isotropic models, the dimensions of the limit sets
coincide at the phase separations. For each anisotropic contact process
at the phase transition, there is an anisotropic branching random
walk with the same Hausdorff dimension of the limit set.
\item[(4)]
\vspace*{-0.2cm}
How to define a suitable symmetry assumption for the contact
process on a homogeneous trees of {\bf odd} degree is a bit 
less obvious. Suppose that, at any vertex, 
associated with each of the $d$ emanating edges
that point away from the root vertex be an infection rate
$\lambda_i,$ $i \in \cAp,$ and associated with the emanating
edge that points back towards the root be an infection rate
$\lambda_{k^{-1}}$ for $k \in \cAp.$ Thus, two infection rates
$\lambda_k$ and $\lambda_{k^{-1}}$ are attached to an edge, one
being in use when the infection ``crosses the edge forwards''
and the other being in use when the infection ``crosses the edge
backwards''. Suppose that the contact process on the tree $\cT_{d+1}$ of 
degree $d+1 > 2$ satisfy the {\em symmetry} assumption
$\lambda_k =\lambda_{k^{-1}}$ for each $k \in \cAp.$ \\*[0.1cm]
{\bf Question 1:}
Which of our results continue to hold for this model ?

We remark that this setup gives rise to nonhomogeneous
vertices, more precisely, the set of infection rates
is not the same at different vertices.
The method developed in this paper relies on the geometric
decay (see Lemma \ref{periodicmaximal}), which in turn invokes
the homogeneity of the vertices.

Preliminary results that we pursue in subsequent papers
indicate that if $\lambda_{k^{-1}} = \lambda_*$ for all $k \in \cAp$
({\bf isotropic backtracking}),
then weak survival exists and, for $\bfl \in \cR_2,$ 
all of Theorems \ref{limitset} and \ref{halfboundary}
hold, where $\theta(\bfl) = \sum_{i \in \cAp} \,  \overb_i(\bfl),$
with the possible exceptions of the second parts in  
(\ref{halfdim1}) and (\ref{halfdim2}). Moreover, the 
discontinuity set consists of all $\bfl$ which satisfy
the equation 
$ \sum_{i \in \cAp} \, \overb_i(\bfl)^2 \, \lambda_i /\lambda_* = 1.$
\\*[0.1cm]
{\bf Question 2:}
Does the contact process on $\cT_{d+1}$ for $d+1 > 2$ 
with $\lambda_* \leq \max_{j \in \cAp} \lambda_j$
always exhibit weak survival, more specifically, is the
transition from $\cR_1$ to $\cR_3$ precluded ?
\end{enumerate}

The rest of the paper is organized as follows.
Section 2 discusses the percolation structure, the Markov property,
and the strict monotonicity of the $\beta_x,$ 
introduces the embedded labelled Galton-Watson trees, and
outlines the strategy of proof (Section 2.5) for the existence
of weak survival. In Section 3, we give some background on the
Thermodynamic Formalism, describe the potential functions and
Gibbs states, and prove when $\eta <1.$
Section 4 defines the first-passage and backscatter matrices and
analyzes the strict monotonicity of their lead eigenvalues.
Section 5 is concerned with the continuity features of the lead
eigenvalues, derives a characterizing
equation at criticality,
establishes weak survival, and includes the proofs of  
Theorems  \ref{etaatboundary}, \ref{properties} and 
Theorems \ref{weaksurvival} through  
\ref{extinctionattransition} on the critical contact process.
Section 5 also includes a result on the distribution of the
limit points of the contact process (Theorem \ref{statprocess}). 
The upper and lower bounds for the Hausdorff dimensions are
taken care of in Sections 6 and 7. Section 8 verifies inequalities
(\ref{halfdim1}) and (\ref{halfdim2}) as an appeal to the Gibbs
Variational Principle, which will complete the proofs 
of Theorems \ref{limitset} and \ref{halfboundary}. 
Section 9 reviews the isotropic case and Section 10 addresses 
the dispersal of the contact process over space-time.


\section{Basic Properties and Embedded Galton-Watson Trees}  
\setcounter{equation}{0}

First we describe the strong Markov property and the visualization of the
contact process as a percolation structure, upon which the verification 
of the {\em strict} monotonicity of the functions $u_x$
along certain directions in the space of infection parameters is based. 
This feature together with the continuity properties 
will be essential to our proof of \mbox{Theorem \ref{weaksurvival}}.
Furthermore, we will describe the Galton-Watson trees embedded 
in the set of vertices ever to be infected.

\subsection{Percolation Structure and Strong Markov Property} 

{\bf Percolation structure.}
The contact process may be constructed via the usual 
{\em percolation structure} on $\cT \times (0, \infty),$
that is, as a system of independent Poisson processes 
attached to vertices and ordered pairs of neighbouring
vertices. For each vertex $ x \in \cT,$ the Poisson process
attached to $x$ has rate $1,$ and determines the recovery
times, specifically, at every occurrence time, site 
$x$ recovers if it is infected. For each ordered pair 
$(x,xi),$ $ i \in \cA,$ of
neighbouring vertices, the Poisson process 
attached to $(x,xi)$ has rate $\lambda_i,$ the occurrence
times being precisely those times when an infection at
$x$ may jump to $xi.$ Occurrences in these Poisson processes
are marked on a system of directed rays $\{ x \} \times
[0,\infty)$ connected to the vertices $x$ of $\cT,$ so that
(A) at each occurrence time $t$ of the Poisson process attached
to $(x,x i)$ an {\em infection arrow} is drawn from $(x,t)$
to $(x i,t),$ and (B) at each occurrence time $t$ of 
the Poisson process attached to $x$ a {\em recovery mark} *
is attached to $(x,t).$ There are no simultaneous occurrences
of infection arrows and/or recovery marks in the percolation
structure. At time $t,$ the contact process now consists of
all those vertices $y$ for which there is a (directed) path 
through the percolation structure, the system of rays and
arrows just described, that begins at the root vertex $1,$
ends at $(y,t),$ and does not pass through any recovery marks *.  
An {\em infection trail} is a connected path in the percolation
structure that does not pass through any recovery marks.

\smallskip
{\bf Strong Markov Property.}
An important property that the contact process enjoys 
is the strong Markov property (see also \cite{lase2}, Section 2.2).
Let $G, F_1, F_2, \ldots , F_k$ be pairwise 
nonoverlapping parts of the tree $\cT,$ and let $S_1, S_2, 
\ldots, S_k$ be stopping times determined by the percolation
structure over $G.$ Then conditional on the percolation structure
over $G,$ the post-$S_i$ portions of the percolation structures
over the sets $F_i$ are independent, and for each $i,$ the
post-$S_i$ percolation structure over $F_i$ has the same 
distribution as the entire percolation structure over $F_i.$


\subsection{Downward Infection Trails}

Fix a generator $a \in \cA,$ and consider the subtree
$\cT^* = \cT - \cT(a^{-1})$ of $\cT,$ each of which
vertices except the root vertex $1$ is represented by a word
$x \in \cG$ beginning with a letter $x_1 \not = a^{-1}.$ 
We arrange the tree $\cT^*$ in {\em levels} $\cL_0, \cL_1,
\cL_2, \ldots,$ where $x \in \cL_n$ if and only if $\vert x
\vert =n.$ There are $(2d-1)^n$ vertices at the $n$th level
$\cL_n.$ Moreover, for every $n \geq 1,$ define $\cL_n^*$
to be the subset of $\cL_n$ containing those vertices whose
word representation $x=x_1 x_2 \ldots x_n$ terminates in the
letter $x_n = a.$ Observe that for any vertex $x \in \cL_n^*,$
the set of nearest neighbours of $x$ in $\cL_{n+1}$ is
$\{ xy: \, y \in \cL_1 \}, $ thus, for $x \in \cL_n^*,$
the tree $\cT(x)$ is the left translate of the tree $\cT^*$
by the group element $x.$

Let $x$ be a vertex in level $\cL_n$ and let $y$ be
a vertex contained in $\cT(x).$ In other words, the word
$x$ is a prefix of the word $y$ and $y$ must lie in a level
$\cL_{n+m}$ at larger distance than $x.$ Define a {\em downward
infection trail} from $x$ to $y$ to be an infection trail
that begins at $x$, does not exit $\cT(x),$ and first reaches
$\cL_{n+m}$ at $y,$ where it terminates.
For every vertex $x \in \cG,$
define
$$
   \cD_x = \{ \exists \mbox{ downward infection trail }
                 \{ \mbox{root}\} \rightarrow x \mbox{ beginning
               at } t=0 \}.
$$
Thus, if we write 
\begin{equation}
   \label{trailprobab}
   w_{x} = P \{ \cD_x \},
\end{equation}
then as an appeal to the monotonicity and the strong Markov property
of the contact process, we have 
 $w_{xy} \geq w_x w_y $ for every $x,y \in \cG$ so that
$ \vert x y \vert = \vert x \vert + \vert y \vert.$ It is 
obvious that $w_x \leq u_x.$
For each {\em periodic} sequence $y_n = x x \ldots x  \in \cG_{nk}$
with $ x \in \cG_k,$ the limit 
       $ \lim_{n \rightarrow \infty}(w_{y_n}(\bfl))^{1/n}$ exists
for every $\bfl$ because of the homogeneity of the process 
at each vertex, and in fact,
by the same arguments as presented
for the isotropic contact process in \cite{lase2}, 
\begin{equation}
  \label{wlimitbeta}
   \lim_{n \rightarrow \infty} (w_{y_n}(\bfl))^{1/n} = \beta_x(\bfl),
\end{equation}
where $\beta_x(\bfl)$ was defined in 
(\ref{uperiodicsubadd}). Moreover, 
$ w_{y_n}(\bfl) \leq \beta_x(\bfl)^n$ 
for every $n \geq 0.$ The non-periodic analogue of 
(\ref{wlimitbeta}) in the anisotropic case
will be derived in Section 5.2. 


\subsection{Labelled Galton-Watson Processes and Trees}

Every Galton-Watson chain has its associated genealogical tree 
$\tau.$ This {\em Galton-Watson tree} may be described as
follows (see also \cite{hawk}). Vertices of $\tau$ are arranged
in levels $V_0, V_1, V_2, \ldots.$ The vertices of level $V_n$
represent the individuals of the $n$th generation of the 
corresponding Galton-Watson chain.  Edges of the tree connect
children and their parent, thus, there are edges only between 
vertices of successive levels.
The limit set $\Lambda_{GW}$
of the tree is the set of {\em ends}, i.e.\ the set of all
infinite paths that start at the root $V_0$ and visit each
of the levels $V_n$ exactly once. If the Galton-Watson chain
dies out, then the limit set is empty. 
For each $\alpha \in (0,1),$
the $d_{\alpha}$-distance between two ends $\gamma, \gamma'$
is defined to be $\alpha^n,$ where $n=n(\gamma, \gamma')$
is the last level $V_n$ where the paths touch. The set
$\Lambda_{GW}$ is a compact metric space with this metric 
$d_{\alpha}.$ 

For any {\em finite} set $\cB,$ a {\em labelled Galton-Watson
process} with label space $\cB$ is determined by a probability
distribution $Q$ on the set of $2^{\cB}$ 
subsets of $\cB.$ Each individual $\zeta,$ regardless of its type,
produces a random set $\cO_{\zeta}$ of offspring, with distribution
$Q,$ and the offspring sets of different individuals are 
conditionally independent, as in an unlabelled Galton-Watson 
process. Thus, a labelled Galton-Watson process is a multitype
Galton-Watson process in which (i) all types have the
same offspring distribution, and (ii) the offspring distribution 
is constrained to allow at most one individual of each label.
If $Z_n$ denotes the cardinality of the $n$th generation, then
the sequence $\{ Z_n \}_{n \geq 0}$ is an ordinary Galton-Watson
process. Moreover, if the expected cardinality of a random set 
chosen according to the distribution $Q$ is bigger than $1,$
then $\{Z_n \}_{n \geq 0}$ is {\em supercritical}.

A {\em labelled Galton-Watson tree} $\tau$ is the genealogical
tree associated with a labelled Galton-Watson process so that
the labels of the corresponding individuals are assigned to the
vertices of $\tau.$ If the underlying Galton-Watson process 
is supercritical, then with positive probability the tree $\tau$
is infinite. On this event, each end of $\tau$ will be naturally
identified with a unique semi-infinite sequence $\omega =x_1 x_2
\ldots, $ where the end crosses the $n$th level through 
$x_n $ (The root need not be labelled).
Hence, the set $\partial \tau$ 
of ends of $\tau$ is naturally embedded
in the sequence space $\cB^{\bf N}.$


\subsection{Embedded Labelled Galton-Watson Trees $\tau_r$}

The construction of embedded Galton-Watson trees $\tau_r$
in the anisotropic contact process resembles the ones 
for anisotropic branching random walks \cite{lahu} 
and isotropic contact processes \cite{lase2}.
The {\em offspring} of a vertex $x \in \tau_r$ will be
vertices $y$ at distance $r$ from $x$ such that there 
is a  downward infection trail  beginning at $x$ 
and ending at $y. $ 
It is apparent that the concatenation of infinitely
many such downward infection trails gives rise to an
end of $\tau_r$ that is contained in the limit set $\Lambda.$

\smallskip
{\bf Definition of $\tau_r.$} Fix an integer $r \geq 1.$ Define
generations $V_n(r)=V_n \subset \cL_{nr}^*$ inductively as follows:
\begin{enumerate}
\item[(i)] $ V_0 = \{1\} .$
\item[(ii)] For each $x \in V_n,$ the {\em offspring} of $x$ 
            are those $y \in \cL_{nr+r}^* $ such that there is
           a downward infection trail $\xi_x$ from $x$ to $y,$
          beginning at time $\cS_n,$ the first time when $x$ is reached,
           and first reaching $y$ at time
          $\cS_{n+1}$ (suppose that $\cS_0=0$), and this infection
          trail $\xi_x$ is the {\em first} downward infection
          trail beginning at $(\cS_n, x)$ to reach $y.$  The offspring
          vertex $y$ has {\em parent} $x.$
\item[(iii)] $V_{n+1}$ is the set of all offspring of vertices in 
             $V_n.$
\end{enumerate}

Note that the random times $\cS_i$ are stopping times.
For each vertex $y \not = 1$ of the tree $\tau_r,$
define the {\em label} of $y$ to be the word $x^{-1} y ,$
where $x$ is the parent of $y.$ Each label is a reduced word
of length $r$ which ends in the letter $a.$ Thus, the set of
labels is finite. 

\begin{lemma}
  $\tau_r $ is a labelled Galton-Watson tree.
\end{lemma}

\begin{proof}
The infection trails $\xi_x$ associated with different
$x \in V_n$ are distinct and the existence of an infection
trail $\xi_x$ depends only on the percolation structure
in $\cT(x)$ after time $\cS_n,$ thus, is independent of the 
pre-$\cL_{nr}$ history of the contact process initiated
with $\xi_{\{1\}}. $ Hence, 
the trails $\xi_x, \,  x \in V_n,$
involve nonoverlapping regions of the percolation structure
and do not overlap the region of the percolation structure
that determines $V_n,$ whence by the strong Markov property,
are mutually independent. Since each vertex of $V_n$ has
final letter $a,$ the contact processes initiated by the trails
$\xi_x, \, x \in V_n,$ are all ``oriented" the same way,
relative to the level structure of the tree $\cT^*.$
This implies that the offspring distributions for vertices 
$x \in V_n$ as probability distributions on the set of labels
are all the same as that one of the initial vertex $1 $ 
in $V_0.$ As a consequence, $\tau_r$ is a labelled Galton-Watson
tree.
\end{proof}

Observe that the relation between the metrics 
$d_{\alpha}^{\tau_r}$ for the tree 
$\tau_r$ and $d_{\alpha}$ for the full tree $\cT$ is 
\begin{equation}
  \label{gwmetric}
     d_{\alpha}^{\tau_r}(x,y)  = d_{\alpha^r}(x,y).
\end{equation}
The mean offspring number $\mu_r$ for the Galton-Watson trees 
$\tau_r$ is, by construction, 
$$
    \mu_r = \sum_{x \in \cL_r^*} w_x.
$$


\subsection{Geometric Decay and Strict Monotonicity of $u_x$}

It is not difficult to show that,
in the isotropic case, $\bfl \in \cR_3$ implies $\beta(\bfl) =1.$
The reverse direction of this statement is more subtle but
has been shown in \cite{lase2}, more precisely, if $\bfl \not \in
\cR_3,$ then $\beta(\bfl) \leq 1/ \sqrt{2d-1},$ with strict inequality
holding for each $\bfl$ in the interior of the weak survival phase and
equality emerging at the transition. 
The feature $\beta < 1$ guarantees exponential decay of
the infection probabilities $u_n$ in $n,$ 
for every integer $n \geq 1.$
It is essential to verifying both, the strict monotonicity and
continuity of the relevant growth variables of the contact process
in the infection parameter. 

In fact, for the anisotropic contact process as well,
$\beta_i(\bfl) < 1$ for each $i \in \cA$ and 
for $\bfl \not \in \cR_3$ would be a desirable property
in proving strict monotonicity and continuity of certain
functions in the infection parameters. Different values
of these functions at $\overR_1 \cap \overline{\cR_1^c}$ and 
$\overline{\cR_3^c} \cap \overR_3$ could then be used to argue
that the region $\cR_2$ has nonempty interior. 
To help the understanding of the rest of the paper,
we here outline our {\em strategy} in establishing the
existence of a weak survival region.
We shall prove the following: 

\smallskip
\mbox{} \\
{\bf Strategy of proof.}
\begin{itemize}
\item[(A)]
\vspace*{-0.2cm}
If each $\beta_i(\bfl) < 1,$ then we have 
continuity and strict monotonicity of certain functions,
\item[(B)]
\vspace*{-0.2cm}
a set $\cD_c$ of discontinuities is located,
\item[(C)] 
\vspace*{-0.2cm}
as $\bfl$ with each $\beta_i(\bfl) < 1$
approaches any point in $\cD_c,$
each $\beta_i$ stays {\em bounded away} \mbox{from $1.$}
\end{itemize}

From (A)---(C) together with several other
considerations, it will then follow
that in fact the set $\cD_c$
coincides with the boundary 
$\overline{\cR_3^c} \cap \overR_3$ 
and that $(\overline{\cR_1^c} \cap \overR_1) \cap 
(\overline{\cR_3^c} \cap \overR_3) = \emptyset.$
It will follow that $\cD_c$ separates
the regions $\cR_1 \cup \cR_2$ and $\cR_3,$ and that
the functions at hand are strictly increasing
and continuous for $\bfl \in \mbox{int} (\cR_1 \cup \cR_2)$  
and left-continuous for $\bfl \in \overR_1 \cup \overR_2$ 
such that each $\lambda_k >0$ for $d>1.$
Sections 5.5 and 5.6 are the principal sections in our
discussion of identifying the three regions
and classifying their boundaries, whereas 
the sections beforehand introduce machinery and
give a full account on the specifics about the functions
which are to be continuous and strictly monotone. Some
of the latter parts may appear a bit technical.

The distinguished set of all $\bfl$ that have each 
$\beta_i(\bfl) <1$ is given a name as follows.
If $\cQ = \{ \bfl \in {\bf R}^d: \mbox{ each } \lambda_k \geq 0 \},$ 
define the sets
\begin{eqnarray*}
   \label{kset}
    \cK & = & \{ \bfl \in \cQ : \, \beta_i(\bfl) < 1 \mbox{ for 
             each } i \in \cA \} ,  \\
      \label{poskset} 
    \cK_0 & = & \{ \bfl \in \cK: \, \lambda_i > 0 \mbox{ for each }
                 i \in \cA \}.
\end{eqnarray*}

{\bf Geometric decay of $u_x.$} We begin with some estimates
on the decay of the infection probabilities. 
We point out that the upper bound for $u_x(\bfl)$ in 
(\ref{ugeometric}) below is uniform in $\bfl$ for $d \geq 2,$
that is, if we assume that there is
some $\epsilon > 0$ such that at least two infection rates
$\lambda_i, \lambda_j \geq  \epsilon $ for two distinct 
$ i, j \in \cAp .$ If $\bfl$ is restricted to those $\bfl \in \cK,$ 
then it can be shown that there is some $\rho = \rho(\epsilon) < 1$ 
such that each $\beta_i(\bfl) \leq \rho < 1.$  

\begin{lemma}
 \label{periodicmaximal}
For every fixed $\bfl,$ each integer $n\geq 0,$ 
and every $x \in \cG_n,$
\begin{equation}
     u_x(\bfl) \leq [ \max_{a \in \cA} \beta_a(\bfl)]^n.
\end{equation}
For each $\bfl \in \cK$ and each $x \in \cG_n,$ 
there is some constant $0 < \gamma < 1$ such that  
\begin{equation}
      \label{ugeometric}
     u_x(\bfl) \leq \gamma^n.
\end{equation}
\end{lemma}

\begin{proof}
The proof proceeds by induction over the distance of the vertices from the
root. Fix $\bfl.$
Without loss of generality, we may assume that
$\lambda_b \geq \lambda_i$ for every $i \in \cA.$ Let
$y_n = b b \ldots b$ with $\vert y_n \vert = n.$ For $n=1,$
clearly, $u_{y_1}(\bfl) = u_b(\bfl) \geq u_i(\bfl)$ for every
$i \in \cA$ by the
definition of the functions  $u_x(\bfl) $ and the one of the rules of
infection of the contact process.
Now assume that $u_{y_k}(\bfl) \geq u_x(\bfl)$ for every $x \in \cG_k$
and every $k=1,2, \ldots, n-1.$ 
Thus, no $x \in \cG_{n-1}$ has larger probability than $y_{n-1}$
to ever be infected. But upon infection of $y_{n-1},$
no vertex in $\cG_n$ has larger probability to
be infected than $y_n$ because of the rules of infection of 
the contact process, the homogeneity of the process at each vertex,
and the fact that, along any path $bbb....,$ 
the neighbourhoods of the vertices look the same.
Hence, $u_{y_n}(\bfl) \geq u_x(\bfl)$
for every $x \in \cG_n.$ Since this argument is valid for every integer
$n > 0,$ the desired results follow from the fact that
$u_{y_n}(\bfl) \leq \beta_b(\bfl)^n$
for every $n \geq 0.$ 
\end{proof}

\medskip
{\bf Strict Monotonicity Properties.}
For each infection parameter $\bfl \in \mbox{int}(\cK),$
there exist directions of increase and decrease, 
respectively, where the infection probabilities $u_x$ {\em strictly}
increase or decrease, respectively. 
It is not apparent, however, whether there
is an easy criterion to decide for each vector pointing away from
some $\bfl,$ whether the
infection probabilities have strict monotonicity properties.
We will construct directions of strict monotonicity for each 
$\bfl$ by ``thinning" the percolation structure of the
contact process. This useful idea to modify the percolation structure
originates from \cite{lal2}, where in the isotropic case
strict monotonicity of the function $\beta(\cdot)$
in the infection parameter was shown. In the anisotropic case, 
the precise statement is more subtle, though. 
For this purpose, we first
describe a {\em modified contact process}
that arises when the percolation
structure is manipulated by a set of Bernoulli-$p$ random variables,
where $p$ is precisely chosen to decrease the {\em smallest}
infection rate $\lambda_c,$ say, to $\tilde{\lambda}_c < \lambda_c.$ 
The next proof will show that tuning 
the smallest infection rate, forces all infection rates to 
decrease, which might possibly be more than necessary but it will
suffice to lead to strict monotonicity of many characteristics of
the contact process that will be discussed later.

Suppose that $\lambda_c = \min_{a \in \cA} \lambda_a$ is the
smallest positive infection rate of the contact process $A_t$ with
infection parameter $\bfl$ and recovery rate $1.$
Fix $\tilde{\lambda}_c < \lambda_c.$
To each infection arrow $\omega$ of the percolation structure, 
there is attached a \mbox{Bernoulli-$p$} random variable 
$\xi_{\omega}.$ These are conditionally independent, given
the realization of the percolation structure.
Choose the value 
$p= P \{ \xi_{\omega}=1 \}$ so that 
$$
    \tilde{\lambda}_c =  p \lambda_c/ ( 1+ q \lambda_c).
$$
where  $q = 1-p \in (0,1).$ Call this assignment of a set of
Bernoulli-$p$ random variables to the percolation structure
a $p$-{\em thinning} of the percolation structure relative 
to $\lambda_c.$ A version $A_t'$ of the contact process with
infection rates  $p \lambda_i $
and with recovery rates $r_i= 1+ q \lambda_i,$ 
$i \in \cA,$
may be constructed using the augmented percolation structure by 
(1) first setting up a modified percolation structure by
changing every infection arrow $\omega$ such that $\xi_{\omega}
=0$ to a recovery mark *, then (2) defining $A_t'$ to be
the set of all vertices $y$ for which there exists
a directed path from $(1,0)$ to $(y,t)$ in the above manipulated
percolation structure that does not pass through any recovery
marks *. The new contact process $A_t'$ is a time-changed
version of a contact process with infection rates 
$\bfl' = ( \{ \lambda_a' \}_{a \in \cAp})$ 
and recovery rates $r'=( \{r_a' \}_{a \in \cAp})$ 
with $\lambda_i' = p \lambda_i/(1+ q \lambda_c) < \lambda_i $
and  $r_i' = (1+ q \lambda_i)/(1 + q \lambda_c).$  
Thus, since
$r_c' =1$ and $r_i' \geq 1$ for each $i \in \cA,$
there is a contact process $A_t^*$
that has the same percolation structure
as $A_t',$ whose recovery rate is $1$ 
and whose infection rates are no larger than $\lambda_i'.$   

Note that it is enough to state the next result for a
contact process that has only positive infection rates $\lambda_a$
(because otherwise the set of generators $\cAp$ may be updated).

\begin{proposition}
  \label{monotone} 
Let $A_t$ be a contact process with infection parameter 
$\bfl $ and recovery \mbox{rate $1.$} Let $\bfl \in \cK_0$ with 
$0 < \lambda_c = \min_{a \in \cA} \lambda_a.$
Assume that its percolation structure is
$p$-{\em thinned} relative to $\lambda_c$ for some $ 0 < p <1$
so that a contact process $A_t^* $ is obtained with infection
parameter $\tilde{\bfl}$ and recovery rate $1$ with
$ 0 < \tilde{\lambda}_c  = p \lambda_c/(1+ q \lambda_c) < \lambda_c $
and $\tilde{\lambda}_i \leq p \lambda_i/(1+ q \lambda_c) < \lambda_i$  
for every $i \in \cA.$ Then there is some constant $0 < \omega < 1$
such that for every $n$ and $ x \in \cG_n,$
\begin{equation}
   \label{udecrease}
          u_x (\tilde{\bfl}) \leq u_x(\bfl) \omega^n .
\end{equation} 
In particular,
    $ \beta_x(\tilde{\bfl}) < \beta_x(\bfl)$
for every $x \in \cG,$ where $\beta_x$ is defined
in (\ref{uperiodicsubadd}). 
\end{proposition}

\begin{proof}
The ideas of the proof are essentially the same as those 
given in \cite{lal2}, Proposition 9, for the isotropic case.
We outline the necessary modifications and refrain
from repeating the entire (not short) proof. The argument rests
on the previously described modification of the percolation
structure of the contact process by $p$-thinning by means of 
Bernoulli-$p$ random variables relative to $\lambda_c.$
Choose the unique value $p$ so that 
$ \tilde{\lambda}_c =  p \lambda_c/ ( 1+ q \lambda_c).$
Assume that $A_t$ is a contact process with infection
parameter $\bfl,$ and recovery rate $1$ and that $A_t^* $
is a contact process with infection
parameter $\tilde{\bfl}$ and recovery rate $1$ with
$ \tilde{\lambda}_c  = \lambda_c' =
   p \lambda_c/(1+ q \lambda_c) < \lambda_c $
and $\tilde{\lambda}_i \leq  \lambda_i' = 
p \lambda_i/(1+ q \lambda_c) < \lambda_i$  
for every $i \in \cA.$ Thus, all infection rates have been decreased
by the thinning process.

To verify (\ref{udecrease}),
it suffices to show that there is some $\omega \in (0,1)$ such that
the contact process $A_t',$ gotten from the $p$-thinned percolation
structure, with infection
parameter $\bfl'$ and recovery rates $r'$ has 
$u_x(\bfl') \leq u_x(\bfl) \omega^n$ for every $x \in \cG_n$ and every
integer $n >0,$ because the contact process $A_t^*$ is
equivalent to $A_t'$ in the sense that
they have the same percolation structure, and thus,
$u_x(\tilde{\bfl}) = u_x(\bfl')$ for every $x \in \cG.$ 

Let $G_x$ be the event that the contact process $A_t$ infects
vertex $x \in \cG_n$ at some finite time,
and let $G_x'$ be the corresponding event for the contact process
$A_t',$ thus, $P \{ G_x \} = u_{x}(\bfl) $ and 
$ P \{ G_x' \} = u_{x}(\bfl')$
(Here we abuse the notation since except for this proof, $u_x$ refers
to a contact process with recovery rate $1$). 
By construction, $G_x' \subset G_x$ because every infection arrow
in the modified percolation structure occurs in the unmodified one as
well, and on the other hand, every recovery mark * in the unmodified
percolation structure is retained in the modified percolation
structure. On the event $G_x,$ there is at least one directed path through
the unmodified percolation structure that leads from $(1,0)$ to
$ x \times (0, \infty),$ in fact, there may be many overlapping
such paths. Call an infection arrow $\omega$ in the unmodified
percolation structure {\em essential} for the event $G_x$ if (1) event
$G_x$ occurs, and (2)  changing $\omega$ from an infection
arrow to a recovery mark would destroy all directed paths from
$(1,0)$ to $ x \times (0, \infty).$ Define $N_x$ to be the number
of essential arrows for the event $G_x,$ when $G_x$ happens,
and $N_x=0$ when $G_x$ does not happen. Since removing any one of
the essential arrows would disconnect $(1,0)$ from
$ x \times (0, \infty),$ in order that event $G_x' $ occur it 
is necessary that $N_x \geq 1$ and that $\xi_{\omega} = 1 $ for
{\em every} essential arrow $\omega.$ This event has conditional
probability $p^{N_x},$ given a realization of the unmodified
percolation structure. Hence,
$$
  u_{x}(\bfl') = P \{ G_n' \} = E p^{N_x} I_{G_x},
$$   
where $I_{\{ \cdot \}}$ denotes the indicator function.
If it were the case that $N_x \geq cn$ on $G_x$ for some positive
constant $c,$ then it would follow that
 $u_{x}(\bfl') \leq u_{x}(\bfl) p^{cn}.$
On the other hand, if for some $\rho < 1 $ it were
the case that $P \{ N_x \leq cn \vert G_x\} \leq \rho^n,$ then it
would follow that $u_{x}(\bfl') \leq u_{x}(\bfl) (p^{cn} + \rho^n).$ 
Therefore in both cases, it would follow that 
there is some $0 < \omega <1$ such that
 $u_{x}(\bfl') \leq u_{x}(\bfl) \omega^n$ and,
 since this argument holds
 for every $n>0$ and $x \in \cG_n,$ prove (\ref{udecrease})
 and, by subadditivity, that  
        $ \beta_x(\tilde{\bfl}) < \beta_x(\bfl).$
Now, verifying Lemma 3 in \cite{lal2} (there stated in the isotropic
case) for the anisotropic contact process will complete the proof:

There exist constants $ c<  \infty $ and $0 < \rho < 1$ such that for
all sufficiently large $n,$
$$ P \{ N_x \leq cn \vert G_x \} \leq \rho^n$$ (\cite{lal2}, Lemma 3).
Note that this statement refers only to the unmodified percolation
structure. In the proof of Lemma 3 given in \cite{lal2},
only features common to both the anisotropic
and isotropic contact processes are relied on except for 
a single passage in the very last paragraph of the paper, which 
requires a concrete estimate of the underlying anisotropic 
contact process, that is, an upper bound $\gamma^n,$ $0 < \gamma <1,$
for the probability $u_x(\bfl)$ for every sufficiently large $n$
and each $x \in \cG_n.$ But, since $\bfl \in \cK_0 \subset \cK,$
this result follows from (\ref{ugeometric}) 
in Lemma \ref{periodicmaximal}. This finishes our proof.
\end{proof}

We say that the directions given by the 
vectors $\bfl -\tilde{\bfl}$ and $\tilde{\bfl} - \bfl$ are
{\em directions of increase} and {\em decrease} for $\bfl,$
respectively, for the contact process with infection 
parameter $\bfl$ (and recovery rate 1), where $\tilde{\bfl}$
is defined in Proposition \ref{monotone} for any $0 < p < 1.$
Note that any function in $\bfl$ that is strictly increasing 
along directions of increase for $\bfl$ necessarily is
strictly increasing if every component of $\bfl$ is increased.


\section{Potential Functions and $\eta < 1 $}
\setcounter{equation}{0}

The proofs of the principal results of this paper rely on classical
results from the theory on Gibbs states as developed in 
\cite{bowe,ruel}. They enable us to ``normalize" the ensemble 
of infection probabilities $u_x$ so that, for any $\varepsilon>0$
and at all fixed sufficiently large distances
$n$ from the root vertex, there exists a  
shift-invariant probability measure
concentrated on the set of vertices in $\cG_k$ with 
$ n(1- \varepsilon) < k < n (1 +\varepsilon)$ 
(see e.g.\ the remarks surrounding  (\ref{azeta}) through 
(\ref{genlength})).

We begin with describing this approach via 
potential functions, quote a result that applies the theory to
some matrices that satisfy a H\"{o}lder condition, and proceed
to apply these tools to the collection of
infection probabilities for the contact process.
Our goal is to define Gibbs states supported by  the geometric
boundary $\Omega,$ which will allow us to keep track of the growth of the 
number of infection trails that wander off to infinity and of 
those that return to the root vertex at some finite time. 
From the properties of Gibbs states, we will derive 
estimates that describe the dispersal of the
contact process on $\cT$ in space-time, as for instance,
the rightmost and leftmost positions 
of the infection at time $t,$ which both will turn out, 
almost surely on the event of survival, to move at a distance 
from the root that is linear in time. 


\subsection{Background: Thermodynamic Formalism}

Recall that $\cA$ denotes the set of generators of $\cG$ and their inverses.
Let $\cAN = \{ \mbox{one-sided }$ 
   \linebreak 
$\mbox{infinite sequences from } \cA \} $ and
   $\cAZ = \{ \mbox{two-sided infinite sequences from } \cA \}. $ Recall
$\Sigma$ to be the set of all doubly infinite reduced words 
 $\xi = (x_n)_{n=- \infty}^{\infty}$ from $\cA$ (reduced means that, for
 every $n,$ $x_{n+1} \not = x_n^{-1}$). The subset $\Sigma$ is a closed
 subset of the coding space $\cAZ,$ thus, is compact in the metric 
 $d_{\alpha}.$
 Let $\sigma: \Sigma \rightarrow \Sigma$ denote the forward shift on $\Sigma,$
 that is, $\sigma(x_k x_{k+1} \ldots) = x_{k+1} x_{k+2} \ldots $ for every $k.$
 Note that $\sigma$ is Lipschitz continuous. For any function
 $f: \Sigma \rightarrow {\bf R},$ define
 $$
     S_n f = f + f \circ \sigma + f \circ \sigma^2 + \ldots + f \circ 
\sigma^{n-1}
 $$
 for every $n \geq 1.$

\bigskip
{\bf Potential Functions, Gibbs States and Pressure.}
Define the $n$-cylinder 
sets
$$
   \Gamma(i) = \Gamma(i_1 i_2 \ldots i_n) =
   \{ \xi = (x_k)_{k=- \infty}^{\infty} \in \Sigma: x_j =i_j, \; 1 \leq j
 \leq n\}
$$  
for every $n$ and $i = i_1 i_2 \ldots i_n \in \cG_n.$ 
In spirit of \cite{bowe}, for any H\"{o}lder continuous function $f$ on 
$\Sigma,$
there is a real constant $P(f),$ some constants $0 < C_1 \leq C_2 < \infty,$
and a unique $\sigma$-invariant probability measure
$\mu_f$ on the Borel sets of $\Sigma$ such that for each
$i = i_1 i_2 \ldots i_n \in \cG_n,$ 
\begin{equation}
    \label{gibbsf}
  C_1 \leq \frac{\mu_{f}(\Gamma(i))}{\exp \{S_n f(j) - n P(f) \}} \leq C_2 
\end{equation}
for every $j \in \Gamma(i).$ The measure $\mu_f$ is called the 
{\em Gibbs state}
with {\em potential function} $f,$ and the constant $P(f)$ is called the 
{\em pressure} of $f.$ 

We list some features of the pressure function that will be useful 
to our subsequent analysis.
Two H\"{o}lder continuous functions $f$ and $g$ are called {\em cohomologous}
if there exists a H\"{o}lder continuous function $h$ such that 
$f-g = h - h \circ \sigma.$ If two functions are cohomologous, then
they have the same pressure and the same Gibbs state. The function
$a \rightarrow P(af)$ is continuous \cite{ruel}. For every H\"{o}lder
continuous $f$ and for every integer $n \geq 1,$ the pressure functional
satisfies $P(S_nf) = n P(f).$  Most importantly, if $f < 0,$ then there
exists a unique constant $\delta > 0$ such that
\begin{equation}
       P(\delta f) = 0.
\end{equation}
The pressure and the measure-theoretic entropy of the Gibbs state are related
to each other by the {\em Variational Principle} (\cite{bowe}, Theorem 1.22)
since prominently the Gibbs state for $t f,$ any $ t > 0,$ 
is the {\em unique} equilibrium state for $- t f.$ Choosing $t$ to be
the $\delta$ that nullifies the pressure yields
\begin{equation}
   \label{variationalprinc}
      h(\mu_{\delta f}) =  P(\delta f)   - \int \delta f \, d \mu_{\delta f}  
                        = - \int \delta f \, d \mu_{\delta f}.
\end{equation}

\bigskip
{\bf Counting problems.}
For $0 < \zeta < 1,$ define
\begin{equation}
     \label{astar}
      \cA^* (\zeta) = \bigcup_{n=1}^{\infty}
                     \{ i \in \cG_n : S_n f(i) \leq 
                \log \zeta
                       \mbox{ and }   S_k f(i) > \log \zeta , \, 
               \forall k < n \}
\end{equation}
to be the set of comparable finite sequences relative 
to $f$ at scale $\zeta.$ The set  $\cA^*(\zeta)$ consists 
of those finite sequences of possibly different 
lengths $n$ such that $S_nf$ takes a value just below $\log \zeta.$ 
Observe that, since $f$ is bounded, for any 
$i \in \cA^*(\zeta),$ $S_nf(i)$ differs from $\log \zeta$
by at most $\vert \! \vert f \vert \! \vert_{\infty} < \infty.$
Note that for every sequence $i \in \Sigma,$ there exists a unique $n$ such 
that the finite sequence $i_1 i_2 \ldots i_n$ is an element of  
 $\cA^*(\zeta).$  The following result, borrowed from \cite{hula2},
is presented along with a proof because the latter offers some  
insights that will be of use in the subsequent discussion. 

\begin{proposition}[Hueter and Lalley \cite{hula2}, Proposition 2.1]
  \label{arhocard}
  Let  $ \delta > 0$ be the unique positive number such that $P(\delta f)=0.$ 
  Then there are some suitable positive finite constants $C_3$ and $C_4$ 
  so that, as $\zeta \rightarrow 0,$ the cardinality of  $\cA^*(\zeta)$
is given by
  \begin{equation}
       C_3 \zeta^{- \delta} \leq \vert \cA^*(\zeta) \vert 
 \leq C_4 \zeta^{- \delta}.
  \end{equation} 
 \end{proposition}

\begin{proof}
Each element $i \in \cA^*(\zeta)$ gives rise to a cylinder set $\Gamma(i). $
Moreover, the cylinder sets $\{ \Gamma(i): i \in  \cA^*(\zeta) \} $
are pairwise disjoint and their union is the entire sequence space $\Sigma.$
Consequently,
$$
     \sum_{i \in  \cA^*(\zeta)} \, \mu_{\delta f}(\Gamma(i)) = 1.
$$      
By the definition of $\delta$ and the Gibbs state, there are constants
$ 0 < C_1 \leq C_2 < \infty$ such that, for every $i \in  \cA^*(\zeta),$ 
the measure of the cylinder set $\Gamma(i)$ satisfies
$$
     C_1 \exp \{ \delta S_n f(i) \} \leq \mu_{\delta f}(\Gamma(i)) \leq
       C_2 \exp \{ \delta S_n f(i) \},
$$
where $n$ is the length of $i.$ But by the defining property of         
$ \cA^*(\zeta),$ there is a positive constant $c_3 \leq 1$ such that
for every $i \in  \cA^*(\zeta), $ 
$$
   c_3 \zeta \leq \exp \{ S_n f(i) \} \leq \zeta. 
$$    
Assembling the last three displayed formulae yields
$$
   c_4  \sum_{i \in  \cA^*(\zeta)} \, \zeta^{\delta} \leq 1 \leq
              C_2 \sum_{i \in  \cA^*(\zeta)}  \, \zeta^{\delta} 
$$              
for a suitable constant $c_4,$ which proves the advertized inequalities. 
\end{proof}                        

``Most" sequences in $\cA^*(\zeta)$ are approximately
$\mu_{\delta f}$-distributed.
Indeed, by the Birkhoff ergodic theorem, 
for every H\"{o}lder continuous function 
$g: \Sigma \rightarrow {\bf R}$ and every $\epsilon > 0,$
$$
\vert \{ i = i_1 i_2 \ldots i_n \in \cA^*(\zeta): \max_{0 \leq t \leq 1}
                   \vert \frac{S_{[nt]} g(i)}{n}\,  - \,
                     t \, \int\,  g \, d \mu_{\delta f}  \vert > \epsilon \}
                      \vert = o(\zeta^{-\delta}).
$$ 
Therefore, if we let $ n_{\zeta} = \log \zeta / \int g \, d\mu_{\delta f},$
then ``most" sequences in  $\cA^*(\zeta)$ have lengths between 
  $n_{\zeta} (1- \epsilon)$ and  $n_{\zeta} (1+ \epsilon).$  
These observations indicate that 
the set $\cA^*(\zeta)$ is nearly a set of sequences of length $n_{\zeta}$
that are approximately ``generic" for the measure $ \mu_{\delta f}.$ 
In light of the Shannon-McMillan-Breiman theorem, the cardinality of the
latter is approximately $\exp \{ h_{\mu_{\delta f}} n_{\zeta} \},$ 
where $ h_{\mu_{\delta f}}$ denotes the
entropy of the measure $ \mu_{\delta f},$ whereas in view of Proposition 
\ref{arhocard}, the cardinality of the former set is of the order
$ \zeta^{-\delta}.$ This confirms the variational principle
in (\ref{variationalprinc}), which implies that
$$
     \zeta^{-\delta} = e^{n_{\zeta} h_{\mu_{\delta f}} } \, .
$$ 

Next, we state a result of \cite{lal0}, Proposition 5.1,
an extension of the Perron-Frobenius theorem.
For any matrix $M,$ let $\vert \! \vert M \vert \! \vert = 
 \sup_{v \not \equiv 0}
( \vert M v \vert / \vert v \vert)$ denote the usual matrix norm.

\begin{theorem}[Lalley \cite{lal0}, Proposition 5.1]
  \label{existencepotential}
Assume that $ M_x $ is a nonnegative aperiodic $n \times n$ matrix and 
$x \rightarrow M_x$ is a H\"{o}lder continuous function (with some
exponent) on $\cAZ.$ Then there exist constants $C < \infty$ and
$0 < \alpha < 1$ and H\"{o}lder continuous functions
$\varphi, \gamma: \Sigma \rightarrow {\bf R}$ and $v,w: \Sigma \rightarrow
\cP_+ = \{ u \in {\bf R}^n : \sum_{i=1}^n u_i = 1 \mbox{ and } u_i >0 
\mbox{ for every } 1 \leq i \leq n \}$ such that for every 
 $\xi = (x_k)_{k=- \infty}^{\infty} \in \Sigma$ and every integer $n > 0,$ 
\begin{equation}
  \label{holderprod}
\vert \! \vert e^{-S_n \varphi(\xi)} M_{x_1} M_{x_2} \ldots M_{x_n}
                  - \gamma(\sigma^n \xi) v(\xi) w(\sigma^n \xi)^t
 \vert \! \vert
                  \leq C \alpha^n,
\end{equation}
where 
\begin{eqnarray}
  \label{matnorm}
      \gamma(\xi) &=& 1/ w(\xi)^t v(\xi),  \\
   \label{potential}
      M_{x_1} v(\sigma \xi) & = & e^{\varphi(\xi)} v(\xi),  \\
\mbox{and} \mbox{} \qquad \qquad \qquad & & \nonumber \\
     \label{lefteigen}
    w(\sigma^{-1} \xi)^t M_{x_1} & = & e^{\varphi( \sigma \xi)}
                            \frac{\gamma(\xi)}{\gamma(\sigma^{-1} \xi)}     
                          w(\xi)^t .
\end{eqnarray}
Both $v(\xi)$ and $\varphi(\xi)$ are functions of the ``forward" coordinates
$x_1, x_2, \ldots$ and $w(\xi)$ is a function only of the ``backward" 
coordinates $\ldots , x_{-1} , x_0.$
\end{theorem}

Thus, the Perron-Frobenius theorem has $x \rightarrow M_x$ 
a constant function. We will consider the special case when the
matrices $M_x$ are $2d \times 2d,$ thus, the size of the matrices being 
fixed and
independent of the number of terms in the matrix product in 
(\ref{holderprod}). 
The proof in \cite{lal0} carries over unmodified. 
Observe that, for each $n,$ $v(\xi)$ and $w(\xi)$ are right and left 
eigenvectors, respectively, of $ M_{x_1} M_{x_2} \ldots M_{x_n} $ 
with associated eigenvalue $\exp \{ S_n \varphi(\xi) \}.$


\subsection{Critical Exponent $r_u$} 

Recall that $\cG_m$ denotes the set of all vertices at distance
$m$ from the root vertex. For every $\bfl,$ 
define
\begin{equation}
   \label{symmetryindex}
     r_u = r_u (\bfl) = \inf \{ r>0: \sum_{m=0}^{\infty} \sum_{x \in \cG_m} 
                      \,  u_x(\bfl)^r < \infty \}.
\end{equation}
This exponent $r_u (\bfl)$ takes some finite value for $\bfl \in 
\mbox{int}(\cK)$ and is nondecreasing in each $\lambda_j$ because
$u_x(\bfl)$ is nondecreasing in each $\lambda_j.$
For instance, it is easy to see that,    
for $\bfl \not \in \cR_1,$ we must have 
$r_u (\bfl) \geq 1.$ Some intuitive values yet subtle to prove
are as follows:
we will ultimately show (Corollary \ref{criticalexponents}) that,
for $d>1,$ we obtain $r_u(\bfl)=1$ 
for $\bfl \in \overR_1 \cap \overR_2$ and $r_u(\bfl)=2$ 
for $\bfl \in \overR_2 \cap \overR_3.$ 
We begin to view some easy but crucial facts.

\begin{lemma}
   \label{dominatedu}
For each $t > r_u(\bfl),$ we have
$$
 \sum_{m=0}^{\infty} \sum_{x \in \cG_m} 
                      \,  u_x(\bfl)^t  < \infty .
$$
\end{lemma}

\begin{proof}
This is readily concluded from the inequality $u_x^t \leq u_x^{r_u}.$
\end{proof}

\begin{proposition}
   \label{rumono}
Let $\bfl \in \cK_0$ with 
$0 < \lambda_c = \min_{a \in \cA} \lambda_a$
and let $\tilde{\bfl}$ be an infection parameter
that corresponds to the $p$-thinned percolation 
structure for some $0 < p < 1$ relative to $\lambda_c$ such that
$ 0 < \tilde{\lambda}_c  = p \lambda_c/(1+ q \lambda_c) < \lambda_c. $
Thus, $\tilde{\bfl}$ lies in a direction of decrease for $\bfl.$ 
Then 
   $$r_u(\tilde{\bfl}) < r_u(\bfl).$$ 
\end{proposition}    

\begin{proof}
Fix some $0 < p < 1,$ let 
$\bfl \in \cK_0,$ and let $\tilde{\bfl}$ be a
parameter that corresponds to the infection parameter
gotten from the $p$-thinned percolation structure
relative to $\lambda_c$ in the
sense that the percolation structure is preserved. Choose
$\tilde{\bfl}$ such that the contact process has recovery rate $1.$
In view of Proposition \ref{monotone},
there is some constant $0 < \omega < 1$
such that for every $n$ and $ x \in \cG_n,$
$u_x (\tilde{\bfl}) \leq u_x(\bfl) \omega^n .$
Moreover by Lemma \ref{dominatedu}, for each $t > r_u(\bfl),$ we have
$
 \sum_{m=0}^{\infty} \sum_{x \in \cG_m} 
                       u_x^t(\bfl)  < \infty.
$
Consider
\begin{eqnarray}
    \sum_{m=0}^{\infty} \sum_{x \in \cG_m} 
            u_x(\tilde{\bfl})^{r_u(\bfl)- s}  & \leq &
              \sum_{m=0}^{\infty} \sum_{x \in \cG_m} 
             (\omega^m)^{r_u(\bfl)- s} \, u_x(\bfl)^{r_u(\bfl)- s}  
                \nonumber \\
        & = &  
          \label{usumone}
            \sum_{m=0}^{\infty} \sum_{x \in \cG_m} 
             [\omega^{m(r_u(\bfl)- s)}/ u_x(\bfl)^{\delta} ] \, \,
               u_x(\bfl)^{r_u(\bfl)- s+ \delta}  
\end{eqnarray}             
for some reals $s, \delta > 0.$ 
In view of parallel arguments, resting on the minimal
positive infection rate $\lambda_c,$ to the ones employed in the proof
of Lemma \ref{periodicmaximal} for the maximal infection rate
along with (\ref{betalimit}), for any $\varepsilon > 0,$ there
exists some $N$ such that for every $m > N $ and $x \in \cG_m,$
we have $u_x  \geq  (1 - \varepsilon)^m [\min_{a \in \cA} \beta_a]^m  
      = (1 - \varepsilon)^m (\beta_c)^m > (\beta_c/2)^m$
with $\beta_c > 0$ because $ \lambda_c >0.$  
Since $\omega < 1,$ for every $s < r_u(\bfl),$
we can choose $\delta$ small enough so that 
$$
     \omega^{m(r_u(\bfl)- s)}/ u_x(\bfl)^{\delta}
      < [\omega^{r_u - s} (2/ \beta_c)^{\delta}]^m \leq \gamma^m
$$
for some $0 < \gamma < 1$ and sufficiently large $m.$    
In addition, for each $ s < \delta,$ we obtain
  $ \sum_{m=0}^{\infty} \sum_{x \in \cG_m} 
            u_x(\bfl)^{r_u(\bfl)- s+ \delta}  < C$ 
for some positive finite constant $C.$ Observe that there exist
sufficiently small $s, \delta > 0$ that can satisfy the former and 
latter conditions.
Hence, all sufficiently advanced terms in the summation
$\sum_{m=0}^{\infty}$ on the righthand side of 
(\ref{usumone}) are bounded
above by $ \gamma^m C $ for some positive finite constant $C.$ 
Therefore, summation over a geometric series leads us to conclude that
there is some positive $s$ such that 
\begin{eqnarray*}
  \sum_{m=0}^{\infty} \sum_{x \in \cG_m} 
            u_x(\tilde{\bfl})^{r_u(\bfl)- s} < \infty.
\end{eqnarray*}
Consequently, $r_u(\tilde{\bfl}) < r_u(\bfl),$ as desired.
\end{proof}

\begin{proposition}
  \label{rubound}
Suppose that $\bfl' \in \overline{\cK} \cap 
\overline{\cK^c}$ and that
$\bfl$ lies in the interior $\mbox{int}(\cK)$ and in a direction
of decrease for $\bfl'.$ Then 
\begin{equation}
   \label{rulessthantwo}
      r_u(\bfl) < r_u(\bfl').
\end{equation}
\end{proposition}        
  
\begin{proof}
Assume that $\bfl$ lies in the interior $\mbox{int}(\cK).$ 
Then there exists a ball $B$ in ${\bf R}^d,$ centered at $\bfl$ 
which is completely inside $\cK.$ 
Since each ball contains a suitable
multiple of each unit vector in ${\bf R}^d,$ for every $0<p<1,$
the ball $B$ contains a line segment $L_p$ 
in the direction of decrease for $\bfl,$
equivalently, increase that contains $\bfl,$ 
defined by the $p$-thinning of the percolation structure relative to
the minimal infection rate.
Hence by Proposition \ref{rumono}, 
each point $\bfl_*$ of $L_p$ 
that has strictly smaller distance
to the boundary $ \overline{\cK} \cap \overline{\cK^c}$ 
than $\bfl$ has $r_u(\bfl) < r_u(\bfl_*).$ 
But $r_u(\bfl_*) \leq r_u(\bfl') $
for $\bfl' \in \overline{\cK} \cap \overline{\cK^c}.$ 
Hence, $r_u(\bfl) < r_u(\bfl_*) \leq r_u(\bfl'),$
as advertised.
\end{proof}


\subsection{The Infection Probabilities satisfy a H\"{o}lder Condition} 
 
Perhaps the most interesting aspect of this and the next section
is how to prove a H\"{o}lder condition for the collection of 
infection probabilities in $\cK_0$ so that the existence of a 
Gibbs state is guaranteed for the contact process and the theory
described in Section 3.1 comes to fruit.
 
\begin{proposition}
  \label{hoelderquot}
There are constants $0< \gamma < 1$ and $0 < C < \infty$    
such that 
for every $\bfl \in \mbox{int}(\cK_0),$  every $\rho > r_u(\bfl),$ 
every integer $k>0,$
and for every sufficiently large integer $n,$  
\begin{equation}
   \sum_{x_{k+1} x_{k+2} \ldots x_{n} \in \cG_{n-k}}
       \, \left(
        \frac{u_{x_1 x_2 \ldots x_n}(\bfl)}{
                  u_{x_1 x_2 \ldots x_k}(\bfl)} 
                  \right)^{\rho} < C \gamma^n.
\end{equation} 
\end{proposition}

\begin{proof}
Fix $\bfl \in \mbox{int}(\cK_0).$ 
Assume that $\lambda_b \geq \lambda_i $ for every $i \in \cA.$
Let $a_{\epsilon}$ denote the probability for the isotropic
contact process with infection parameter $\epsilon>0$ and initial infection
at the root vertex that between time $0$
and $1$ there is no recovery mark * 
at the root $1$ and that, for some $i \in \cA,$ the root vertex
infects its neighbour $i.$  
Thus, $ a_{\epsilon} = e^{-1}(1-e^{-\epsilon}) > 0.$ Fix  $\epsilon$
such that each $\lambda_i \geq \epsilon.$ Then for
every $x \in \cG_k,$ we obtain
$$ 
   u_x (\bfl)  \geq (a_{\epsilon})^k. 
$$ 
Fix some integers $k,M>0.$ 
Since $\rho > r_u(\bfl)=r_u,$ we can choose some $\delta >0$
such that $ \rho - r_u - \delta >0.$ 
As a consequence, by Lemma \ref{periodicmaximal}, for each  
$x \in \cG_k,$ every $n \geq kM,$ 
and each $y \in \cG_n,$
\begin{eqnarray*}
   \frac{u_y(\bfl)^{\rho - r_u - \delta}}{u_x(\bfl)^{\rho}}      
          & \leq &  \frac{\beta_b(\bfl)^{n(\rho -r_u - \delta)}}
                         {(a_{\epsilon})^{k \rho}} \\
          & \leq &  ( \frac{\beta_b(\bfl)^{\rho -r_u - \delta}}
                         { (a_{\epsilon})^{\rho/M}}
                         )^n.
\end{eqnarray*}
Since, by choosing $M$ sufficiently large,
$ \beta_b(\bfl)^{\rho -r_u - \delta}/  (a_{\epsilon})^{\rho/M} 
                    \leq \gamma < 1$ 
for some $\gamma,$ and thus,                        
the righthand side of the last display is bounded above by $\gamma^n.$
Combining this with Lemma \ref{dominatedu} yields,
for every $k M \leq n,$ 
\begin{eqnarray*}
 \sum_{x_{k+1} \ldots x_{n} \in \cG_{n-k}}
       \, ( \frac{u_{x_1 x_2 \ldots x_n}(\bfl)}{
                  u_{x_1 x_2 \ldots x_k}(\bfl)} )^{\rho} 
           & = & 
             \sum_{x_{k+1}  \ldots x_{n} \in \cG_{n-k}}
       \,  \frac{(u_{x_1 x_2 \ldots x_n}(\bfl))^{\rho-r_u - \delta}}{
                 ( u_{x_1 x_2 \ldots x_k}(\bfl))^{\rho}}  \, 
             \, \, u_{x_1 x_2 \ldots x_n}(\bfl)^{r_u + \delta} \\*[0.14cm] 
           & \leq &  \gamma^n C
 \end{eqnarray*}
 for some positive finite constant $C.$ This ends the proof.           
\end{proof}


\subsection{Transition Matrices and Potential Functions}

The H\"{o}lder continuity allows us to define a Gibbs state
by means of certain matrices. Whereas there are numerous 
reasonable choices of definition for transition matrices 
in the context of the contact process on $\cT,$ the one we
propose is natural and has some nice interpretations for some
exponents.  
If nothing else is said, we shall assume that
each $\lambda_k >0$ and that the generators with zero infection rates 
have been eliminated from $\cAp.$ 

For each $x_1 x_2 \ldots x_{k-1} \in \cG_{k-1}$ 
and for all $i,j \in \cA$ with $i \not = x_{k-1}^{-1},$ 
define 
\begin{equation}
   \label{findexset}
    \cF_{ij}= \cF(x_1 \ldots x_{k-1},i,j)  =
        \{ x_{k+1} \ldots x_{n-1} \in \cG_{n-k-1}: \,
              x_{k+1} \not = i^{-1}, \, x_{n-1} \not = j^{-1} \} 
\end{equation}
and let $ \cF(x_1 \ldots x_{k-1},i,j)$ be the empty set if 
$i  = x_{k-1}^{-1}.$
Then, for each real $\rho > 0,$ 
all integers $n-1 > k \geq  1,$ and 
each $x_1 x_2 \ldots x_{k-1} \in \cG_{k-1},$ 
define the $2d \times 2d$ matrix
$H_{\rho} (n ; x_1 x_2 \ldots x_{k-1}; \bfl)$ by 
\begin{equation}
   \label{hmatrices}
    (H_{\rho} (n ; x_1 x_2 \ldots x_{k-1}; \bfl))_{ij} = 
         \sum_{\textstyle x_{k+1} x_{k+2} \ldots x_{n-1} \in \cF_{ij}}
             \, \left(\frac{u_{x_1 x_2 \ldots x_n}(\bfl)}{
                  u_{x_1 x_2 \ldots x_{k-1}}(\bfl)} \right)^{\rho} 
\end{equation}                       
for all $i,j \in \cA,$ thus, equal zero for $i = x_{k-1}^{-1}.$
Note that for $\rho=1,$ $k=1,$ the $(i,j)$-entry of 
$H_1 (n;1; \bfl) $ equals the expected number of vertices 
at distance $n$ from
the root $1$ that are ever to be infected whose word 
representation begins in $i$ and ends in $j.$
Similarly,  for $\rho=2,$ $k=1,$ the $(i,j)$-entry of 
$H_2 (n; 1; \bfl) $ is a lower bound for 
the expected number of vertices  at distance $n$ from
the root whose word representation begins in $i$ and ends in $j$
that are ever infected and upon infection are ever to send 
the infection back to the root. 
Moreover for each real $\rho > 0,$
$ {\bf 1}^t H_{\rho}(n; 1; \bfl) {\bf 1} = \sum_{x \in \cG_n}  u_x^{\rho},$
where  ${\bf 1}$ denotes the $2d$-vector of all ones.  
Importantly, we recover
\begin{eqnarray}
    \label{addingentries}
 {\bf 1}^t H_1 (n;1; \bfl) {\bf 1} & =&  \sum_{x \in \cG_n}  u_x, \\
  {\bf 1}^t H_2 (n;1;\bfl) {\bf 1} & =&  \sum_{x \in \cG_n}  u_x^2.
    \nonumber
\end{eqnarray}
By subadditivity, it is an easy observation that
$$ \vert \! \vert 
      H_{\rho} (n+k-1 ; x_1 x_2 \ldots x_{k-1}; \bfl) \vert \! \vert  
     \geq  \vert \! \vert  H_{\rho} (n;1; \bfl) \vert \! \vert 
$$
for every $k \geq 1,$ in other words,
there might be positive probability that, 
by the time of infection of vertex $ x_1  \ldots x_{k-1},$ some
other infected vertices line up in the complement of 
$\cT(x_1 \ldots x_{k-1})$ to reinforce the infection.
Observe that for every $\bfl \in \cK_0$ and each $\rho> 0,$ 
the matrix $H_{\rho}(n; x; \bfl)$ is
an aperiodic, irreducible, and nonnegative matrix, thus, a Perron-Frobenius
matrix. Therefore, the Perron-Frobenius theorem tells us that
$H_{\rho}(n; x; \bfl)$ 
has a largest positive eigenvalue. We will see in Proposition 
\ref{leadeva} below that 
the lead eigenvalue equation can be solved explicitly.  
In light of Proposition \ref{hoelderquot}, for $\rho > r_u(\bfl),$
the matrix $ H_{\rho}(n;x; \cdot)$ satisfies a H\"{o}lder condition.
The mapping $x \rightarrow H_{\rho}(n;x; \cdot)$ 
is given by $x_1 x_2 \ldots x_{k-1}
\rightarrow H_{\rho} (n; x_1 x_2 \ldots x_{k-1}; \cdot ).$ 
 
\begin{proposition}
   \label{hhoelder}
For every $\bfl \in \mbox{int}(\cK_0),$ 
and each real $\rho > r_u(\bfl),$ 
the matrix $x \rightarrow H_{\rho}(n;x; \cdot)$ is H\"{o}lder
continuous (with some exponent).
\end{proposition}  

\begin{proof}
Fix such $\bfl$ and any $\rho > r_u(\bfl).$ 
Since the matrix norm $\vert \! \vert A- B \vert \! \vert$
of the difference of two nonnegative matrices $A$ and $B$
is bounded above by the maximum $\max(\vert \! \vert A \vert \! \vert,
\vert \! \vert B \vert \! \vert)$ of the individual matrix norms,
and in turn, the norm of any matrix $\vert \! \vert A \vert \! \vert$
is bounded by the maximum of the absolute values of the $A$-matrix
entries times its dimension, it suffices to verify
that the maximal absolute values
of the entries of the differences 
  $ H_{\rho} (n; x_1 x_2 \ldots x_{k-1}; \cdot ) 
     - H_{\rho} (n; x_1 x_2 \ldots x_{k+l-1}; \cdot )$ and 
      $ H_{\rho} (n; x_1 x_2 \ldots x_{k-1}; \cdot ) 
      -  H_{\rho} (n+l; x_1 x_2 \ldots x_{k-1}; \cdot )$ 
are bounded above by $C \gamma^n$ for some constants $ 0 < \gamma < 1$
 and $ 0 < C < \infty,$ for all integers $ k, l > 0$ and
for every sufficiently large $n.$        
But it is an elementary exercise by the results in Proposition
\ref{hoelderquot} to show that each entry of the latter four
matrices is bounded above by  $C \gamma^n,$  where $C$ and $\gamma$
are independent of $ k, l,$ and $n.$ This finishes our proof.
\end{proof}

Recall that $\Sigma$ denotes the set of doubly infinite reduced words
from $\cA.$ Notice that in the subsequent discussion, we shall not
always be careful to distinguish $\Omega$ and $\Sigma.$ 
Now we are ready to define the {\em potential function} 
$\varphi_{\rho}(x) =\varphi_{\rho; \bfl}(x)  $
for every $\rho > r_u(\bfl)$ by applying Theorem \ref{existencepotential} 
with $ M_{x_1} \cdots M_{x_{n-k}} =
 H_{\rho}(n; x_1 \ldots x_{k-1}; \cdot),$ $V=v$ and
$W=w$ with $V,W: \Sigma \rightarrow \cP_+ = \{ u \in {\bf R}^{2d} : 
\sum_{i=1}^{2d} u_i = 1 \mbox{ and } u_i >0 
\mbox{ for every } 1 \leq i \leq 2d \}.$ Then for every
$x= \ldots x_1 x_2 \ldots \in \Sigma,$ every integer $k>0,$ 
and sufficiently large $n,$ define
\begin{equation}
   \label{hpotential}
     H_{\rho}(n; x_1 x_2 \ldots x_{k-1}; \bfl) 
      V(\sigma^{n-k} x ) =  e^{\textstyle S_{n-k} \varphi_{\rho}(x)} V(x).
\end{equation}
A similar equation represents the left eigenvector $W(x)$ associated
with  $ H_{\rho}(n; x_1 \ldots x_{k-1}; \bfl)$ (for more details, 
see Proposition 5.1, \cite{lal0}). 
The H\"{o}lder continuity of
$\varphi_{\rho} = \varphi_{\rho; \bfl}$ 
immediately follows from the H\"{o}lder continuity of $V$ and the map 
$ x \rightarrow  H_{\rho} (n; x_1 \ldots x_{k-1}; \bfl).$
To facilitate notation,
we will often just write $\varphi_{\rho}.$
Observe that $H_{\rho}(n; x_1 x_2 \ldots x_{k-1}; \bfl)$ has
eigenvalue  $ \exp \{ S_{n-k} \varphi_{\rho}(x) \} $ with corresponding
right and left eigenvectors $V(x)$ and $W(x).$ Also,
note that for $\rho > r_u(\bfl),$ we have $\varphi_{\rho} < 0.$
Combining (\ref{holderprod}) with definitions (\ref{hmatrices}) and
(\ref{hpotential}) yields
\begin{equation}
   \label{convestimate}
\sum_{x_{k+1}  \ldots x_{k+n} \in \cG_{n}} 
         \left( \frac{u_{x_1 x_2 \ldots x_{k+n}}(\bfl)}{
                  u_{x_1 x_2 \ldots x_{k}}(\bfl)} \right)^{\rho} 
           = C \exp \{ S_n \varphi_{\rho}(x) \} (1+ O(\alpha^k))
\end{equation}
for some $0< \alpha <1,$
for every $x= \ldots x_1 x_2 \ldots \in \Sigma$, all integers $k$ and $n,$ 
and every $\rho>r_u(\bfl),$ where the implicit bound in the $O(\alpha^k)$
term is uniform in $x.$ The constant $C$ may be bounded by 
$0 < C_1 < C < C_2 < \infty$ with $C_i$ independent of $x,$
$k$ and $n.$

For every $x \in \Sigma,$ the functions 
$S_n \varphi_{\rho}(x),$ and thus, 
$H_{\rho}(n+k; x_1 \ldots x_{k-1}; \cdot),$  
depend on the relative frequencies of the generators in the reduced
word $x_1 x_2 \ldots x_{k-1}$ 
and the {\em order} of the letters as well.

\smallskip
{\bf Matrix Entries $b_i.$}
Next we regard to express each matrix 
$ H_{2}(n; x_1 x_2 \ldots x_{k-1}; \bfl)$
 as the $(n-k)$-th power of a Perron-Frobenius matrix 
$B_{2, n, x_1 \ldots x_{k-1}}(\bfl) $ which
 (a) preserves the lead eigenvalue, and 
 (b) has a certain form that specifies the allowed transitions as 
the matrix does in the following result. 

\begin{proposition}
  \label{leadeva}
For any positive number $\rho,$ define the matrix 
$A_{\rho}$ to be the $2d \times 2d$ matrix, indexed by elements
of $\cA,$ whose entries are given by 
\begin{eqnarray}
  \label{amatrix}
            (A_{\rho})_{ij} & = & a_j^{\rho} \qquad  \, \; \, \,
                               \mbox{if }  j \not = i^{-1}, \\
                            & = & 0 \qquad \, \, \quad  
                        \mbox{if }  j = i^{-1} \nonumber
\end{eqnarray}
for some real numbers $a_j \geq 0.$
Assume that $a_j = a_{j^{-1}}$ for each $j \in \cA$ and that
$A_{\rho}$ is a Perron-Frobenius matrix. Then the
lead eigenvalue of $A_{\rho}$ is the unique positive
solution $ \alpha_{\rho}$ of the equation
\begin{equation}
  \label{alphaequation}
        \sum_{i \in \cA} \, \frac{a_i^{\rho}}{
                    \alpha_{\rho} + a_i^{\rho}} = 1.
\end{equation}
Moreover, for each $\rho > 0,$ the function
$ a_i \rightarrow \alpha_{\rho}$ is strictly increasing
for each $i \in \cA.$ 
\end{proposition}

\begin{proof}
Identity (\ref{alphaequation}) is shown by rearranging the eigenvalue
equation $A_{\rho} u = \alpha_{\rho} u .$  If $\alpha_{\rho}$
is the lead eigenvalue, then the vector $u$ must have nonnegative entries,
not all zero. Let $0 < s =  \sum_{i \in \cA} a_i^{\rho} u_i .$
This identity combined with the eigenvalue equation
$ \alpha_{\rho} u_j = \sum_{i \in \cA} a_i^{\rho} u_i - a_j^{\rho} u_j $
yields
$$ 
   u_j = s/ ( \alpha_{\rho} + a_j^{\rho} ).
$$
Multiplying both sides by $a_j^{\rho}, $ adding over all
$i \in \cA,$ and dividing both sides by $s$
gives (\ref{alphaequation}). 
In addition, the strict monotonicity of $ \alpha_{\rho}$ in each $a_i$
is an immediate consequence of the form of the equation 
(\ref{alphaequation}) and the facts that the map
 $x \rightarrow x/(\alpha_{\rho} +x)$ is continuously differentiable and
 has strictly positive derivative for positive $x.$ 
\end{proof}

Now, for any positive number $\rho$ and $\bfl \in \cK,$ 
define a $2d \times 2d$ matrix 
$ B_{\rho} = B_{\rho,n, x_1 \ldots x_{k-1}}(\bfl), $
indexed by elements of $\cA,$ so that   
(a) their entries $b_j(\bfl)^{\rho}= b_j(\bfl, 
x_1 x_2 \ldots x_{k-1})^{\rho} \geq 0$ come in the form
(\ref{amatrix}), that is,
\begin{eqnarray}
  \label{bmatrix}
            (B_{\rho}(\bfl))_{ij} & = & b_j(\bfl)^{\rho}   
                                    \quad 
                               \mbox{if }  j \not = i^{-1}, \\
                            & = & 0 \qquad \, \, \quad  
                        \mbox{if }  j = i^{-1} \nonumber
\end{eqnarray}
and  
(b) the $b_j = b_j(\bfl)$  satisfy (\ref{alphaequation}) 
 with  
\begin{eqnarray}  
   \label{equationsforbj}
    a_j^2 & = & b_j^2, \, \rho=2, \mbox{ and } \alpha_2 
     = \exp \{ S_{n-k} \varphi_2(x)/(n-k) \} . 
\end{eqnarray}        
Hence, the $b_j,$ $j \in \cA,$
link $\varphi_2$ and the matrices $B_{\rho}.$
In this construction, the power $\rho$ in 
(\ref{equationsforbj}) needs to be chosen sufficiently large
to guarantee the existence of $\varphi_{\rho}$ and to assure
that the $b_j$ are well defined. This is possible since one
shows that $r_u(\bfl) < \infty$ for $\bfl \in \cK.$ 
Our choice $\rho=2$ is the smallest possible because
it will turn out later \mbox{(Corollary \ref{valuesattrans})} 
that $r_u(\bfl) <2$ for $\bfl \in \mbox{int}(\cK).$
Thus, $ B_2^{n-k} $ has lead eigenvalue 
$ \exp \{ S_{n-k} \varphi_2(x)\}.$ Both $B_1$ and $B_2$ are 
Perron-Frobenius matrices. The entries $b_j$ of $B_1$ and 
$b_j^2$ of $B_2,$ respectively, 
depend on $n$ and the sequence $x_1 x_2 \ldots x_{k-1}.$
We shall not always be careful to indicate these 
dependencies and assume that they are clear from the context.
Observe that, if $\beta_i(\bfl)=1,$ then for each $x \in \cG_{k-1},$
we can choose $b_i(\bfl, x) =1,$ due to the construction of the
matrix $B_1,$ equation (\ref{alphaequation}), and the fact that
$u_y$ is independent of the choice of the beginning finite segment 
$x_1 x_2 \ldots x_{k-1}$ of an infinite word 
$y = \ldots x_1 x_2 \ldots. $ 
An analogous argument shows that, if $\beta_i(\bfl)=0,$ we may
choose $b_i(\bfl, x)=0$ for each $x \in \cG.$ Of course,
the converse also holds.
The following observation is a summary.

\begin{lemma}
   \label{bilessthanone}
If  $b_i(\bfl,x) < 1$ for each $x \in \cG$ and $i \in \cA,$
then $\bfl \in \cK.$ Moreover, $b_i(\bfl, x) > 0$ for
each $x \in \cG$ and $i \in \cA$ if and only if $\beta_i(\bfl)>0.$
\end{lemma}

The next result shows that the functional $ S_{n-k} \varphi_{2}$
is compatible with the construction of the matrix $B_2 .$

\begin{proposition}
   \label{additivefct}
For every sufficiently large $n-k > 0,$ 
the functional $S_{n-k} \varphi_2$ is additive, 
that is, the function 
$\varphi_2: \Sigma \rightarrow {\bf R}$ satisfies
$$
     S_{n-k} \varphi_2 = \varphi_2 + \varphi_2 \circ \sigma
      + \varphi_2 \circ \sigma^2 + \ldots + \varphi_2 \circ 
       \sigma^{n-k-1},
$$
where $\sigma$ denotes the forward shift on $\Sigma.$ 
\end{proposition}

\begin{proof}
Let $x= \ldots x_1 x_2 \ldots \in \Sigma.$
Since $\alpha_{2}^{n-k}$ is the lead eigenvalue of 
 $ H_{2}(n; x_1 x_2 \ldots x_{k-1}; \bfl)$ and
 $  B_{2}^{n-k},$ and thus, the right eigenvector $u$ of
  $  B_{2}^{n-k}$ is dominated by
a positive finite constant multiple of the right eigenvector
of  $ H_{2}(n; x_1 x_2 \ldots x_{k-1}; \bfl)$ and vice versa, 
the (right) eigenvectors of $  B_{2}^{n-k}$ and  
     $ H_{2}(n; x_1 x_2 \ldots x_{k-1}; \bfl)$ differ by at
most a positive finite constant multiple. Additionally,
in the eigenvalue equation
 $  B_{2}^{n-k} u = \alpha_{2}^{n-k} u,$ the same eigenvector $u$
 is relied on for every $n-k.$  Hence, iteration brings
\begin{eqnarray*}
   \exp \{ S_{n-k} \varphi_{2}(x)\}  V(x) c_1
        & = &  H_{2}(n; x_1 x_2 \ldots x_{k-1}; \bfl)      
                               \, V(\sigma^{n-k} x ) c_1\\
         & \leq &  B_{2}^{n-k} u  \\
               & = & \alpha_{2} \,  B_{2}^{n-k-1}  u   \\
               & = & \alpha_{2}^2 \,  B_{2}^{n-k-2}  u   \\  
                & = & \ldots \ldots \\     
                & = & \alpha_{2}^{n-k}   u   \\ 
                & = & \exp \{ 
                       \sum_{i=0}^{n-k-1} 
                       \varphi_{2}(\sigma^{i} x) \} u. 
\end{eqnarray*} 
Similarly, the reverse inequality holds if the constant $c_1$
is replaced by a constant $c_2.$ 
Since for sufficiently large $n,$ the multipliers
$c_i$ are negligible as compared to the argument of
the exponential function, this completes the proof.                
\end{proof}


\subsection{Invariant Measures}

Next we recollect some facts from Section 3.1.
An elementary exercise shows
that $r_u(\bfl) < \infty$ for $\bfl \in \cK.$ 
We have seen that  $\varphi_{\rho} < 0$
for each $\rho > r_u(\bfl).$ 
Thus, there is a unique shift-invariant probability measure
$\mu_{\varphi_{\rho}}$ on the Borel sets of $\Sigma$ such
that there is some real constant $P(\varphi_{\rho})$ and some
constants $0 < C_3 \leq C_4 < \infty$ such that, for each $m,$
$i \in \cG_m,$ and cylinder set $ \Gamma(i),$ we have
\begin{equation}
  \label{invmeasrho}
 C_3 \exp \{ S_m \varphi_{\rho}(j) - m P(\varphi_{\rho}) \}
   \leq \mu_{\varphi_{\rho}}(\Gamma(i)) \leq C_4 
        \exp \{ S_m \varphi_{\rho}(j) - m P(\varphi_{\rho}) \}
\end{equation}          
for every $j \in \Gamma(i).$ 
From (\ref{astar}), recall
\begin{equation}
    \label{azeta}
      \cA_{\rho}^*(\zeta) = \bigcup_{n=1}^{\infty} \{ i \in \cG_n:
       S_n \varphi_{\rho}(i) \leq 
                \log \zeta
                       \mbox{ and } S_k \varphi_{\rho}(i) > \log \zeta , \, 
               \forall k < n \} 
\end{equation}
for every $1 > \zeta >0$ 
and that, for $\rho > r_u(\bfl),$ the sum
$ \sum_{i \in \cA_{\rho}^*(\zeta)} 
   \mu_{\delta \varphi_{\rho}}(\Gamma(i))=1,$
where $\delta = \delta_{\rho} > 0$ is such that the pressure 
$P (\delta \varphi_{\rho}) = 0.$  
Define the expectation of $\varphi_{\rho}$ 
with respect to the probability measure $\mu_{\delta \varphi_{\rho}}$ by
\begin{equation}
   \label{phimeans}
 \overline{\varphi}_{\rho} = \int \varphi_{\rho} \, 
                           d \mu_{\delta \varphi_{\rho}}.
\end{equation}
The Birkhoff ergodic theorem implies that, for every $\epsilon > 0,$
 \begin{equation}
     \label{birkhoff2} 
     \vert \{ i = i_1 i_2 \ldots i_n \in
           \cA_{\rho}^*(\zeta): \max_{0 \leq t \leq 1}
                   \vert \frac{S_{[nt]} \varphi_{\rho}(i)}{n}\,  - \,
                     t \,  \overline{\varphi}_{\rho} 
                     \vert > \epsilon \}
                      \vert = o(\zeta^{-\delta}), 
\end{equation}
thus, ``most" sequences in $ \cA_{\rho}^*(\zeta)$ have lengths
between  $ n_{\zeta}(1 - \epsilon)$ and $ n_{\zeta} (1+ \epsilon),$
where 
\begin{equation}
   \label{genlength}
    n_{\zeta} = \frac{\log \zeta}{ \overline{\varphi}_{\rho}}. 
\end{equation} 

Next, consider the limit points $ \omega = \omega_1 \omega_2  \ldots 
 \in \Omega$ in the weak survival phase.
 Define the distribution $\mu_n$
 under $P$  of the process  $\omega_n , \omega_{n+1}, \ldots,$ 
 that is, for any Borel set $V \subset \Omega,$ let
 \begin{equation}
   \label{distprocess}
    \mu_n(V) = P \{ (\omega_n , \omega_{n+1}, \ldots) \in V \} .    
\end{equation} 
Now by definitions (\ref{hmatrices}), 
(\ref{hpotential}), and (\ref{bmatrix}) along with the remark at
the end of Section 3.4, 
there is some constant  $0 < \alpha <1$ such 
that for every sufficiently large $m,$
$$
 P \{ \omega_j = x_j \mbox{ for each } 1 \leq j \leq m \} 
                       = C (\prod_{j=1}^m b_{x_j}) (1 + O(\alpha^m))
$$  
for every $x \in \Sigma$ and some constant $C$ that may be bounded by 
$0 < c_1 < C < c_2 < \infty$ with $c_i$ independent of $x$
and  $m.$
If we define the function $\varphi : \Sigma \rightarrow \bR$
to be $ \varphi (\ldots x_1 x_2 \ldots ) = \log b_{x_1}, $ then for every
$x \in \Sigma,$
\begin{equation}
   \label{stataddfct}
  P \{ \omega_j = x_j \mbox{ for each } 1 \leq j \leq m \} 
                  = C \exp \{ S_m \varphi(x) \}  (1 + O(\alpha^m)). 
\end{equation} 
Note that this function $\varphi$ is different from the functions
$\varphi_{\rho}$ defined in (\ref{hpotential}) above.
Since $\varphi <0,$ there exists a unique $\delta >0$ such that
the pressure $P(\delta \varphi) =0.$ Also, by (\ref{gibbsf}),
there are some constants $0 < C_1 \leq C_2 < \infty$
and a unique $\sigma$-invariant probability measure
$\mu_{ \varphi} $ on the Borel sets of $\Sigma$ such that
\begin{equation}
    \label{gibbszerop}
  C_1 \leq \frac{\mu_{ \varphi}(\Gamma(x_1 x_2 \ldots x_m ))}
      {\exp \{S_m  \varphi(j) -m P(\varphi) \}} \leq C_2 
\end{equation}
for every $j$ in the cylinder set $ \Gamma(x_1 x_2 \ldots x_m).$  

\begin{proposition}
   \label{pstatprocess}
For $\bfl \in \cR_2 \cap \cK,$
and every $n \geq 1,$ the measure $\mu_n$ 
is absolutely continuous with respect to $\mu_{\varphi}$ and
$\mu_n \stackrel{\cD} {\rightarrow} \mu_{ \varphi}$
as $n \rightarrow \infty.$ Furthermore, $P (\varphi) =0.$ 
\end{proposition}

\begin{proof}
Our proof is much the same as the one for Theorem 5.4, \cite{lal0}.
For each $m \geq 1,$ the sum over all $x_1 x_2 \ldots x_m$ in $\cG_m$
of the probabilities  
$P \{ \omega_j = x_j \mbox{ for each } 1 \leq j \leq m \}$
is \mbox{equal $1.$}  Since $C$ in (\ref{stataddfct})
is positive and finite, by (\ref{gibbszerop}), $ \mu_{ \varphi}$ 
is a probability measure and the pressure
$P (\varphi) =0.$ Since 
(\ref{stataddfct}) and (\ref{gibbszerop}) hold for all cylinder
sets and these generate the Borel $\sigma$-algebra, it follows
that $\mu_1 \ll \mu_{ \varphi}$ and that the Radon-Nikodym
derivative $h = ( d \mu_1 / d \mu_{ \varphi}) $ is bounded
away from $0$ and $\infty.$ Next consider the 
restrictions of $\mu_1$ and $\mu_{ \varphi}$ 
to the $\sigma$-algebra  $\cF_n$ generated by the coordinate
functions $x_j,$ $ j \geq n.$ Since the tail field $\cF_{\infty} =
\bigcap_{n \geq 1} \cF_n$ is $0-1$ under $\mu_{\varphi}$
( $\mu_{\varphi}$ is mixing, see e.g. \cite{bowe}),
we obtain
$$
  \left( \frac{d \mu_1 \vert \cF_n}{d \mu_{\varphi} \vert \cF_n}
     \right) 
     = E_{\mu_{\varphi}}( h \vert \cF_n) \longrightarrow
          E_{\mu_{ \varphi}} (h \vert \cF_{\infty}) =1.
$$          
Hence, it follows that $\mu_n \ll \mu_{ \varphi}$ and 
$\mu_n \stackrel{\cD} {\rightarrow} \mu_{ \varphi}.$
\end{proof}

Thus, the stochastic process $\omega_1, \omega_2, \ldots $
is asymptotically stationary, that is, the joint distribution
of $\omega_n , \omega_{n+1}, \ldots$ converges to that one
of a stationary process as $n \rightarrow \infty.$ The limiting
process is a Gibbs state, hence, isomorphic to a Bernoulli shift,
but in general non-Markovian.
From (\ref{stataddfct}), we see that each stationary distribution
decays exponentially in the distance $m$ from the root. The rate
at which the measure of a cylinder set $\Gamma(x_1 x_2 \ldots x_m )$
decays depends on the relative frequency of the generators 
in the reduced word $ x_1 x_2 \ldots x_m $ and the order of the letters.
It is easy to show that the distribution is {\em spherically symmetric}
if and  only if the contact process is isotropic.                     


\subsection{Proof that $\eta < 1$ away from $\cK^c$}

The crux of this section consists of finding bounds for the collection
of ``reinfection" probabilities. In this section, we restrict 
ourselves to $\bfl \in \mbox{int}(\cK_0).$
Let $R_{x_{k+1} x_{k+2} \ldots x_{k+m}}$
denote the event that there is an infection trail that begins 
at $x_1 x_2 \ldots x_{k}$ at some time, reaches vertex
$ x_1 x_2 \ldots x_{k+m} $ at some later time, 
and returns to $x_1 x_2 \ldots x_k$ in finite time.
Thus, $R_{x_1 x_2 \ldots x_m}$
denotes the event that there is an infection trail that begins 
at the root at time $0,$ reaches vertex $ x_1 x_2 \ldots x_m $
at some time, and ever returns to the root from $ x_1 x_2 \ldots x_m.$
First, we note the following about the ``infection probabilities"
$b_j.$ Thanks to constant updating while the infection moves off 
to the boundary of the tree, the $b_j$ in the limit
approximate the actual infection probabilities to any desired degree.
Specifically, for large $k$ and fixed segment $x_1 x_2 \ldots x_{k-1}$
of some infinite word in $\Sigma,$
the probability that $x_1 x_2 \ldots x_{k-1} \ldots x_n$ is ever
infected, given $x_1 x_2 \ldots x_{k-1}$ is infected, is
approximated by $\prod_{j=k}^n b_{x_j}$ as $n \rightarrow \infty,$
where each $b_{x_j}$ depends on the initial infinite string
$\ldots x_1 x_2 \ldots x_{k-1}.$ 
The products  $\prod_{j=k}^n b_{x_j} $ approximate the 
conditional probabilities 
$u_{x_1 x_2 \ldots x_n} / u_{x_1 x_2 \ldots x_{k-1}}$
in the sense that, for large $k,$
the sum on the righthand side of (\ref{hmatrices}) 
is the same as the very same sum with  $\prod_{j=k}^n b_{x_j} $
replacing $u_{x_1 x_2 \ldots x_n} / u_{x_1 x_2 \ldots x_{k-1}}.$
In addition, for two distinct vertices $x= \ldots x_1 x_2 \ldots$ and 
$y = \ldots y_1 y_2 \ldots$ in $ \Sigma$ 
and for all sufficiently large $m,$ the probability of the event 
that $x_1 x_2 \ldots x_m$ or $ y_1 y_2 \ldots y_m $ 
is ever infected is approximately equal 
$\prod_{j=1}^m  b_{x_j} +  \prod_{j=1}^m b_{y_j},$ 
where the $b_{x_j}$ depend on $x \in \Sigma$ and the $b_{y_j}$
depend on $y \in \Sigma.$ 

If we recall the symmetry assumption of the contact process,
(\ref{hmatrices}) for $\rho=2,$ definition (\ref{equationsforbj})
of the $b_j,$ 
and the facts that the $b_j$ are based on the information
contained along an infinite path, whereas we view a finite piece
of the path, a similar reasoning as above applies to the 
probabilities $  P \{ R_{x_{k+1} x_{k+2} \ldots x_{k+m}} \}.$ 
Thus, the probabilities $  P \{ R_{x_{k+1} x_{k+2} \ldots x_{k+m}} \}$ 
of reinfection are approximated by the products 
$\prod_{j=k+1}^{k+m} b^2_{x_j}$ (depending on $x_1 \ldots x_k$)
in such a way that the entire collection of probabilities 
$  P \{ R_{x_1 x_2 \ldots x_m} \} $ is approximated 
to any desired degree in the limit. Hence,
from (\ref{convestimate}), for each $\bfl $ such that $r_u(\bfl)<2,$ 
for every $x = \ldots x_1 x_2  \ldots \in \Sigma,$
and all  integers $m,k  >0,$ 
\begin{equation}
   \label{retprobandpotential}
  P \{  \bigcup_{x_{k+1}  \ldots x_{k+m} \in \cG_m}  
     R_{x_{k+1}  \ldots x_{k+m}} \} = 
             C \exp \{ S_m \varphi_2(x) \} (1+ O(\alpha^k))
\end{equation}
for some constants $0 < \alpha < 1$ and  $0 < C < \infty.$ 
We point out that an alternative approach (not pursued here), 
which leads to the same estimates, is based on 
approximating the set of vertices ever
to be infected by labelled Galton-Watson trees. 

\begin{proposition}
   \label{returnspace}
For each $\bfl \in \mbox{int}(\cR_2 \cap \cK)$ such that $r_u(\bfl)<2,$
there exist constants
$ 0 < C_1, C  < \infty$ such that
for every $\epsilon > 0,$ sufficiently small $\zeta >0,$ 
and fixed sufficiently large $k,$
\begin{equation}
   \label{returnprbound}
 C_1  \exp \{ n_{\zeta} 
                 (\overline{\varphi}_2 - \epsilon)\}   
    \leq  P \{ \bigcup_{m=k}^{\infty} \bigcup_{x_1  \ldots x_m \in \cG_m}  
       R_{x_1 \ldots x_m} \} \leq 
    C  \epsilon n_{\zeta} \,
                 \exp \{ n_{\zeta} 
                 (\overline{\varphi}_2 + \epsilon)\}, 
\end{equation}
where $  n_{\zeta}=\log \zeta / \overline{\varphi}_{2}.$
\end{proposition}

\begin{proof}
First observe that, by subadditivity, 
  $ P \{ R_{x_1 \ldots x_m} \} \leq  
    P \{ R_{x_{k+1} x_{k+2} \ldots x_{k+m}} \}$ 
(because $u_{x_1 \ldots x_m} \leq u_{x_1 \ldots x_{k+m}} / u_{x_1 
\ldots x_k}$). Fix such $\bfl,$
some $\epsilon > 0,$ and some sufficiently small $\zeta > 0.$
Recall that every generic $x \in \cA_{2}^*(\zeta)$ satisfies 
\begin{equation}
  \label{typicalpot}
      \exp \{ n_{\zeta} 
    (\overline{\varphi}_2 - \epsilon)\} \leq
   \exp \{ S_{n_{\zeta}} \varphi_2(x) \} \leq \exp \{ n_{\zeta} 
    (\overline{\varphi}_2 + \epsilon)\}. 
\end{equation} 
Next fix some sufficiently large integer $k.$ 
Then by the previous remarks, by Proposition \ref{arhocard}, by
(\ref{birkhoff2}),  
 (\ref{retprobandpotential}), and  (\ref{typicalpot}), we obtain 
\begin{eqnarray*}
  P \{ \bigcup_{m=k}^{\infty} \bigcup_{x_1  \ldots x_m \in \cG_m} 
       R_{x_1 \ldots x_m} \}  
     & \leq & 
         P \{ \bigcup_{m=1}^{\infty} 
            \bigcup_{x_{k+1} \ldots x_{k+m} \in \cG_m} 
      R_{x_{k+1}  \ldots x_{k+m}} \} \\*[0.1cm]
        & \leq & 
       C_2 \, \sum_{m=1}^{\infty} \; \, 
        \sum_{\mbox{one } x_{1} \ldots x_{m} \in  
                   \cA_{2}^*(\zeta)} 
            \, \exp \{ S_m \varphi_2(x) \}  \\*[0.1cm]
        & \leq &  
            C_2 c \, \sum_{m:
      \vert 1 - m/n_{\zeta} \vert \leq \epsilon} \;  \;
              \sum_{\mbox{one } x_1 x_2 \ldots x_m \in  
                   \cA_{2}^*(\zeta)} 
               \, \exp \{ S_m \varphi_2(x) \}  \\*[0.1cm]
         & \leq &  C_2 c \, (2 \epsilon n_{\zeta}) \, 
                 \exp \{ n_{\zeta} 
                 (\overline{\varphi}_2 + \epsilon)\} 
\end{eqnarray*}  
for some constants $0 < C_2, c < \infty.$

The reverse direction is straightforward. Indeed, 
by Proposition \ref{arhocard}, by  
(\ref{retprobandpotential}) and  (\ref{typicalpot}),
 
\begin{eqnarray*}
  P \{ \bigcup_{m=k}^{\infty} \bigcup_{x_1  \ldots x_m \in \cG_m} 
      R_{x_1 \ldots x_m} \}  
     & \geq &  P \{ \bigcup_{x_1 x_2 \ldots x_m \in  
                   \cA_{2}^*(\zeta)} 
                     R_{x_{m+1} \ldots x_{2m}} \}  \\*[0.1cm]
              & \geq & C_1 \, \, \sum_{\mbox{one } x_1 x_2 \ldots x_m \in  
                   \cA_{2}^*(\zeta)} 
                   \, \exp \{ S_m \varphi_2(x) \}  \\*[0.1cm] 
          & \geq &  C_1 \, \exp \{ n_{\zeta} 
                 (\overline{\varphi}_2 - \epsilon)\} 
\end{eqnarray*}  
for some constant $0 < C_1 < \infty.$
This completes the proof of (\ref{returnprbound}). 
\end{proof}

{\bf Time-dependent Infection Probabilities.}
Next, we turn to discuss the time-dependent infection 
probabilities
\begin{equation}
   \label{timedepinfprob}
   u_{x,t}(\bfl) = P  \{ x \in A_t \}
\end{equation}
for $x \in \cG$ and every real $t >0.$
Again, the strong Markov property and 
the monotonicity and homogeneity properties of the process 
imply that $ u_{xy,s+t}(\bfl) \geq u_{x,s}(\bfl) u_{y,t}(\bfl)$
for all $s,t >0$ and $x,y \in \cG$ such that 
$\vert x  y \vert = \vert x \vert + \vert y \vert.$
Clearly, $u_{x,t}(\bfl) \leq u_x(\bfl).$ 
A subadditivity argument shows that, for every $x \in \cG_k$
and $y_n = x x\ldots x \in \cG_{nk},$ the limit 
\begin{eqnarray}
  \label{ulimit}
   \lim_{n \rightarrow \infty}
   u_{y_n,kns}(\bfl)^{1/n} & = & U_{x,s} =  U_{x,s}(\bfl) 
\end{eqnarray}
exists, and that $u_{y_n, kns}(\bfl) \leq U_{x,s}(\bfl)^n$ 
for all integers $n \geq 0,$ in particular, for $x=i \in \cA,$ we have 
$u_{y_n, ns}(\bfl) \leq U_{i,s}(\bfl)^n.$ Obviously, since 
$u_{x,t}(\bfl) \leq u_x(\bfl),$ we have  
\begin{equation}
    \label{compubeta}
    U_{x,s}(\bfl) \leq \beta_x (\bfl)
\end{equation}
for every $x \in \cG$ and each $s > 0.$ 
Observe that the functions $U_{x,s}(\cdot)$ share the same basic
properties with the functions $u_x(\cdot),$ which were essential
to prove their geometric decay in the distance from the root,
for instance,
the H\"{o}lder conditions, and the convergence of their sums. 
Therefore, a parallel analysis may be carried out to define
potential functions and the corresponding Gibbs states as described
for the $u_x(\bfl).$ To save space, we
shall omit the details with the exception of 
a key ingredient to the proof of $\eta < 1,$
which is the scaling properties between the
space variable $x \in \cG$ and the time variable
$s$ in $U_{x,s}(\bfl).$
To wit, let $\Phi_{\rho,s}(x) = \Phi_{\rho,s; \bfl}(x) $
be the analogue of $\varphi_{\rho}(x) = \varphi_{\rho;\bfl}(x),$
where in definition (\ref{hmatrices}), on the righthand side,
the ratio $ (u_{x_1 x_2 \ldots x_n}(\bfl)/
                  u_{x_1 x_2 \ldots x_{k-1}}(\bfl) )^{\rho}$ is
 replaced by  $(u_{x_1 x_2 \ldots x_n,ns}(\bfl) /
                  u_{x_1 x_2 \ldots x_{k-1}}(\bfl) )^{\rho},$ 
that is,
\begin{equation}
   \label{htimematrices}
    (H_{\rho,s} (n ; x_1 x_2 \ldots x_{k-1}; \bfl))_{ij} = 
         \sum_{ x_{k+1} \ldots x_{n-1} \in \cF_{ij}}
               \, \, \left(\frac{u_{x_1 x_2 \ldots x_n,ns}(\bfl)}{
                  u_{x_1 x_2 \ldots x_{k-1}}(\bfl)} \right)^{\rho} 
\end{equation}                       
for every $i,j \in \cA.$ 
Then the same machinery leads to the potential function  
$\Phi_{\rho,s}$ and its average 
with respect to the probability measure $\mu_{\delta \Phi_{\rho,s}},$ 
\begin{equation}
   \label{bigphimeans}
 \overline{\Phi}_{\rho,s} = \int \Phi_{\rho,s} \, 
                           d \mu_{\delta \Phi_{\rho,s}}.
\end{equation}
 
For every $s>0,$ let $R_{x_1 x_2 \ldots x_m}(ms)$
denote the event that there is an infection trail that begins 
at the root at time $0,$ reaches vertex $ x_1 x_2 \ldots x_m $
at some time, returns to the root from  $ x_1 x_2 \ldots x_m $ 
prior to time $ms$ and stays
without recovery mark * up to time $ms.$ By
the same token as (\ref{returnprbound}) was derived, 
for every $\epsilon > 0,$  sufficiently small $\zeta >0,$ 
and fixed sufficiently large $k,$ 
we obtain
\begin{equation}
   \label{returntime} 
  c_3  \exp  \{ n_{\zeta}( \overline{\Phi}_{2,s} - \epsilon) \}    
       \leq   P \{ \bigcup_{m=k}^{\infty} 
         \bigcup_{x_1  \ldots x_m \in \cG_m}  
       R_{x_1 \ldots x_m}(ms) \}  
         \leq c_4 \epsilon n_{\zeta} \, 
         \exp  \{ n_{\zeta}( \overline{\Phi}_{2,s} + \epsilon) \}
\end{equation}
for some constants $0 < c_3, c_4 < \infty$ and
where  $ n_{\zeta} = n_{\zeta,s} =
 \log \zeta/ \overline{\Phi}_{2,s}.$ Therefore by definition, for
every $K>0,$ 
$$
  c_3  \exp  \{ n_{\zeta}( \overline{\Phi}_{2,Ks} - \epsilon) \}    
       \leq  P \{ 
        \bigcup_{m=k}^{\infty} \bigcup_{x_1  \ldots x_m \in \cG_m}  
       R_{x_1 \ldots x_m}(msK) \}  
         \leq c_4 \epsilon n_{\zeta} \, 
             \exp  \{ n_{\zeta}( \overline{\Phi}_{2,Ks} + \epsilon) \}.
$$
Let $t_{\zeta}$ denote the time for which the 
expected number of infection trails returning to the root
from distance $n_{\zeta}$ from the root and staying
without recovery mark * up to time $t_{\zeta}$ is {\em maximal}. 
Then in view of $ n_{\zeta^K} = K n_{\zeta}$ and 
by the additivity of the functional $S_n \Phi_{2,s},$ 
$$
  c_3  \exp  \{ n_{\zeta^K}( \overline{\Phi}_{2,t_{\zeta}} - \epsilon) \}    
       \leq  P \{ 
         \bigcup_{m=k}^{\infty} \bigcup_{x_1  \ldots x_m \in \cG_m}  
       R_{x_1 \ldots x_m}(m t_{\zeta} K) \}  
         \leq c_4 \epsilon n_{\zeta} \, 
         \exp  \{ n_{\zeta^K}( \overline{\Phi}_{2,t_{\zeta}}
         + \epsilon) \}.
$$
Combining the last two displays brings 
$ K \overline{\Phi}_{2,t_{\zeta}} =  \overline{\Phi}_{2,K t_{\zeta}}.$ 
In other words, distance and time are scaled by the same factor
in the mean functionals $ \overline{\Phi}_{2,t_{\zeta}},$ 
that is, there is some constant $0 < \Delta(\bfl) < \infty$ such that,
for each $\zeta > 0,$ 
\begin{equation}
   \label{linearscale} 
      t_{\zeta} = \Delta(\bfl) n_{\zeta}.
\end{equation}
Observe that, due to the form of (\ref{genlength}), for each $s>0,$
\begin{equation}
  \label{twoaverages}
   n_{\zeta,s} \overline{\Phi}_{2,s} = n_{\zeta} \overline{\varphi}_2
   = \log \zeta.
\end{equation}

\begin{proposition}
   \label{etalessthanone}
For each $\bfl$ in the interior of $\cK \cap \cR_2$ 
with $r_u(\bfl) < 2,$ there is some  \linebreak $0 < \Delta = 
\Delta(\bfl) < \infty$ such that
$\eta^{\Delta} = \eta(\bfl)^{\Delta} = 
\exp \overline{\varphi}_{2; \bfl},$ 
in particular,
$$ 
     \exp \overline{\varphi}_{2; \bfl} <  1  \,
     \mbox{ if and only if } \,  \eta(\bfl) < 1.
$$
Furthermore, for each $\bfl$ in the interior of $\cK$
with $r_u(\bfl) < 2,$ we have $\eta(\bfl) < 1.$
\end{proposition} 

\begin{proof}
First recall that for  $\bfl$ in the interior 
$\mbox{int}(\cR_2 \cap \cK)$ with $r_u(\bfl) < 2,$ 
we have $\overline{\varphi}_2 < 0.$  
Since the function $\eta(\cdot)$ is
nondecreasing in each variable $\lambda_j,$
it remains to be shown that 
$\eta(\bfl)^{\Delta} = \exp \overline{\varphi}_{2; \bfl}$ 
for $ \bfl \in \mbox{int}(\cR_2 \cap \cK).$ 

If we fix some $\epsilon >0,$ some sufficiently small $\zeta >0,
$ and some sufficiently large $k,$ by (\ref{linearscale})
and Proposition \ref{returnspace}, 
\begin{eqnarray*}
    \eta = \lim_{t \rightarrow \infty} (P \{ 1 \in A_t \} )^{1/t} 
         & = &  \lim_{t \rightarrow \infty}
              \left(  P \{ \bigcup_{m=k}^{\infty} 
                \bigcup_{x_1 x_2 \ldots x_m \in \cG_m}  
                 R_{x_1 \ldots x_m}(t) \} \right)^{1/t}  \\*[0.1cm]
             & \leq &  \lim_{t_{\zeta} \rightarrow \infty} 
                  [ (C \epsilon n_{\zeta})
                 \exp \{ n_{\zeta} 
                 (  \overline{\varphi}_{2} + \epsilon)\}
                  ]^{1/t_{\zeta}} \\*[0.1cm]
            &  = &  \lim_{n_{\zeta} \rightarrow \infty} 
                  [ (C \epsilon n_{\zeta}) 
                 \exp \{ n_{\zeta} 
                 (  \overline{\varphi}_{2} + \epsilon)\}
                  ]^{1/\Delta n_{\zeta}} \\*[0.1cm]
          & = &   \exp  \{  \overline{\varphi}_{2} + 
                       \epsilon \} ^{1/ \Delta}    \\
                & = &  \exp  \{ ( \overline{\varphi}_{2} + 
                       \epsilon)/ \Delta \},           
\end{eqnarray*} 
where the employed constant $C$ is positive and finite. 

Next, it is not difficult to see that 
Proposition \ref{returnspace} has a time-dependent
version (the proof being the same), namely,
$$
C_1  \exp \{ n_{\zeta,s}  ( \overline{\Phi}_{2,s} 
                  - \epsilon) \}   
    \leq   P \{ \bigcup_{m=k}^{\infty} 
         \bigcup_{x_1 \ldots x_m \in \cG_m}  
         R_{x_1 \ldots x_m}(ms) \} \leq 
     C  \epsilon n_{\zeta,s} \,
                 \exp \{ n_{\zeta,s} 
                 ( \overline{\Phi}_{2,s} + \epsilon)\}, 
$$
for every $s>0.$ 
Then the lower inequality follows from similar arguments in combination
with the product structure of the probabilities 
$\prod_{j=k+1}^{k+m} b_{x_j}^2,$
on which the approximation (\ref{returnprbound}) is based, and again
the equality $ n_{\zeta,s} \overline{\Phi}_{2,s} = 
n_{\zeta} \overline{\varphi}_2,$ that is,
\begin{eqnarray*}
    \eta = \lim_{t \rightarrow \infty} (P \{ 1 \in A_t \} )^{1/t} 
         & = &  \lim_{t \rightarrow \infty}
              \left(  P \{ \bigcup_{m=k}^{\infty} 
                \bigcup_{x_1  \ldots x_m \in \cG_m}  
                 R_{x_1 \ldots x_m}(t) \} \right)^{1/t}  \\[0.1cm]
             & \geq &  \lim_{t_{\zeta} \rightarrow \infty} 
                  [  C_1  \exp \{ n_{\zeta,t_{\zeta}} 
                   ( \overline{\Phi}_{2,t_{\zeta}} 
                  - \epsilon) \} 
                  ]^{1/t_{\zeta}} \\[0.1cm]
               & = &  \lim_{t_{\zeta} \rightarrow \infty} 
                  [  C_1  \exp \{ n_{\zeta}  ( \overline{\varphi}_{2} 
                  - \epsilon \gamma) \} 
                  ]^{1/t_{\zeta}} \\[0.1cm]
              & = &  \lim_{n_{\zeta} \rightarrow \infty} 
                  [  C_1  \exp \{ n_{\zeta}  ( \overline{\varphi}_{2} 
                  - \epsilon \gamma) \} 
                  ]^{1/\Delta n_{\zeta}} \\[0.1cm]
            & = &  \exp  \{ ( \overline{\varphi}_{2} - 
                       \epsilon \gamma)/ \Delta \},     
\end{eqnarray*} 
where $\gamma =  \overline{\varphi}_{2} /  
\overline{\Phi}_{2,t_{\zeta}},$
the  constant $C_1$ is positive and finite and $k$ 
was chosen suitably large. Since $\epsilon >0$ and $\zeta >0$
were both arbitrary and $\overline{\varphi}_2 < 0,$
the inequality $\eta < 1$ follows together with all other claims.
This finishes our proof.
\end{proof}


\section{First-Passage and Backscatter Matrices} 
\setcounter{equation}{0}

Recall from (\ref{bmatrix}) that, for fixed $x, \rho,$ and $\bfl,$  
$B_{\rho}(\bfl)$ is a matrix with entries $b_j(\bfl)^{\rho},$
where the $b_j(\bfl)^{2}$ are associated with $\varphi_{2 ; \bfl}.$
Now, if we let $\overb_j^{2} = \overb_j(\bfl)^{2}$ be the matrix
entries associated with the function 
$\overline{\varphi}_{2 ; \bfl},$ 
where $\overline{\varphi}_{2 ; \bfl}$ was defined in  
(\ref{phimeans}) for $\rho=2,$ then, 
for any positive number $\rho,$ define the matrix 
$M_{\rho} = M_{\rho} (\bfl)$ to be the $2d 
\times 2d$ matrix, indexed by elements of $\cA,$ whose entries 
are given by 
\begin{eqnarray}
  \label{rhomatrix}
            (M_{\rho} (\bfl))_{ij} & = & \overb_j(\bfl)^{\rho} \quad
                        \mbox{if }  j \not = i^{-1}, \\
                            & = & 0 \qquad \, \, \quad  
                        \mbox{if }  j = i^{-1}. \nonumber
\end{eqnarray}
We will refer to $M_1(\bfl)$ and $M_2(\bfl)$ as the
{\em first-passage matrix} and the {\em backscatter matrix},
respectively. These two matrices will play a distinguished role,
for instance, they are important in the proofs of Theorems 
\ref{limitset} and \ref{halfboundary}, stated in the Introduction.
Note that $M_2$ coincides with the matrix in  
(\ref{backscattermatrix}).
Since by construction, for $\rho > 0 $ and $ \bfl \in \cK_0,$
$M_{\rho}(\bfl) $ is an aperiodic, irreducible, and
nonnegative matrix, thus, a Perron-Frobenius matrix,
the Perron-Frobenius theorem lets us conclude that 
$M_{\rho}(\bfl)$ has a largest positive
eigenvalue $ \theta(\rho; \bfl). $ Henceforth, we will recur to the
shorthands $\theta_{\rho}$ for $ \theta_{\rho}(\bfl) = \theta(\rho;\bfl), $
in particular, we write
\begin{eqnarray}
 \label{eigenvaluestheta}
     \theta_1 & = & \theta_1(\bfl) = \theta(1;\bfl) = \theta(\bfl) \\ 
      \theta_2 & = &  \theta_2(\bfl) = \theta(2; \bfl). \nonumber
\end{eqnarray}
Note that $\theta_2 = \exp \overline{\varphi}_2 $ and that
the eigenvalue $\theta$ in Theorem \ref{limitset} coincides
with $\theta_1.$ In addition, it is not difficult to 
show that $\theta_1 = \exp \overline{\varphi}_1 $ 
for $1 > r_u(\bfl)$ and
$\theta_{\rho} = \exp \overline{\varphi}_{\rho}$ 
for $\rho > r_u(\bfl).$ 
By definitions (\ref{bmatrix}) and (\ref{rhomatrix}) 
and by Proposition \ref{leadeva}, the entries of $M_{\rho}$
satisfy the equation           
\begin{equation}
  \label{mequation}
        \sum_{i \in \cA} \, \frac{\overb_i^{\rho}}{
                    \theta_{\rho} + \overb_i^{\rho}} = 1.
\end{equation}


\subsection{Strict Monotonicity of the Lead Eigenvalues}

\begin{proposition}
  \label{potmonotone}
Let $\bfl \in \cK_0$ and let $\tilde{\bfl}$ be an infection parameter
in the direction of decrease for $\bfl,$ as defined in Section 2.5.
Then for $\rho > r_u(\bfl)$ and every $s>0,$  
\begin{eqnarray}
         \label{fmeansmono}
       \overline{\varphi}_{\rho; \tilde{\bfl}} 
                   & < &  
        \overline{\varphi}_{\rho; \bfl},  \\ 
          \label{smeansmono} 
      \overline{\Phi}_{\rho,s; \tilde{\bfl}} & < & 
            \overline{\Phi}_{\rho, s; \bfl}  .
\end{eqnarray}
In particular, 
\begin{eqnarray}
   \label{eigenvmono}
      \theta_{1}(\tilde{\bfl}) & < & \theta_{1}(\bfl), \\
       \theta_{2}(\tilde{\bfl}) & < & \theta_{2}(\bfl). \nonumber        
 \end{eqnarray}             
Moreover for each $s>0,$ each $\bfl \in \cK_0,$ and 
$t > \rho > r_u(\bfl),$ we have  
\begin{eqnarray}
   \label{tsmono}
     \overline{\varphi}_{t; \bfl} & < & 
       \overline{\varphi}_{\rho; \bfl} , \\
  \overline{\Phi}_{t,s; \bfl} & < &  \overline{\Phi}_{\rho,s; \bfl},
      \nonumber
\end{eqnarray}
thus, $\theta_t(\bfl) < \theta_{\rho}(\bfl),$ in particular,
$ \theta_2(\bfl) < \theta_1(\bfl).$ 
\end{proposition}
  
\begin{proof}
Only the proof of (\ref{fmeansmono}) is carried out. The proof of
(\ref{smeansmono}) runs in parallel whereas the second line of
(\ref{eigenvmono}) will then follow immediately from the
definition of $\theta_2$ and the first line of (\ref{eigenvmono}) 
will follow by recalling the definitions of the $\overb_j$
and $\theta_1.$ The verification of (\ref{tsmono}) is deferred 
to the end of this proof. 
Fix $\bfl \in \cK_0$ and $\rho > r_u(\bfl).$
Write
\begin{equation}
  \label{littleh}
    h_{\rho}(\bfl) = h_{\rho}(\bfl,x,k) 
      = \left \{ \sum_{x_{k+1} x_{k+2} \ldots x_{n-1} \in \cG_{n-k-1}}
         \, \left(\frac{u_{x_1 x_2 \ldots x_n}(\bfl)}{
                  u_{x_1 x_2 \ldots x_{k-1}}(\bfl)} \right)^{\rho} 
              \right \}^{1/(n-k)} 
\end{equation} 
for all integers $n-1 > k>0,$  $x \in \Sigma,$
and $\bfl \in \cK_0.$ Let $\tilde{\bfl}$ be an 
infection parameter in the direction of decrease for $\bfl.$
From (\ref{hpotential}),
we see that there is some $0 < C(x,\bfl) < \infty$
so that for all sufficiently large $n,$ 
$$ 
   W(x )^t 
     H_{\rho}(n; x_1 x_2 \ldots x_{k-1}; \bfl) 
      V(\sigma^{n-k} x ) = C(x,\bfl)\, 
      \exp \{ S_{n-k} \varphi_{\rho; \bfl}(x)\}.
$$
It suffices to prove that, for fixed $x \in \Sigma,$
\begin{equation}
   \label{hlessthanh}
   \lim_{n,k  \rightarrow \infty}  h_{\rho}(\tilde{\bfl},x,k) 
       <   \lim_{ n,k  \rightarrow \infty}  h_{\rho}(\bfl,x,k)  
\end{equation}
(where the limits exist in view of our assumption $\rho >r_u$).
The meaning of the limit $  \lim_{ n,k  \rightarrow \infty}$
is that we take both variables to $\infty,$ yet $n$ much faster than
$k.$ 
Once (\ref{hlessthanh}) has been shown to hold, 
it will follow that, in particular, it is valid 
for all $x$ such that $  S_{n} \varphi_{\rho; \bfl}(x)/n \rightarrow
\overline{\varphi}_{\rho; \bfl}$ 
as $n\rightarrow \infty,$ consequently, 
 $  \overline{\varphi}_{\rho; \tilde{\bfl}}
              <  \overline{\varphi}_{\rho; \bfl}.$ 

Fix $x = \ldots x_1 x_2 \ldots \in \Sigma.$ 
From (\ref{udecrease}), recall
that there is some constant $0 < \omega < 1$ such that
     $ u_{x_1 \ldots x_n} (\tilde{\bfl}) \leq
      u_{x_1 \ldots x_n}(\bfl) \omega^n.$
Choose $n$ sufficiently large so that 
$$
   \frac{ u_{x_1 x_2 \ldots x_{k-1}}(\bfl)}
        { u_{x_1 x_2 \ldots x_{k-1}}(\tilde{\bfl})} 
          < \frac{1}{\omega^{n/2}}.
$$           
Now combine these to 
\begin{eqnarray*}
   \frac{ u_{x_1 x_2 \ldots x_n}(\tilde{\bfl})}
        { u_{x_1 x_2 \ldots x_{k-1}}(\tilde{\bfl})}
         & = &   \frac{ u_{x_1 x_2 \ldots x_n}(\bfl)}
         { u_{x_1 x_2 \ldots x_{k-1}}(\bfl)} \cdot
           \frac{ u_{x_1 x_2 \ldots x_n}(\tilde{\bfl})}
        { u_{x_1 x_2 \ldots x_n}(\bfl)} \cdot
          \frac{ u_{x_1 x_2 \ldots x_{k-1}}(\bfl)}
        { u_{x_1 x_2 \ldots x_{k-1}}(\tilde{\bfl})} \\
       & < & \frac{ u_{x_1 x_2 \ldots x_n}(\bfl)}
         { u_{x_1 x_2 \ldots x_{k-1}}(\bfl)} \,  
          \omega^n \, \omega^{-n/2} \\ 
      & = & \frac{ u_{x_1 x_2 \ldots x_n}(\bfl)}
         { u_{x_1 x_2 \ldots x_{k-1}}(\bfl)} \, \omega^{n/2} . 
 \end{eqnarray*} 
Hence,
\begin{eqnarray*}
    h_{\rho}(\tilde{\bfl},x,k) & < &  \omega^{\rho n /(2 (n-k))}
                  \,  h_{\rho}(\bfl,x,k),    \\
\mbox{and} \qquad \qquad \qquad \qquad \qquad & & \\
    \lim_{n,k  \rightarrow \infty}  h_{\rho}(\tilde{\bfl},x,k) 
       & < &  \lim_{n,k \rightarrow \infty}  
             \omega^{\rho n /(2 (n-k))} h_{\rho}(\bfl,x,k) \\
       & < & \omega^{\rho /2}  \,                       
             \lim_{n,k \rightarrow \infty} h_{\rho}(\bfl,x,k) \\
       & < &  \lim_{n,k \rightarrow \infty} h_{\rho}(\bfl,x,k),        
\end{eqnarray*}
which proves (\ref{fmeansmono}) since $\rho > r_u(\bfl)$ was arbitrary.
Consequently, Proposition \ref{leadeva} implies that the function 
$b_j(\bfl), $ $ j \in \cA,$ may be chosen
as strictly decreasing function along directions of
decrease, and similarly, the entries of $M_{\rho} (\bfl).$ 
Therefore, we also verified the two inequalities in (\ref{eigenvmono})
in view of the definitions of $\theta_1$ and $\theta_2.$

Finally, the first statement in (\ref{tsmono}) (and similarly, 
its second statement) rests on the following observation.
By arguments along the lines used in the proof of 
Proposition \ref{hoelderquot}, there is a constant
$ 0 < \tilde{\gamma} <1$ so that for each $k>0$ and
 sufficiently large $n,$ we have
 $u_{x_1 \ldots x_n}(\bfl) / u_{x_1 \ldots x_{k-1}}(\bfl)
   \leq \tilde{\gamma}^n.$
Hence,
\begin{eqnarray*}
    (\frac{u_{x_1 x_2 \ldots x_n}(\bfl)}
         { u_{x_1 x_2 \ldots x_{k-1}}(\bfl)})^{t}
         & = & (\frac{u_{x_1 x_2 \ldots x_n}(\bfl)}
         { u_{x_1 x_2 \ldots x_{k-1}}(\bfl)})^{t- \rho} \, 
               (\frac{u_{x_1 x_2 \ldots x_n}(\bfl)}
         { u_{x_1 x_2 \ldots x_{k-1}}(\bfl)})^{\rho} \\
          & < &
         \gamma^n   (\frac{u_{x_1 x_2 \ldots x_n}(\bfl)}
         { u_{x_1 x_2 \ldots x_{k-1}}(\bfl)})^{\rho}
\end{eqnarray*}
for some constant $ 0<\gamma <1$ because $t > \rho.$
Proceeding then along the above route to prove (\ref{fmeansmono}) 
provides the desired result and ends our proof.
\end{proof}

Observe that this proof also shows that the functions 
in (\ref{fmeansmono})--(\ref{eigenvmono}) are nondecreasing 
in each variable $\lambda_j$ because, clearly, for all 
$\tilde{\bfl}$ and $\bfl$  in $\cK_0$ with all infection rates
$\tilde{\lambda}_i = \lambda_i$ the same except for
one $\tilde{\lambda}_j < \lambda_j,$ 
for sufficiently large $n,$ the ratio
 $( u_{x_1  \ldots x_n}(\tilde{\bfl}) /
        u_{x_1 \ldots x_{k-1}}(\tilde{\bfl}) ) \cdot
          ( u_{x_1 \ldots x_{k-1}}(\bfl)/
         u_{x_1  \ldots x_n}(\bfl)) \leq 1.$


\section{Continuity Properties of the Infection Probabilities}
\setcounter{equation}{0}

This section is concerned with the properties of continuity and
discontinuity of the lead eigenvalues of the matrices $M_{\rho}$
and the functions $\beta_x$ and ends with the proofs of the
results about the weak survival phase (Theorems \ref{weaksurvival}
through \ref{etaatboundary}) and of Theorem \ref{properties}.


\subsection{Continuous Potential Functions and Lead Eigenvalues}

If for every $ \epsilon >0$ and  
for every subset $A \subset {\bf R}^d,$
we let $A_{\epsilon}$ be the $\epsilon$-neighbourhood of $A,$
i.e.\ consisting of all those $z \in {\bf R}^d$ that have
(Euclidean) distance to $A$ less than $\epsilon,$ and 
let $\cK^c$ denote the complement of $\cK$ in ${\bf R}^d,$ then
define $\cK_I = \cK \setminus (\cK^c)_{\epsilon}$ to be the set of points
in $\cK$ which are at a distance to $\cK^c$ at least $\epsilon.$
The next result along with Proposition \ref{secondtrans} below are
crucial moments in proving \mbox{Theorem \ref{weaksurvival}.}
 
\begin{proposition}
  \label{contpotential}
For each $\rho > r_u(\bfl),$ $s>0,$ and all $i, j \in \cA,$ 
the functions
\begin{eqnarray}
    \label{threefctcont}
      \lambda_j & \rightarrow &  \overb_i(\bfl)  \\  
       \lambda_j &  \rightarrow & 
      \overline{\varphi}_{\rho}= \overline{\varphi}_{\rho; \bfl}
        \nonumber \nonumber \\
         \lambda_j &  \rightarrow & 
      \overline{\Phi}_{\rho,s}= \overline{\Phi}_{\rho,s; \bfl}
          \nonumber 
\end{eqnarray} 
are continuous for each $\bfl \in \mbox{int}(\cK),$ 
with $ \overb_i(\bfl)$ being defined in (\ref{rhomatrix}). 
The meaning of the first statement is that each $\overb_i(\bfl)$
can be chosen as a continuous function in each variable $\lambda_j.$
In particular, 
\begin{equation} 
  \label{etacont}
   \lambda_j \rightarrow \eta(\bfl)
\end{equation}
is a continuous function for each $\bfl \in \mbox{int}(\cK \cap \cR_2)$ 
with $r_u(\bfl) < 2.$  
\end{proposition}

\begin{proof}
Since the verification is the same for both functions
$\overline{\varphi}_{\rho}$ and $\overline{\Phi}_{\rho,s},$ 
we only present the proof for the former. As explained in the
proof of Proposition \ref{potmonotone}, we can conclude
that the property holds for $ \overb_i(\bfl)$ as well, 
which will accomplish the first line of (\ref{threefctcont}).
Also, by virtue of Proposition \ref{etalessthanone}, 
claim (\ref{etacont}) will follow immediately.

Let $\rho > r_u(\bfl).$ 
Let $ \epsilon >0.$ Observe that for $\bfl \in \cK_I,$ 
each $\lambda_k >0.$ For each $T>0,$ define
\begin{eqnarray*}
  u^T_x & = & u^T_x(\bfl) 
            = P \{ x \in A_t \mbox{ for some } t \in [0,T] \}, \\
  u^{(T, \infty)}_x & = &  u^{(T, \infty)}_x (\bfl) =
         P \{ x \in A_t \mbox{ for some } t > T \}. 
\end{eqnarray*}
Denote  $ h_{\rho}(\bfl) = h_{\rho}(\bfl,x,k)$ as it has already been
defined in (\ref{littleh}) and define
$$
    g_{\rho}(\bfl) = g_{\rho}(\bfl,x,k) = h_{\rho}(\bfl,x,k)^{n-k}.
$$ 
Write $h^T_{\rho}(\bfl,x,k)$ and  $g^T_{\rho}(\bfl,x,k),$ respectively,
for the functions that result when all $u_x$ 
are replaced by $u^T_x$
in the expressions for $h_{\rho}(\bfl,x,k)$ and  
$g_{\rho}(\bfl,x,k),$ respectively. 
With these agreements along with the observation relying on 
(\ref{hpotential}), as described in the proof of 
Proposition \ref{potmonotone}, 
in order to prove that the function $\lambda_j \rightarrow  
\overline{\varphi}_{\rho}= \overline{\varphi}_{\rho; \bfl}$ is 
continuous, it is sufficient to show that, for fixed $x \in \Sigma,$
for $\rho > r_u(\bfl),$ and
each $j \in \cA,$ the function
$$
   \lambda_j \rightarrow  \lim_{n,k \rightarrow \infty}
       h_{\rho}(\bfl,x,k) =  \lim_{n,k \rightarrow \infty} 
                             \lim_{T \rightarrow \infty}
                               h^T_{\rho}(\bfl,x,k)
$$     
is continuous for each $\bfl \in \cK_I.$ 
Here, $\lim_{n,k \rightarrow \infty}$ means that $n$ and $k$
are taken to $\infty,$ $n$ much faster than $k.$
Since $\epsilon >0$ 
is arbitrary and in view of Proposition \ref{leadeva},
our first two statements in (\ref{threefctcont}) will then follow.

Key ingredients to the proof are the following
three items that we will verify below:
\begin{enumerate}
\item[(A)]
 For each $j \in \cA,$  the function
 $ \lambda_j \rightarrow h^T_{\rho}(\bfl,x,k)$ is continuous for
 $ \bfl \in K_I$ and fixed $T.$ 
\item[(B)] 
For fixed $\bfl,$ the convergence of $h^T_{\rho}(\bfl,x,k)$ to
$h_{\rho}(\bfl,x,k),$ as $T \rightarrow \infty,$
can be controlled uniformly for $\bfl \in \cK_I.$
\item[(C)] 
For fixed $\bfl,$ the convergence of $h_{\rho}(\bfl,x,k)$ 
to its limit, as $n, k \rightarrow \infty,$ happens in a uniform
fashion as well for $\bfl \in \cK_I.$ 
\end{enumerate}
These three claims have some consequences. From the first and second claim,
it will follow that the function
$\lambda_j \rightarrow h_{\rho}(\bfl,x,k)$ is continuous. The
third claim will imply that the limit $\lambda_j \rightarrow 
\lim_{n,k \rightarrow \infty} h_{\rho}(\bfl,x,k)$ is a continuous
function for $\bfl \in \cK_I.$   

Indeed, the first claim is straightforward since it is obvious
that $u^T_x(\cdot)$ is a continuous function in each $\lambda_j$
and that, because each $\lambda_k >0,$ we have $u^T_x(\cdot) >0.$ 
Hence, the function $h^T_{\rho}(\bfl,x,k)$ inherits these properties
from the $u^T_x(\cdot).$ 

To show the second claim, note that by monotonicity,
for each $\bfl \in \cK_I,$
\begin{eqnarray}
   \label{gandhstrict}
   h_{\rho}(\bfl,x,k) & \leq &  
        \sup_{\tilde{\bfl} \in \cK \cap (\cK^c)_{\epsilon} } \, 
         \max_{x \in \cG_{k}}  \, h_{\rho}(\tilde{\bfl},x,k) \\
     g_{\rho}(\bfl,x,k) & \leq & 
     \sup_{ \tilde{\bfl} \in \cK \cap (\cK^c)_{\epsilon} } \, 
      \max_{x \in \cG_{k}}  \, g_{\rho}(\tilde{\bfl},x,k) .
        \nonumber
\end{eqnarray}                
Next, since $ u^T_{x_1 x_2 \ldots x_{k-1}}(\bfl) \leq 
u_{x_1 x_2 \ldots x_{k-1}}(\bfl),$ we obtain by Taylor expansion
of $ (1 - \frac{u^{(T, \infty)}_{x_1  
               \ldots x_n}
                }{ u_{x_1  \ldots x_n}})^{\rho},$ 
\begin{eqnarray*}
 0 \leq  g_{\rho}(\bfl,x,k) -  g^T_{\rho}(\bfl,x,k) 
       & \leq &   
       \sum_{x_{k+1} x_{k+2} \ldots x_{n-1} \in \cG_{n-k-1}}
         \,  \frac{(u_{x_1 x_2 \ldots x_n}(\bfl))^{\rho}
               - (u^T_{x_1 x_2 \ldots x_n}(\bfl))^{\rho} }{
                 (u_{x_1 x_2 \ldots x_{k-1}}(\bfl))^{\rho} } \\
        & = & 
           \sum_{x_{k+1} x_{k+2} \ldots x_{n-1} \in \cG_{n-k-1}}
         \,  \frac{(u_{x_1 x_2 \ldots x_n}(\bfl))^{\rho}
                }{ (u_{x_1 x_2 \ldots x_{k-1}}(\bfl))^{\rho} } \cdot          
             \{ \rho \,  \frac{u^{(T, \infty)}_{x_1 x_2 
               \ldots x_n}(\bfl)
                }{ u_{x_1 x_2 \ldots x_n}(\bfl) } \\
                & & \mbox{} \; + o( \rho \frac{u^{(T, \infty)}_{x_1 x_2 
               \ldots x_n}(\bfl)
                }{ u_{x_1 x_2 \ldots x_n}(\bfl)}) \} \\*[0.1cm]
         & \leq & C_T  \, g_{\rho}(\bfl,x,k), 
\end{eqnarray*}
the reason for the last line being the convergence of the Taylor series
for sufficiently large $T.$ Clearly,
$C_T \rightarrow 0$ as $T \rightarrow \infty.$
Continuing and using Taylor expansion again brings
\begin{eqnarray*}
  0 \leq  h_{\rho}(\bfl,x,k) -  h^T_{\rho}(\bfl,x,k)
         & = &  (g_{\rho}(\bfl,x,k))^{1/(n-k)} -  
                (g^T_{\rho}(\bfl,x,k))^{1/(n-k)} \\*[0.1cm]
          & = &  (g_{\rho}(\bfl,x,k))^{1/(n-k)} - 
                    (g_{\rho}(\bfl,x,k))^{1/(n-k)} \\
                    & & \mbox{} \; \cdot 
                   [1 - \frac{g_{\rho}(\bfl,x,k) - g^T_{\rho}(\bfl,x,k)} 
               { g_{\rho}(\bfl,x,k)}]^{1/(n-k)} \\*[0.1cm] 
         & = &   (g_{\rho}(\bfl,x,k))^{1/(n-k)} \, \{
                \frac{1}{n-k} \, 
                \frac{g_{\rho}(\bfl,x,k) - g^T_{\rho}(\bfl,x,k)} 
               { g_{\rho}(\bfl,x,k)} \\*[0.1cm] 
             & & \mbox{} \; + o(  \frac{1}{n-k} \, 
                \frac{g_{\rho}(\bfl,x,k) - g^T_{\rho}(\bfl,x,k)} 
               { g_{\rho}(\bfl,x,k)}) \}  \\*[0.1cm] 
          & \leq &  h_{\rho}(\bfl,x,k) \, \{ \frac{1}{n-k} C_T +
                    o( \frac{1}{n-k} C_T ) \}, 
\end{eqnarray*}
where in the last line the upper bound of the preceding calculation
provided. In view of (\ref{gandhstrict}), the obtained upper bound
in the last display is uniform for $\bfl \in \cK_I.$
This completes the proof of the second claim.

Finally, to see the validity of
the third claim, recall from (\ref{convestimate}) 
that there is some constant $ 0 < \alpha < 1$
so that $$
            g_{\rho}(\bfl,x,k)  =
    C \exp \{ S_{n-k} \varphi_{\rho}(x) \} (1+ O(\alpha^{k}))
$$
for every $x = \ldots x_1 x_2 \ldots \in \Sigma,$
where the implicit bound in the $ O(\cdot)$ term is uniform in $x$
and the constant $C$ may be bounded
by $ C_1 < C < C_2$ with $C_i$ independent of $x$ and $n-k.$
Hence, $  S_{n-k} \varphi_{\rho; \bfl}(x)/(n-k) \rightarrow
\varphi_{\rho; \bfl}(x)$ 
as $n, k \rightarrow \infty.$  
Considering those $x$ such that 
$ S_{n} \varphi_{\rho; \bfl}(x)/n \rightarrow
\overline{\varphi}_{\rho; \bfl}$ 
as $n, k \rightarrow \infty,$ and combining the last display with 
(\ref{gandhstrict}) completes our proof. 
\end{proof}
 
Observe that for the proof of the statements in (\ref{threefctcont}),
it was not necessary to use $ \eta < 1.$
The proof shows that, for $\rho > r_u(\bfl),$ the function
$\lambda_j \rightarrow b_i(\bfl)$ is continuous as well 
for each $\bfl \in \mbox{int}(\cK)$ and $i,j \in \cA.$

\begin{corollary}
  \label{rucont}
For each $j \in \cA,$ the function $\lambda_j \rightarrow
r_u(\bfl),$ defined in (\ref{symmetryindex}), is continuous 
for each $\bfl \in \mbox{int}(\cK).$    
\end{corollary}
                     
\begin{proof}
First recall that to each $\rho >0$ there is associated
a matrix $M_{\rho}(\bfl)$ whose entries $\overb_k^{\rho}(\bfl)$
may be chosen as continuous
functions in each variable $\lambda_j$ by 
Proposition \ref{contpotential}. In that case, all matrix entries
$\overb_k(\bfl)$ are less than one. From these facts, it easily follows
that $\rho$ can be varied continuously so that the norm, say,
of $M_{\rho}(\cdot)$ stays constant as $\lambda_j$ is changed
to $\lambda_j'$ in a neighbourhood of $\lambda_j,$ and thus,
$r_u(\bfl)$ must be continuous in $\lambda_j.$ 
\end{proof}


\subsection{Mean Offspring Number $\mu_k$ of the Galton-Watson Trees}

Recall $w_x= P \{ \cD_x \}$ from (\ref{trailprobab}).
Since $w_x(\bfl) \leq u_x(\bfl),$ it follows that
the collection of functions $w_x$ satisfy a H\"{o}lder
condition as do the functions $u_x.$ Consequently,
with each function $w_x(\bfl)$ there is associated 
a well-defined potential function $\varphi^w_{\rho; \bfl}(x),$ 
constructed by the same means as described
to obtain $\varphi_{\rho; \bfl}(x)$ corresponding to $u_x.$
Next we show that the two potential functions $\varphi^w_{\rho}(x)$
and $\varphi_{\rho}(x)$ coincide for every $\rho> r_u(\bfl).$

\begin{proposition}
  \label{wmeanlimit}
For every $\bfl \in \cK_0$ and $\rho >r_u(\bfl),$
the two potential functions  $\varphi^w_{\rho}$ and $\varphi_{\rho}$ 
coincide, that is, for every $x \in \Sigma,$
\begin{eqnarray}
    \label{wusamepot}
     \varphi^w_{\rho}(x) & = & \varphi_{\rho}(x),
\end{eqnarray}      
in particular, 
\begin{eqnarray} 
 \label{wumeans} 
   \overline{\varphi}^w_{\rho} =  \overline{\varphi}_{\rho}. 
\end{eqnarray}
\end{proposition}

\begin{proof}
Let $\bfl \in \cK_0.$ 
In order to facilitate the presentation, we may carry out the proof
for $\rho=1 > r_u(\bfl)$ and look to bound 
$ \sum_{x_{k+1} \ldots x_n \in \cG_{n-k}} w_{x_1 \ldots x_n}$ instead
of $ \sum_{x_{k+1} \ldots x_n \in \cG_{n-k}} w_{x_1 \ldots x_n}^{\rho}.$ 
Minor modifications establish claim (\ref{wusamepot})
for arbitrary $\rho> r_u(\bfl).$
By definitions (\ref{hmatrices}) and (\ref{hpotential}), we collect
for all integers $k<n$ and each $x \in \Sigma,$
\begin{eqnarray*}
        \sum_{x_{k+1} \ldots x_n \in \cG_{n-k}} w_{x_1 \ldots x_n}
              & \leq & 
           \sum_{x_{k+1} \ldots x_n \in \cG_{n-k}} u_{x_1 \ldots x_n}  
               \\ & \leq &  
           ( \max_{ y \in \cG_k} u_y )  \,
                 \sum_{x_{k+1} \ldots x_n \in \cG_{n-k}} 
                 \frac{u_{x_1 \ldots x_n} }{u_{x_1 \ldots x_k}} \\
                  & \leq & ( \max_{ y \in \cG_k} u_y )  \, 
                   C \exp  \{ S_{n-k} \varphi_1(x) \} 
\end{eqnarray*}
for some finite constant $C.$ 
Therefore, for fixed $k$ and each $x \in \Sigma,$ we obtain
\begin{eqnarray} 
   \label{upperwupot}
     \limsup_{n \rightarrow \infty} \, 
      (\sum_{x_{k+1} \ldots x_n \in \cG_{n-k}} w_{x_1 \ldots x_n})^{1/n}
              & \leq   \limsup_{n \rightarrow \infty} \,
                \exp  \{ S_{n-k} \varphi_1(x)/n \}.  
\end{eqnarray}
It remains to be shown that, for all sufficiently large $k$ and
$ x \in \Sigma,$
\begin{eqnarray*}  
         \liminf_{n \rightarrow \infty} \,   
          (\sum_{x_{k+1} \ldots x_n \in \cG_{n-k}}
            w_{x_1 \ldots x_n})^{1/n} 
          & \geq &  
            \liminf_{n \rightarrow \infty} \, 
                \exp  \{ S_{n-k} \varphi_1(x)/n \}.
\end{eqnarray*} 
This will imply claim (\ref{wusamepot}) because $
\lim_{n \rightarrow \infty} S_{n-k} \varphi_1(x)/n  = \varphi_1(x)$
exists  for $ 1 > r_u(\bfl).$
                  
We shall mimic the construction of the proof of Proposition 1 
in \cite{lase2}. 
For each $x \in \cG_m,$ let $v_{\{x, k\}}$
denote the probability that there is an infection trail from the
root $1$ to $x$ that remains within distance $k$ of the geodesic
segment from $1$ to $x.$ Clearly, $ v_{\{x, k\}} \leq u_x $
and, as $k \rightarrow \infty,$ $v_{\{x, k\}} \uparrow u_x,$
that is, $ \lim_{k \rightarrow \infty} v_{\{x, k \}} =
u_x.$ 

Thus, by definitions (\ref{hmatrices}), 
(\ref{hpotential}), and (\ref{bmatrix}), if
we write $ u_{x_1 \ldots x_m} =  u_{x_1 \ldots x_{k}} 
(u_{x_1 \ldots x_m} / u_{x_1 \ldots x_{k}}),$ then for any
$\varepsilon >0,$ for all $m$ and $k$ sufficiently large,
\begin{equation}
  \label{vb}
   v_{\{x, k\}} \geq   (1 - \varepsilon) u_{x_1 \ldots x_k} 
          \prod_{j=k+1}^m b_{x_j},
\end{equation}
where each $ b_{x_j}= b_{x_j}(\bfl, x_1 \ldots x_{k})$ depends on
$x_1 \ldots x_{k}$ and $m.$ 

For every vertex $z \in \cG_k$ at distance $k$ from the root $1,$
define $\alpha_z$ to be the probability that there is a {\em
direct} infection trail from the root $1$ to $z,$ that is,
a trail which follows the geodesic segment from $1$ to $z.$
It is apparent that $\alpha_z > 0. $

Now choose $x=x_1 x_2 \ldots x_n \in \cG_n $
so that $n = 2k +j +Nm $ for some $ 0 \leq j \leq m-1$ and
let $y_0, y_1, \ldots , y_N$ be the vertices on the geodesic
segment from the root $1$ to $x$ such that $y_i \in \cG_{k+im}.$
Suppose that all of the following events take place:
(A) There is a direct infection trail from $1$ to $y_0$
that reaches $y_0$ at a stopping time $S_0,$ (B) for each
$i= 0,1, \ldots, N-1,$  there is an infection trail from
$y_i$ to $y_{i+1}$ that begins at time $S_i$ and ends at
time $S_{i+1},$ which remains within distance $k$ of the geodesic
segment from $y_i$ to $y_{i+1},$
(C) there is a direct infection trail from
$y_N$ to $x,$ beginning at time $S_N. $ (If we let $S_i$
be the first time after $S_{i-1}$ that such a path arrives
at $y_i,$ then the random times $S_i$ are stopping times.)
Observe that concatenating the infection trails (A), (B),
and (C) constitutes a downward infection trail from the root
$1$ to the vertex $x.$  Because the $S_i$ are
stopping times and the events (A), (B), and (C) occur on 
nonoverlapping parts of the percolation structure, by the 
strong Markov property and the monotonicity properties,
we have
$$
    w_x \geq \alpha_{y_0} \, ( \prod_{i=0}^{N-1}
               v_{\{y_i^{-1} y_{i+1},k\}} ) \, 
               \alpha_{y_N^{-1}x}.
$$
Notice that $ y_0 (y_0^{-1} y_{1}) (y_1^{-1} y_{2}) \ldots 
           (y_{N-1}^{-1} y_{N}) = y_N .$
Hence by (\ref{vb}), we have shown that
$$
   w_x \geq \alpha_{y_0}  \alpha_{y_N^{-1}x} \, (1 - \varepsilon)^N
                    (u_{x_1 \ldots x_k})^N \, 
                     \prod_{j=k+1}^{n-k-j} b_{x_j}  . 
$$              
If we write $\alpha_1 = \min_{y_ \in \cG_k} \alpha_y,$
     $\alpha_2 = \min_{y \in \cG_{k+j}} \alpha_y,$ 
     and  $u_* = \min_{x \in \cG_k } u_{x},$  by (\ref{hpotential}), 
     we obtain
\begin{eqnarray*}
 \sum_{x_{k+1} \ldots x_n \in \cG_{n-k}} w_{x_1 \ldots x_n} 
    & \geq &   \alpha_1 \alpha_2 (1 - \varepsilon)^N 
                  (u_*)^N \, 
              \sum_{x_{k+1} \ldots x_n \in \cG_{n-k}} 
            \prod_{j=k+1}^{n-k-j} b_{x_j} \\*[0.1cm]
           & \geq &  \alpha_1 \alpha_2 (1 - \varepsilon)^N 
                  (u_*)^N \ C'  \exp  \{ S_{n-2k-j} \varphi_1(x) \} 
\end{eqnarray*}
for some positive constant $C'.$ Note that the employed constant $C'$
(as well as $C$ earlier) takes care of the discrepancy in the index 
sets of the summations. Then
$$ 
 \liminf_{ n \rightarrow \infty} \,      
 (  \sum_{x_{k+1} \ldots x_n \in \cG_{n-k}} w_{x_1 \ldots x_n})^{1/n} 
     \geq   [(1 - \varepsilon) 
                  u_*]^{1/m} \,  
              \liminf_{ n \rightarrow \infty} \,  
               \exp  \{ S_{n-2k-j} \varphi_1(x)/n \}.  
$$
As $\varepsilon > 0$ is arbitrary and as $m$ can be chosen sufficiently
large for suitable $k,$ the desired result (\ref{wusamepot}) now follows. 
Finally,  (\ref{wumeans}) is an immediate consequence of definition
(\ref{phimeans}). 
\end{proof} 

\begin{lemma}
  \label{meanfromroot}
As $n \rightarrow \infty,$
\begin{equation}
   \label{limmeanfromroot}
   ( \, \sum_{x \in \cG_n} u_x(\bfl) \, )^{1/n}  
   \rightarrow  \exp \overline{\varphi}_{1; \bfl} = \theta_1(\bfl).
\end{equation}   
\end{lemma}

\begin{proof}
To simplify our exposition, we will present the proof in terms of
the $u_x$ and $\varphi_{1; \bfl},$ rather than lifting the
calculation to $u_x^2$ and $\varphi_{2; \bfl}.$ However, observe
that the two calculations are related via the functions $b_j.$
For every integer $k>0$ and $x  \in \Sigma,$
by subadditivity and (\ref{hpotential}), 
\begin{eqnarray*}
    \sum_{x \in \cG_n} u_x(\bfl) & \leq & 
              \sum_{x_{k+1} x_{k+2} \ldots x_{k+n}
                     \in \cG_n} \frac{u_{x_1 x_2 \ldots x_{k+n}}(\bfl)}  
                          {u_{x_1 x_2 \ldots x_{k}}(\bfl)} \\
           & \leq &  C \exp \{ S_n \varphi_{1; \bfl}(x) \}
\end{eqnarray*}
for some finite constant $C.$ 
Hence, for every $x \in \Sigma,$
$$
   \limsup_{n \rightarrow \infty}  \,
   (  \sum_{x \in \cG_n} u_x(\bfl))^{1/n}
         \leq   \limsup_{n \rightarrow \infty} \, 
                    \, \exp \{ S_n \varphi_{1; \bfl}(x)/n \}.
$$
Since this inequality holds for every $x \in \Sigma,$ it must
hold for those $x \in \Sigma$ with $S_n \varphi_1(x)/n \rightarrow
\overline{\varphi}_1$ as $n \rightarrow \infty.$
Thus, letting $k \rightarrow \infty,$ it is easily derived that
$$
\limsup_{n \rightarrow \infty} \,
   (  \sum_{x \in \cG_n} u_x(\bfl))^{1/n}
         \leq  \exp \overline{\varphi}_{1;\bfl}.     
$$
 
To show the reverse direction, fix $\varepsilon >0,$ and
define
$$ 
Q_{\varepsilon}  =
     \{ x_1 x_2 \ldots x_k \in \cG_k : \, 
       \sum_{x_{k+1} \ldots x_n \in 
       \cG_{n-k}} u_{x_1 x_2 \ldots x_n} (\bfl) \geq 
        \exp \{ (\overline{\varphi}_{1; \bfl}
           - \varepsilon)(n-k)\}  \} .
$$ 
It is obvious that for sufficiently large $n$ and fixed sufficiently large
$k,$ the set $ Q_{\varepsilon} $ is nonempty. Therefore,
\begin{eqnarray*}
     \sum_{x \in \cG_n} u_x(\bfl) & \geq & \sum_{x_1 \ldots x_k \in 
         Q_{\varepsilon} }  u_{x_1 x_2 \ldots x_n}(\bfl) \\
        & \geq &   \exp \{ (\overline{\varphi}_{1; \bfl}
           - \varepsilon)(n-k) \}.
\end{eqnarray*}
Whence,
\begin{eqnarray*}
    \liminf_{n \rightarrow \infty} \,
      ( \sum_{x \in \cG_n} u_x(\bfl))^{1/n}
        & \geq &   \exp  \{ \overline{\varphi}_{1; \bfl}
           - \varepsilon \}.
\end{eqnarray*}
As $\varepsilon >0$ was arbitrary, this verifies claim
(\ref{limmeanfromroot}). 
\end{proof}   

\begin{corollary}
  \label{exactmean} 
If we recall $\cL_n$ from Section 2.2, then
\begin{equation}
    \label{almostgwmeans} 
\lim_{n \rightarrow \infty} 
     (\sum_{x \in \cL_n}  w_{x}(\bfl) )^{1/n}
 = \lim_{n \rightarrow \infty} 
        (\sum_{x \in \cG_n} w_{x}(\bfl))^{1/n} = 
        \exp  \overline{\varphi}_{1; \bfl} = \theta_1(\bfl).
\end{equation}
\end{corollary}

\begin{proof}
Again, our proof is in terms of the $u_x$ and $\varphi_{1; \bfl}$ 
(see the remark at the outset of the proof of
Lemma \ref{meanfromroot}).
In view of Proposition \ref{wmeanlimit}
and Lemma \ref{meanfromroot}, the leftmost equality in 
(\ref{almostgwmeans}) remains to be verified.  
This may be accomplished by showing that
  $$ \lim_{n \rightarrow \infty} 
        (\sum_{x \in \cL_n} u_{x})^{1/n} 
       =\lim_{n \rightarrow \infty} 
     (\sum_{x \in \cG_n}  u_{x} )^{1/n}.
$$             
But in light of (\ref{addingentries}),
\begin{eqnarray}
   \label{udifferentstart}  
  \sum_{x \in \cG_n}  u_x  & =& {\bf 1}^t H_1 (n;1; \cdot) {\bf 1}, \\
  \sum_{x \in \cL_n}  u_x  & = & {\bf u}_a^t H_1 
                       (n; 1; \cdot) {\bf 1},   \nonumber 
\end{eqnarray}
where ${\bf 1}$ is the vector all of whose entries are $1$ and 
${\bf u}_a$ is the vector with entry $0$ in the $a^{-1}$ slot
and all other entries $1.$ Since $ H_1 (n;1; \cdot)$ is a
Perron-Frobenius matrix with lead eigenvalue  
   $\exp S_{n} \varphi_{1}(1),$ taking the $n$-th root and  
 $n \rightarrow \infty$ on both sides in both lines of
 (\ref{udifferentstart}) in combination with Lemma 
 \ref{meanfromroot} proves claim (\ref{almostgwmeans}).
\end{proof}
 
\begin{corollary}
 \label{almostmean} 
The mean offspring numbers $\mu_k$ for the Galton-Watson trees
$\tau_k,$ as defined in Section 2.4, satisfy
\begin{equation}
   \label{gwmeans}
         \lim_{k \rightarrow \infty} \mu_k^{1/k} = 
         \exp  \overline{\varphi}_{1; \bfl}  = \theta_1(\bfl) . 
\end{equation}
\end{corollary}

\begin{proof}
The mean offspring number $\mu_k$ for the Galton-Watson tree
$\tau_k$ is, by construction, 
$$
    \mu_k = \sum_{x \in \cL_k^*} w_x.
$$
This sum  differs from the one on the lefthand side of the left equality 
in (\ref{almostgwmeans}) in that 
the smaller index set $\cL_k^*$ replaces $\cL_k,$ thus, the
vertices $x$ have word representations ending in the letter
$a.$ Consider a vertex $x \in \cL_k.$ Clearly, there is positive
probability $\rho,$ independent of $x,$
that a vertex $y \in \cL_{k+2}^* \cap \cT(x) $ is
infected within two time units. This brings
$$
    \sum_{x \in \cL_{k+2}^*} w_x \geq \rho \sum_{x \in \cL_k} w_x.
$$
Corollary \ref{exactmean} thus finishes the proof of (\ref{gwmeans}).
\end{proof}

For all integers $r,n \geq 1,$ define $Z_n(r) = Z_n = \vert V_n(r) \vert$ 
to be the cardinality of the $n$-th generation of $\tau_r.$
For $n=1$ for instance, $EZ_1(r) =  \sum_{x \in \cL_r^*} w_x.$

\begin{lemma}
   \label{infmeansupercrit} 
Suppose that $\theta_1(\bfl) >1 .$   
Then for any $1 < \gamma_* < \theta_1(\bfl),$ we have
$$
  \liminf_{k \rightarrow \infty}
       P \{ Z_1(2^k) > (\gamma_*)^{2^k} \}
              = \rho > 0.  
$$              
\end{lemma}

\begin{proof}
This proof is essentially the same as the one in \cite{lase2}
(Corollary 3).
We shall put it in our context.
Let $x \in \cG_{2^k}.$ Fix some $1 \leq m \leq k,$ and let
$x_0 = 1, x_1 , \ldots ,x$ be the vertices along the geodesic
segment from $1$ to $x$ at distance $i 2^m$ for $i=0,1, \ldots, 
2^{k-m}.$ If for every $i,$ there is a downward infection trail
$\xi_i$ from $x_i$ to $x_{i+1}$ that begins at the time of
termination of $\xi_{i-1},$ then there is a downward infection
trail from $1$ to $x,$ thus, $x \in V_1(2^k).$ Therefore,
$x \in V_{2^{k-m}}(2^m)$ implies $x \in V_{1}(2^k).$ Hence,
$Z_{1}(2^k) \geq  Z_{2^{k-m}}(2^m).$ 
By construction, $V_n(2^m)$ is a Galton-Watson process with
mean offspring number $EZ_n(2^m) =  \sum_{x \in \cL_{n 2^m}^*} w_x.$
A standard result from the theory of Galton-Watson processes
tells us that
$$
    \lim_{n \rightarrow \infty} \frac{Z_n(2^m)}{EZ_n(2^m)} = Z
$$     
exists and, because the offspring distribution has finite
support, $Z > 0$ almost surely on the event of nonextinction
(see \cite{atne}, Theorem 2, Section 6). This means
that $P \{ Z_n(2^m) > ((\gamma_*)^{2^m})^n \mbox{ eventually} \} >0,$
which together with $Z_{1}(2^k) \geq  Z_{2^{k-m}}(2^m) $
concludes the proof.
\end{proof}

\bigskip
{\bf Proof that $r_u \leq 1 $ at the Transition to Survival.} 
Consider the set $\cY_m$ of all vertices in $\cG_m$
ever to be infected. 
In this paragraph and again when we address the
behaviour of the contact process at the first phase transition, 
we will exploit the fact that, on homogeneous trees of degree 
larger than $2,$
whenever $  E \vert \cY_m \vert \rightarrow \infty$ as $m \rightarrow
\infty,$ a (labelled) Galton-Watson tree may be embedded in the
set of vertices ever to be infected, which gives rise to a 
supercritical Galton-Watson process. This Galton-Watson process
grows without bound with positive probability, which implies
that, with positive probability, $\vert A_t \vert \rightarrow \infty.$

Clearly,
\begin{eqnarray*}
    E \vert \cY_m \vert & = &
                   E \sum_{x \in \cG_m} I_{\{x \mbox{ is ever
                               infected} \}} \\
                & = & \sum_{x \in \cG_m} E I_{\{x \mbox{ is ever
                               infected} \}} \\
                 & = & \sum_{x \in \cG_m} u_x  ,
\end{eqnarray*}
where $I_{ \{ . \}}$ denotes the indicator function.
In view of Lemma \ref{meanfromroot}, we know that
for $\bfl \not \in \cR_1,$ the lead eigenvalue
$\theta_1(\bfl) \geq 1,$ because $\theta_1(\bfl) < 1$ would imply that  
$ \sum_{m=1}^{\infty} E \vert \cY_m \vert < K$
for some constant $K < \infty,$ and thus, the contact process
would not survive with positive probability, which contradicts
our assumption that $\bfl \not \in \cR_1.$ As a consequence,
$r_u(\bfl) >1 $ for $\bfl \in \mbox{int}(\cR_2 \cup \cR_3).$ 
 
\begin{lemma}
 \label{meanone}
For $\bfl \in \overR_1 \cap \overline{\cR_1^c} \cap \mbox{int}(\cK),$ 
\begin{eqnarray}
   \label{rulessthanone}
   \theta_1(\bfl) &  \leq & 1, \\
           r_u(\bfl) & \leq & 1. \nonumber
\end{eqnarray}
\end{lemma}

\begin{proof}
Throughout the proof, let $\bfl \in \mbox{int}(\cK).$
It is enough to prove the first claim.
For $d=1,$ claim (\ref{rulessthanone}) is obvious because
$1 \geq \beta(\bfl) = \theta_1(\bfl).$ Thus, let $d>1.$
We suppose that $\theta_1(\bfl) > 1$ for 
$\bfl \in \overR_1 \cap
\overline{\cR_1^c} $ and proceed by contraposition.
 By Proposition \ref{contpotential},
for each $j \in \cA,$ the function 
$\lambda_j \rightarrow \theta_1(\bfl)$ is continuous.
Furthermore, by Proposition \ref{potmonotone}, for any $\bfl_*$
in a direction of decrease for $\bfl,$ we collect
 $\theta_1(\bfl_*) < \theta_1(\bfl).$ 
By a combination of these two properties, because 
$\theta_1(\bfl) > 1,$ there exists a $\bfl_*$
in the interior of $\cR_1$ such that 
$\theta_1(\bfl_*) > 1.$ Thus, for this  $\bfl_*,$ 
$ E \vert \cY_m \vert \rightarrow \infty$ as $m \rightarrow \infty.$
Now, by Lemma \ref{infmeansupercrit}, for sufficiently large $k$
 (as $k$ runs through powers of $2$), the Galton-Watson process
$\{ Z_n(k) \}_{n \geq 0}$ is supercritical and the corresponding
Galton-Watson tree
has positive probability to grow to infinity. Because the Galton-Watson
process is dominated by the number of vertices
ever to be infected, it follows that, with positive probability,
the contact process survives. However, this contradicts our assumption
that $\bfl_* \in \mbox{int}(\cR_1) \subset \cR_1.$  Hence,
it must be the case that $\theta_1(\bfl) \leq 1.$
\end{proof}


\subsection{Asymptotics of the $u_x$ and Related Functions}

Inspired by some ideas in \cite{sch1},
we collect some bounds between various key
functions that are equationally interrelated. 
In fact, our Lemmata \ref{ubar}, \ref{ixlemma}, 
Proposition \ref{uinequalities} and their proofs 
are adjusted versions of the reasoning in \cite{sch1}.
We will prove the results in this section, Section 5.4 and in
Propositions \ref{betarecursion} and \ref{equalbeta} under the
following Standing Hypothesis, which will shortly turn out to be
extra (for more on this, see Corollary \ref{valuesattrans}).

\medskip
{\bf Hypothesis I. }
Assume that, for every $\epsilon >0,$ the exponent $r_u(\bfl) < 2$ 
for each $\bfl \in \cK \cap (\cK^c)_{\epsilon}.$
 
\medskip 
A consequence of Hypothesis I is that, 
for every $\epsilon >0,$ we have $r_u(\bfl) <2$ for each
$\bfl \in \cK_I,$ due to the monotonicity properties of the
contact process, where  
$\cK_{I} = \cK \setminus (\cK^c)_{\epsilon}$ 
denotes the set of points
in $\cK$ which are at a distance to $ \cK^c$ at least $\epsilon.$
 
Next define 
\begin{equation}
   \label{uover}
         \overline{u}_x = \overline{u}_x(\bfl) = \int_0^{\infty}
              P \{x \in \at \} dt.
\end{equation}
Obviously, $u_x \leq \overline{u}_x$ for every vertex $x \in \cG.$

\begin{lemma}
  \label{ubar}
Let $\epsilon >0.$ Under Hypothesis I, there exists 
some constant $0 < c < \infty$ such that for every
$\bfl \in \cK_{I}$ and $x \in \cG_n,$
$$
 \overline{u}_x  \leq c (n+1) u_x.
$$  
\end{lemma}

\begin{proof}
Let $a_{\epsilon}$ denote the probability for the isotropic
contact process with infection parameter $\epsilon>0$ and initial infection
at the root vertex that between time $0$
and $1$ there is no recovery mark * 
at the root $1$ and that, for some $i \in \cA,$ the root vertex
infects its neighbour $i.$  
Thus, $ a_{\epsilon} = e^{-1}(1-e^{-\epsilon}) > 0.$ 
Fix $\bfl \in \cK_I$ and fix $\epsilon_*>0$ such that each $
\lambda_k \geq  \epsilon_*$ ($\epsilon_* = \epsilon / \sqrt{2d}$
should suffice). Then for every $x= x_1 x_2 \ldots
x_n \in \cG_n$ and $ y=y_1 y_2 \ldots y_t  \in \cG_t$
 such that $\vert xy  \vert = n+t,$ it follows that
\begin{eqnarray}
     u_x (\bfl) & \geq & (a_{\epsilon_*})^n  \nonumber \\
  \label{xinfectlow}
   P \{ 1 \in A_{t+n}  \} & \geq & 
          (a_{\epsilon_*})^n \, P \{ x \in A_t \} .
\end{eqnarray}
Fix $\tilde{\bfl} \in \cK \cap (\cK^c)_{\epsilon}$ in a direction of
increase for $\bfl$
(so that each $ \beta_j(\tilde{\bfl}) \geq \beta_j(\bfl) >0$). 
Thus, $r_u(\tilde{\bfl}) <2.$
Recall that $\eta(\cdot)$
is nondecreasing in each variable $\lambda_j$ and
from Proposition \ref{etalessthanone} that $\eta(\bfl) < 1$
for $\bfl \in \mbox{int}(\cK).$ Thus,
$ P \{ 1 \in A_t \} \leq \eta^t(\bfl) \leq 
\eta^t(\tilde{\bfl}) < 1$ for every $t>0.$
This observation combined with (\ref{xinfectlow}) gives
\begin{equation}
  \label{xinfectup}
    P \{ x \in A_t \} \leq a_{\epsilon_*}^{-n} 
                                \eta(\tilde{\bfl})^{t+n} 
   \leq    a_{\epsilon_*}^{-n} \eta(\tilde{\bfl})^{t}.
\end{equation}
Recall that $u_x \leq \overline{u}_x.$
Therefore, if $K > 0,$ by (\ref{xinfectup}), 
\begin{eqnarray*}
     \overline{u}_x(\bfl)  & = & 
             \int_0^{\infty} P \{ x \in A_t \} dt \\
              & \leq& 
             u_x(\bfl) Kn + \int_{Kn}^{\infty} P \{ x \in A_t \} dt \\
            & \leq &  u_x(\bfl) Kn + a_{\epsilon_*}^{-n} 
                          \eta (\tilde{\bfl})^{Kn} /
       \log(1/\eta (\tilde{\bfl})).
\end{eqnarray*}
For  $K$ sufficiently large, in particular, such that 
$\eta(\tilde{\bfl})^K/ a_{\epsilon_*} < 1,$ we obtain
$$
    \overline{u}_x(\bfl)  \leq 2Kn u_x(\bfl) \leq c(n+1) u_x(\bfl)
$$
for some $0 < c < \infty$ and for each
$n > 0.$ The case $n=0$ is an easy instance of $\eta(\bfl) \leq
\eta(\tilde{\bfl}) < 1.$ This ends our proof.
\end{proof}

For any $x \in \cG,$ let $X_x$ denote the total number of infection trails 
(infection arrows in the percolation structure without
recovery marks) leading to $x$ and let
\begin{eqnarray*}
    I_1 & = &  I_1(\bfl) = E X_1 +1 \\
    I_x & = &  I_x(\bfl) = E X_x.
\end{eqnarray*}
It is apparent that $ u_x \leq I_x$ for every $x \in \cG.$

\begin{lemma}
 \label{ixlemma}
For every $x, y \in \cG $ such that $ \vert x y \vert = \vert x \vert + 
\vert y \vert, $
\begin{equation}
 \label{totalinfect}
   I_{x y} \leq I_x I_y. 
\end{equation}
\end{lemma}

\begin{proof}
Clearly, (\ref{totalinfect}) holds for $x = 1 $ or $y = 1 $
(recall that $1 $ denotes the root). 
As was done in \cite{sch1}, we argue by comparing the underlying
process to a multitype contact process that evolves as follows.
Think of an infected vertex as a site hosting
a particle. Infected vertices can carry infections of types $0,1,2,
\ldots, $ thus, particles can be of type $0,1,2, \ldots,$ with
the initial particle at the root being of type $0.$
Particles of type $1,2, \ldots$ evolve as independent contact processes,
independent of particles of type $0,$ in particular, particles of
different types can coexist at the same vertex at the same time.
Fix a vertex $x \in \cG.$ Particles of type $0$ evolve as a contact process,
except for the fact that they cannot infect vertex $x.$
The first time a $0$-particle attempts to infect $x,$ a particle
of type $1$ is placed at $x.$ The second time a $0$-particle
attempts to infect vertex $x,$ a particle of type $2$ is placed
at $x.$ And, so forth. 

This multitype contact process dominates the contact process
in the sense that, when there is a particle of the contact process  
at some vertex, then there is at least one particle of the multitype 
contact process at the same vertex. Let $X_x'$ denote the total
number of attempts by a $0$-particle to infect $x$ (thus, the number
of types of particles distinct from $0$), and let
$X_{xy}''$ denote the total number of attempts by a particle of any type
to infect the vertex $xy,$ where $y$ is such that $\vert xy \vert = 
\vert x \vert + \vert y \vert.$ 
Thus, the chain of inequalities
\begin{eqnarray*}
   I_{xy} & \leq & E(X_{xy}'') \\
             & = & E(E(X_{xy}'' \vert X'_x)) \\
            & = & E(X'_x E(X_y)) \\
            &  = & E(X'_x) E(X_y) \\
            & \leq & E(X_x) E(X_y) = I_x I_y
\end{eqnarray*}
finishes the proof of (\ref{totalinfect}).
\end{proof}

Similarly as for $u_x,$
a subadditivity argument shows that, for any integer $k,$ 
 $x \in \cG_k,$
and each $y_n = x x \ldots x \in \cG_{nk},$ the limit
\begin{equation}
  \label{ilimit}
   \lim_{ n  \rightarrow \infty} I_{y_n}^{1/n} = \tilde{\beta}_x = 
              \tilde{\beta}_x(\bfl)
\end{equation}
exists and that $I_{y_n} \geq \tilde{\beta}_x^n$ for each integer
$n \geq 0.$ 

Recall that an infected vertex $x \in \cG$ attempts to infect its
nearest neighbour $xa$ at the rate $\lambda_a$ for $a \in \cA.$
Therefore, enlarging on equation (3.2) in \cite{sch1}, yields the
following recurrence relation,  
\begin{equation}
  \label{recursion}
       I_x = \sum_{a \in \cA}  \lambda_a \overline{u}_{xa}.
\end{equation}

\begin{proposition}
  \label{uinequalities}
Let $\epsilon >0$ and assume Hypothesis I.
Then there are some positive finite
constants $ c_1(\epsilon)$ and $ c_2(\epsilon)$ such that
 for each $\bfl \in \cK_{I},$ each integer $k>0,$ and every 
$x  \in \cG_k,$
\begin{equation}
   \label{chainineq}
     I_{x}(\bfl)
       \leq  c_1(\epsilon) 
        \overline{u}_x(\bfl)
        \leq  c_2(\epsilon) (n+1) u_x(\bfl)
         \leq  c_2(\epsilon) (n+1) I_x(\bfl).
\end{equation}
For each
 $\bfl \in \mbox{int}(\cK),$ each $ x \in \cG_k,$ and each
$y_n =x x \ldots x \in \cG_{nk},$ we have
\begin{equation}
  \label{uilimitsbeta}
 \lim_{n \rightarrow \infty} \overline{u}_{y_n}(\bfl)^{1/n}
            = \lim_{n \rightarrow \infty} I_{y_n}(\bfl)^{1/n}
       = \tilde{\beta}_x(\bfl)
        =   \lim_{n \rightarrow \infty} u_{y_n}(\bfl)^{1/n}
          =\beta_x(\bfl).
\end{equation} 
Furthermore, there are some constants 
$0 <  C_i(\epsilon) < \infty,$ $i=1,2,3,$
such that for $\bfl \in \cK_{I},$  for all integers $n,k >0,$
every $x \in  \cG_k,$ and $y_n=x x \ldots x \in \cG_{nk},$
\begin{eqnarray}
   \label{ubarineq}
  \frac{\beta_x(\bfl)^n}{C_1(\epsilon) (n+1)} & \leq & u_{y_n}(\bfl) \leq
                    \beta_x(\bfl)^n \\*[0.1cm]
     \frac{\beta_x(\bfl)^n}{C_2(\epsilon) }
      & \leq & \overline{u}_{y_n}(\bfl) \leq
                   C_3(\epsilon) (n+1)  \beta_x(\bfl)^n \nonumber \\*[0.1cm] 
      \beta_x(\bfl)^n & \leq & I_{y_n}(\bfl) \leq
                   C_1(\epsilon) (n+1)  \beta_x(\bfl)^n. \nonumber 
\end{eqnarray}
\end{proposition}

\begin{proof}
Let $\epsilon > 0$ and $\bfl \in \cK_{I}.$
Again fix some $\epsilon_*>0$ such that each $\lambda_i \geq \epsilon_*$
($\epsilon_* = \epsilon / \sqrt{2d}$ will do).  
Then by (\ref{xinfectlow}), for every
$x=x_1 x_2 \ldots x_k \in \cG_k$ and any $x_{k+1} \in \cA,$
we have
\begin{eqnarray*}
    \overline{u}_{x_1 x_2 \ldots x_k}(\bfl)
       & \geq & \int_1^{\infty} P \{ x_1 x_2 \ldots x_k
                         \in A_t \} dt \\
                   & \geq &  \int_1^{\infty}
                        a_{\epsilon_*} P \{ 
                       x_1 x_2 \ldots x_{k+1} \in A_{t-1} \} dt \\
                    & = & \int_0^{\infty}  a_{\epsilon_*} P \{ 
                       x_1 x_2 \ldots x_{k+1} \in A_t \} dt \\
                    & = & a_{\epsilon_*}\, 
  \overline{u}_{x_1 x_2 \ldots x_{k+1}}(\bfl)
\end{eqnarray*}
(Note that $x_1 x_2 \ldots x_{k+1}$ may be in $\cG_{k-1}$).
This inequality combined with Lemma \ref{ubar} and relation
(\ref{recursion}) implies
$$
     I_{x}(\bfl)
       \leq  c_1(\epsilon) 
        \overline{u}_x(\bfl)
        \leq c_2(\epsilon) (n+1) u_x(\bfl)
         \leq  c_2(\epsilon) (n+1) I_x(\bfl)
$$
for each $x \in \cG_k$ 
and for some constants $c_1(\epsilon), c_2(\epsilon),$
both depending on $\epsilon,$ as advertized in (\ref{chainineq}).
Next apply (\ref{chainineq}) with $x=y_n,$ 
take the $n$-th root in each expression and take
limits. Since $\epsilon > 0$ is arbitrary, we conclude that,
for $\bfl \in \mbox{int}(\cK),$
$$
  \tilde{\beta}_x(\bfl) = \beta_x(\bfl)
$$   
and obtain (\ref{uilimitsbeta}).
Combining this with (\ref{uperiodicsubadd}) 
and (\ref{ilimit}) yields all of (\ref{ubarineq}), as desired.
\end{proof}

\begin{corollary}
   \label{equalpotentials}
Let $\epsilon >0.$ Under Hypothesis I, there is some positive finite
constant $ C=C(\epsilon)$ such that
 for each $\bfl \in \cK_{I},$ all integers $n-1 > k>0,$ and every 
$x= \ldots x_1 x_2 \ldots  \in \Sigma,$
\begin{eqnarray*}
   \frac{1}{C k} \, \sum_{x_{k+1}  \ldots x_{n-1} \in \cG_{n-k-1}}
         \,  \frac{u_{x_1  \ldots x_n}(\bfl)}{
                  u_{x_1  \ldots x_{k-1}}(\bfl)} 
                & \leq & 
                 g^I , \, g^{\overline{u}}   
                 \leq  
     C (n+1) \, \sum_{x_{k+1}  \ldots x_{n-1} \in \cG_{n-k-1}}
         \,  \frac{u_{x_1  \ldots x_n}(\bfl)}{
                  u_{x_1  \ldots x_{k-1}}(\bfl)}  ,                 
\end{eqnarray*}
where
\begin{eqnarray}
     g^I = g^I (\bfl, x, k) & = & 
         \sum_{x_{k+1}  \ldots x_{n-1} \in \cG_{n-k-1}}
                \,  \frac{I_{x_1  \ldots x_n}(\bfl)}{
                  I_{x_1  \ldots x_{k-1}}(\bfl)}  \nonumber \\*[0.15cm]
    g^{\overline{u}} = g^{\overline{u}} (\bfl, x, k) & = & 
            \sum_{x_{k+1}  \ldots x_{n-1} \in \cG_{n-k-1}}
         \,  \frac{\overline{u}_{x_1  \ldots x_n}(\bfl)}{
                  \overline{u}_{x_1  \ldots x_{k-1}}(\bfl)} .
\end{eqnarray}
\end{corollary}

\begin{proof}
This is an easy exercise thanks to (\ref{chainineq}). For instance,
there is some positive finite constant such that for all integers
$n > k >0$ and every  $x \in \cG_n,$  
$$   \frac{1}{C k} \, \frac{u_{x_1  \ldots x_n}(\bfl)}
           { u_{x_1  \ldots x_{k-1}}(\bfl)} 
         \leq  
    \frac{I_{x_1  \ldots x_n}(\bfl)}{ I_{x_1  \ldots x_{k-1}}(\bfl)} 
  \leq  C (n+1) \, 
  \frac{u_{x_1  \ldots x_n}(\bfl)}{ u_{x_1  \ldots x_{k-1}}(\bfl)}.   
$$
Analogously, the ratio for $\overline{u}_x$ can be squeezed in.
\end{proof}

An interesting consequence of Corollary \ref{equalpotentials}
is that under Hypothesis I
the results, obtained for the functions $u_x$
in Section 3, can be derived for the functions
$I_x$ and $\overline{u}_x,$ more precisely, H\"{o}lder continuous
matrices may be constructed 
associated with $I_x$ and $\overline{u}_x,$ respectively, as well as
potential functions $\varphi^I_{\rho}(x)$ and
$\varphi^{\overline{u}}_{\rho}(x).$ But in fact, it is an instance
of Corollary \ref{equalpotentials} that these two functions 
are the same and coincide with the function
$\varphi_{\rho}(x)$ corresponding
to $u_x,$ that is, for each $x \in \Sigma,$
\begin{equation} 
     \label{equalbmatrices} 
 \varphi^I_{\rho; \bfl}(x) = \varphi^{\overline{u}}_{\rho; \bfl}(x)
        = \varphi_{\rho; \bfl}(x). 
\end{equation}
Thus, the lead eigenvalues of the corresponding matrices 
$B_{\rho}$ coincide for the three functions
$u_x,$ $ I_x,$ and $\overline{u}_x$ and their matrix entries 
may be chosen the same.


\subsection{Continuous Limits via Subadditivity}

In the isotropic case, continuity of the function $\beta(\cdot)$
in $\bfl$ on $(0, \lambda_2)$ has been proven in \cite{sch1}.
An analogue holds in the anisotropic case with an analogous proof.

\begin{proposition}
  \label{continuity}
Under Hypothesis I,   
for all $j \in \cA$ and $x \in \cG,$
the function $\lambda_j \rightarrow  
\beta_x(\bfl)$ is continuous for every $\bfl \in \mbox{int}(\cK)$
and left-continuous for $\bfl \in {\bf R}_+^d.$
\end{proposition}

\begin{proof}
The proof for the isotropic contact process, given in \cite{sch1},
carries over with minor adaptations. Since the reasoning is short, we 
present the details here.

It suffices to show that, for all $ j \in \cA$ and $x \in \cG,$ \, 
 (a) $\lambda_j \rightarrow \beta_x(\bfl)$ is left-continuous
     on  ${\bf R}_+^d,$ and \,
 (b) $\lambda_j \rightarrow \beta_x(\bfl)$ is right-continuous
     for every $\bfl \in \mbox{int}(\cK).$
     
To show the first of these claims, note that, if for each $ x \in \cG_k,$
 $y_n = x x \ldots x \in \cG_{kn},$ and every $T>0,$ we set
$
  u_{y_n}^T = u_{y_n}^T(\bfl) =
    P \{ y_n \in A_t \mbox{ for some } t \in [0,T] \},
$
then
$$
  \beta_x(\bfl) = \sup_{ n \geq 1} (u_{y_n}(\bfl))^{1/n} 
                = \sup_{n \geq 1} \sup_{T \geq 0} (u_{y_n}^T(\bfl))^{1/n}.
$$ 
It is obvious that $u_{y_n}^T(\cdot)$ is a continuous function in each
$\lambda_j,$ and thus, $\beta_x(\cdot)$ is lower-semi-continuous
in each $\lambda_j.$ Since the function $\beta_x(\cdot)$ is 
nondecreasing in each $\lambda_j,$ it is left-continuous in each
$\lambda_j.$ 

To verify the second claim, for each $x \in \cG_k$ and for
$y_n = x x  \ldots x \in \cG_{kn},$ we may define
$$
   \overline{u}_{y_n}^T(\bfl) = \int_0^{T} P \{ y_n \in A_t \} \, dt. 
$$
Fix $\epsilon > 0$ and $\epsilon_*>0$ such that each $\lambda_j \geq 
\epsilon_*.$     
For fixed $\bfl,$ as $T \rightarrow \infty,$ the convergence of
$ \overline{u}_{y_n}^T(\bfl)$ to $\overline{u}_{y_n}(\bfl)$ can be 
controlled uniformly for  $\bfl \in \cK_I,$ where $\cK_I$ was
described earlier.
To see this, fix $\bfl \in \cK_I$ and 
    $\tilde{\bfl} \in \cK \cap (\cK^c)_{\epsilon}$ in a direction of
    increase for $\bfl,$ and 
recall $ a_{\epsilon} = e^{-1}(1-e^{-\epsilon}) > 0$ and
(\ref{xinfectup}). Then
\begin{eqnarray*}
  0 \leq \overline{u}_{y_n}(\bfl) -  \overline{u}_{y_n}^T(\bfl) & = &
         \int_T^{\infty} P \{ y_n \in A_t \} \, dt \\
          & \leq & ( a_{\epsilon_*})^{-kn} \,  \int_T^{\infty} 
          \eta(\tilde{\bfl})^{t+kn} \, dt.
\end{eqnarray*}      
Since, by Proposition \ref{etalessthanone} and Hypothesis I, 
we have $\eta(\bfl) < 1$ for $\bfl \in \mbox{int}(\cK),$
the expression on the lefthand side
of the last display vanishes as $ T \rightarrow \infty.$ Thus,
$ \overline{u}_{y_n}^T(\cdot)$ converges to $\overline{u}_{y_n}(\cdot)$ 
uniformly on $\cK_I.$
It is apparent that, for each $j \in \cA,$
 the function $\lambda_j \rightarrow
 \overline{u}_{y_n}^T(\bfl)$ is continuous, and thus,  
 the function $\lambda_j \rightarrow
 \overline{u}_{y_n}(\bfl)$ is continuous as well on      
$\cK_I.$ From (\ref{chainineq}) and the left inequality in
the second line of (\ref{ubarineq}), 
we collect $(C_2(\epsilon) \overline{u}_{y_n}(\bfl))^{1/n} 
\geq \beta_x(\bfl)$
and, also, from Proposition  \ref{uinequalities}, 
               $ \lim_{n \rightarrow \infty}
   (C_2(\epsilon) \overline{u}_{y_n}(\bfl))^{1/n} = \beta_x(\bfl).$
This implies that 
$$ 
    \beta_x(\bfl) = \inf_{n \geq 1}
        (C_2(\epsilon) \overline{u}_{y_n}(\bfl))^{1/n}. 
$$
Hence, the function $\lambda_j \rightarrow \beta_x(\bfl)$ is 
upper-semi-continuous on $\cK_I.$     
Since the function $\lambda_j \rightarrow \beta_x(\bfl)$ is 
nondecreasing, it is also right-continuous on this region. 
Since $\epsilon > 0$ is arbitrary, the advertized claim follows. 
\end{proof}


\subsection{Discontinuity}

From recursion (\ref{recursion}), we collect the following
fundamental relation.

\begin{proposition}
 \label{betarecursion}
Let $\epsilon >0.$ If Hypothesis I holds, then
for every $\bfl \in \cK_I,$ every $ i \in \cA,$
and $ x = x_1 x_2 \ldots x_{n-1} i \in \cG_n,$
there is some $D_{\bfl}(x)$ (independent of $i$) so that 
\begin{eqnarray}
   \label{recurrence}
   b_i(\bfl) & = &  \lambda_i D_{\bfl}(x) 
                  + b_i(\bfl)  D_{\bfl}(x) \, \sum_{j \in \cA -\{ i\}}
                   \, \lambda_j \, b_j(\bfl), 
\end{eqnarray}
where the $b_j= b_j(\bfl) $ are defined in (\ref{bmatrix})
and  $ C_1 \leq  D_{\bfl}(x) \leq C_2 (n+1) $ 
for some positive finite constants
$C_1$ and $C_2,$ independent of $x, n,$ and $\bfl.$  
\end{proposition}

\begin{proof}
Fix $\bfl \in \cK_I.$ Suppose that $\lambda_i >0,$ thus,
$\beta_i(\bfl) >0,$ as otherwise 
(\ref{recurrence}) is swiftly verified.
Assume that  $x = x_1 x_2 \ldots x_{n-1} i  \in \cG_n $ for some 
$ i \in \cA.$ Identity (\ref{recursion}) may be rewritten as 
\begin{equation}
    \label{quotident}
    I_x / \overline{u}_{x i^{-1}} = 
          \lambda_i + \sum_{j \in \cA - \{i \}}
           \lambda_j \overline{u}_{xj}/ \overline{u}_{x i^{-1}},
\end{equation}
which in light of (\ref{chainineq}), the definition
of the matrix $B_1(\bfl),$ identity (\ref{equalbmatrices})
and the remark following (\ref{equalbmatrices}) 
may be restated as
$$
  b_i(\bfl) / D_{\bfl}(x) = \lambda_i + 
                b_i(\bfl) \, \sum_{j \in \cA -\{ i\}}
                         \lambda_j  b_j (\bfl),  
$$              
where  $1/c_1(\epsilon) \leq D_{\bfl}(x) \leq 
           (n+1) c_2(\epsilon)/c_1(\epsilon)$ because,
by (\ref{chainineq}), 
$ c_1(\epsilon)/ ((n+1) c_2(\epsilon))  \leq I_x/\overline{u}_x 
        \leq c_1(\epsilon).$ Since the constants $c_1(\epsilon)$ 
and $c_2(\epsilon)$ can be bounded uniformly on $\cK_I,$ 
this ends the proof.
\end{proof}

Observe that the system of identities (\ref{recurrence})
may differ for different $ \bfl$ since the coefficients
$ D_{\bfl}(x)$ are not necessarily the same for different $\bfl.$
For fixed $\bfl$ and $x \in \cG,$  
those are $2d$ equations of which the $d$ equations indexed
by $j \in \cAm $ are redundant since $ \lambda_i = \lambda_{i^{-1}}$
and $ b_i = b_{i^{-1}}.$ Hence, system (\ref{recurrence}) may be 
regarded as a system of $d$ equations in $2d$ variables
or in $d$ variables.

\begin{proposition}
 \label{equalbeta}
Let $\epsilon >0$ and assume that Hypothesis I holds. Then  
for every $\bfl \in \cK_I,$  
we have $\beta_i(\bfl) =
\beta_j(\bfl) $ for all $i,j \in \cA$ if and only if 
$ \lambda_i = \lambda_j $ for all $i,j \in \cA$ (that is,
the contact process is isotropic). Equivalently,
$\overb_i(\bfl) = \overb_j(\bfl) $ 
for all $i,j \in \cA$ if and only if 
$ \lambda_i = \lambda_j $ for all $i,j \in \cA.$ 
\end{proposition}

\begin{proof}
It is clear that the statement in terms of the $\beta_i$ is
equivalent to the statement in terms of the $\overb_i.$ 
Thus in this proof, we will
restrict our attention to the collection of $\overb_i.$

Pick $\bfl \in \cK_I$ and fix $x \in \cG.$ Write $D = D_{\bfl}(x)>0$
and $ C_{\bfl} = \sum_{j \in \cA} \lambda_j.$
First, suppose that $b_i(\bfl)= b(\bfl)$ for each $i \in \cA.$
Then in view of Proposition \ref{betarecursion}, 
relation (\ref{recurrence}) reads
$$
   b(\bfl) = \lambda_i  D
           -  \lambda_i D b(\bfl)^2
                 +  D b(\bfl)^2 \, \sum_{j \in \cA} \lambda_j  
        = \lambda_i D (1- b(\bfl)^2) + D b(\bfl)^2  C_{\bfl}  
$$
for every $i \in \cA.$ But, because this holds for all $i \in \cA,$
it must be the case that, for any two
distinct indices $i, j \in \cA,$
$$ 
   \lambda_i D (1- b(\bfl)^2) = \lambda_j D (1- b(\bfl)^2). 
$$
However, if for some pair $i,j \in \cA,$ it was true that
$\lambda_i   \not = \lambda_j ,$
then this would mean that
$b(\bfl)^2 = 1, $ thus, $b(\bfl)=1.$ But this contradicts
our assumption that $\bfl \in \cK.$ Therefore, we conclude
that all infection rates $\lambda_i$ are identical, thus,
the contact process is isotropic.

To show the reverse direction, suppose that the contact process
is isotropic, that is, $\lambda_i = \lambda_* >0$ for all $i \in \cA.$
Then  equation (\ref{recurrence}) assembles as 
$$ 
   b_i(\bfl)  =  \lambda_* D  - \lambda_* D b_i(\bfl)^2 +
                 \lambda_* D b_i(\bfl) \, \sum_{j \in \cA}
                      b_j(\bfl). 
$$
Equivalently, if we write $C_b =  \sum_{j \in \cA} b_j(\bfl),$ 
$$
     b_i (\bfl)^2 + b_i(\bfl) \, 
           \frac{1 - C_b \bfl_* D}{ \bfl_* D} -1 = 0.
$$
Thus, we find an explicit expression for $b_i(\bfl)$ by means
of the quadratic formula, which is the same for all $i \in \cA$
(because the coefficients of the polynomial are independent of $i$). 
Hence, $b_i(\bfl) =b(\bfl)$ for all $i \in \cA,$ as required.
\end{proof}

The following Proposition is crucial
to prove Theorems \ref{weaksurvival},  \ref{limitset}
\mbox{and \ref{halfboundary}.} 

\begin{proposition}
 \label{secondtrans}
There exists a set $\cD_c \subset \partial \cK$
such that, for each $\bfl \in \cD_c,$ 
the $\{ b_i(\bfl)= b_i(\bfl,x) \}_{i \in \cA}$ satisfy 
\begin{equation}  
  \label{singidentity}
       \sum_{i \in \cA} \, \frac{b_i(\bfl)^2}{
     1 + b_i(\bfl)^2} = 1.
\end{equation}
\end{proposition}

\begin{proof}
In \cite{lahu}, it was shown that, for nonnegative constants 
$p_i < 1,$ not all zero, with all $p_{i}=p_{i^{-1}},$ 
the system  
\begin{equation}
    \label{brwsystem}
    F_i(z) = z p_i - z p_i F_i^2(z) + F_i(z) 
                    \sum_{j \in \cA} z p_j F_j(z)
\end{equation} 
of $2d$ equations is satisfied by {\em analytic} functions 
$\{ F_i(z) \}_{i \in \cA}$ within the largest circle not containing 
any singularities, and, that the system  
has a singularity for $z=R,$ where $R$ is the 
radius of convergence of the power series representation 
of $F_i$ for each $i \in \cA.$ Furthermore, each of the functions
$F_i(z)$ is an {\em algebraic} function, i.e.\ it satisfies a 
polynomial equation with coefficients in the ring ${\bf C}[z].$ 
(See \cite{lahu}, Section 3.2, for a discussion on an 
elimination algorithm and resultants to obtain $F_i.$)
But system (\ref{brwsystem}) is the same as system
(\ref{recurrence}), where now the $b_i(\bfl)$ play 
the role of the $F_i(z)$ and each $\lambda_j D_{\bfl} $ 
the one of $p_j z.$  Therefore, because there does not
exist any continuous choice for the functions $b_i(\bfl)$ 
at the singularity of the $F_i(z),$ it is clear that
if the system is singular for the $F_i,$ it must be singular
for the $b_i.$ Because the argument is an important piece
in our approach and is short, we will derive the equation that
characterizes the singularity.

Let $\cD_c$ denote the {\em set of singularities}, that is, the set 
of $\bfl$ such that system (\ref{recurrence}) 
is singular. If we use the short notation $D = D_{\bfl},$
we recall the equation 
\begin{equation}
   \label{shortrecurrence}    
    b_i(\bfl)  =   \lambda_i D + b_i(\bfl)^2 D \lambda_i
                  + 2 b_i(\bfl)  D  \, \sum_{j \in \cAp -\{ i\}}
                   \, \lambda_j \, b_j(\bfl). 
\end{equation}
The following derivation is in terms of the $b_i$ instead of the $F_i$
to avoid switching notation, i.e.\ set $b_i(\bfl)  = F_i(z).$ The
derivatives are symbolically justified by (\ref{brwsystem}).  
Next we rely on the Implicit Function theorem to locate the 
discontinuity.
If $ \{ H_i( \{b_k \}_{k \in \cAp}) \}_{i \in \cAp}$
denotes the difference between the lefthand and righthand sides of 
(\ref{shortrecurrence}), then the Jacobian matrix of the $d \times d$ system
of equations is given by $ ( dH_i / db_j )_{i,j \in \cAp},$
and whence, by the (complex) Implicit Function theorem
(see e.g.\ \cite{kend}), must be singular at every $\bfl \in \cD_c.$
Since
\begin{eqnarray}
   \label{betadiff}
    1 & = & d b_i / d b_i =  2 \lambda_i D b_i  +
               2  \sum_{k \in \cAp - \{i\}} \lambda_k D b_k  \\
    0 & = & d b_i / d b_j =  2 \lambda_j D b_i,  
       \hspace*{4.3cm} (j \not = i)
    \nonumber
\end{eqnarray}
the Jacobian matrix may be written as $I- J(\bfl),$
where $J$ has entries
\begin{eqnarray}
  \label{jmatrix}
      J (\bfl)_{ij} & = &  2 \sum_{k \in \cAp} \lambda_k D b_k 
               \qquad \qquad \mbox{if $j=i,$} \\
           & = & 2  \lambda_j D b_i  
 \hspace*{2.6cm} \mbox{if $j \not =i$.} \nonumber
\end{eqnarray}
Since these entries are nonnegative, the spectrum of $J(\bfl)$
is contained in the closed disk with radius $\gamma(\bfl),$
where $\gamma(\bfl)$ denotes the lead eigenvalue of $J(\bfl).$
In \cite{lahu}, we showed that the function $\gamma(\cdot)$ is
analytic, nondecreasing and continuous in its arguments.
For $\bfl$ with all $\lambda_i$ small,
the entries of $J(\bfl)$ are small since they are linear combinations
of the $\lambda_i,$ which implies that, for $\bfl$ sufficiently close to
the origin, $\gamma(\bfl) < 1.$ Therefore, the $ 
\bfl $ of smallest distance to the origin so that
$I - J(\bfl)$ is singular must be all those $\bfl,$ 
where $\gamma(\bfl)=1.$ In other words, 
  $
     \cD_c = \{ \bfl: \gamma(\bfl)=1 \}. 
  $ 
   
Next we analyze the corresponding eigenvalue equation
$ J(\bfl) v = v$ for $ \bfl \in \cD_c.$
As a consequence of the Perron-Frobenius theorem, the
vector $v$ has all entries nonnegative and at least one strictly
positive. Thus, set $  s = 2 \sum_{j \in \cAp} \lambda_j D v_j > 0.$
The eigenvalue equation $ J(\bfl) v = v$ together with equations
(\ref{jmatrix}) may be rewritten as
\begin{eqnarray*}
  v_i & = & v_i [ 2 \sum_{j \in \cAp } \lambda_j D b_j ]
          + 2 b_i \sum_{j \in \cAp - \{i\}} \lambda_j D v_j \\
     & = & v_i [ 2 \sum_{j \in \cAp } \lambda_j D b_j
            - 2  \lambda_i D b_i ] 
            + 2 b_i \sum_{j \in \cAp} \lambda_j D v_j \\
      & = &  2 v_i \sum_{j \in \cAp - \{ i \}} \lambda_j D b_j
            + b_i s.
\end{eqnarray*}
Multiplying each side by $ b_i $ and substituting the arrangement
of relation (\ref{recurrence}), i.e.\ 
$$
  b_i - \lambda_i D  - \lambda_i D b_i^2 =  
      2 b_i \, \sum_{j \in \cAp - \{i\}} 
                   \, \lambda_j \,  D \, b_j,
$$
gives
\begin{eqnarray*}
   v_i b_i & = & v_i [ b_i - \lambda_i D - 
       \lambda_i D b_i^2]  + b_i^2  s  \\
        \Longleftrightarrow \qquad 
        v_i D \lambda_i & = & s b_i^2 /( 1+ b_i^2).
\end{eqnarray*}
Recalling that $s =  2 \sum_{j \in \cAp} \lambda_j D v_j$ 
together with the symmetry of the contact process immediately 
provides relation (\ref{singidentity}). 

It remains to be shown that $\cD_c \subset \partial \cK.$
Observe that, in view of the continuity 
(Proposition \ref{contpotential}) 
and monotonicity 
properties of the functions $b_i,$ the 
expression in (\ref{singidentity}) {\em uniquely}
determines the singularity in the reals   
and the set $\cD_c$ cannot lie inside the region $\cK.$ 
On the one hand, if $d=1,$ equation (\ref{singidentity})
says that $b_i(\bfl) = \beta_i(\bfl) = \beta(\bfl)=1,$ thus,
obviously $\bfl \in \partial \cK.$ On the other hand, if
$d>1,$ first note that, for $\bfl \in \cK^c,$ there exists
some $a \in \cAp$ with $\beta_a(\bfl)=1, $ 
thus by Lemma \ref{bilessthanone}, with $b_a(\bfl) =1.$    
Now choose some $\bfl \in \mbox{int}(\cK).$ Since $d >1,$
we can assume that $b_k(\bfl), b_a(\bfl) >0$ 
for some $a,k \in \cA.$ Pick such $a$ and $k$ in $\cA.$
Suppose that $2 b_k(\bfl)^2 / (1 + b_k(\bfl)^2) = \delta >0.$
Therefore, due to the monotonicity of the function
$b_k(\bfl),$ for each $\bfl \in \cD_c,$
\begin{eqnarray*}
        1  =  \sum_{i \in \cA} \, \frac{b_i(\bfl)^2}{
     1 + b_i(\bfl)^2} & \geq & 2 \frac{b_a(\bfl)^2}{
                          1 + b_a(\bfl)^2} + \delta \\
   \Longrightarrow \qquad 
        1  > (\frac{1 - \delta}{1 + \delta})^{1/2} & \geq & b_a(\bfl).                  
\end{eqnarray*}   
By symmetry, $b_k(\bfl) < 1,$ too.
As a consequence, since $a$ and $k$ in $\cA$ were arbitrary,
$\bfl \in \cK.$ Combining these gives $\cD_c \subset
\cK$ and $\cD_c \subset \partial \cK.$  Also note that each
$b_a(\cdot)$ experiences a jump discontinuity for each $\bfl \in
\cD_c.$  
This completes our proof.
\end{proof} 

We point out that, perhaps surprisingly, identity 
(\ref{singidentity}) does not involve the $D_{\bfl}$
that appear in relation (\ref{recurrence}). In other words,
the form of equation (\ref{singidentity}) is independent of 
the $D_{\bfl}.$ Recall our strategy of proof 
from Section 2.5. We now carry out the last few steps. 
The following observation continues the discussion
in the last part of the proof of 
\mbox{Proposition \ref{secondtrans}.}
\begin{eqnarray*}
   \label{characboundary}
  \cD_c \subset \cK  & &  \mbox{ if and only if } \, 
    \partial \cK \subset \cK, \\
\mbox{equivalently,} \qquad \qquad \qquad & &  \\
  \cD_c \subset \cK  & &  \mbox{ if and only if each } \, b_i(\bfl) < 1
\mbox{ for } \bfl \in \cD_c.
\end{eqnarray*} 
We conclude that $\cD_c \subset \cK $ for $d>1,$ and, that
                 $\cD_c \subset \cK^c$ for $d=1$ because
                  $\bfl \in \cD_c$ implies $\bfl \not \in \cK.$
We end this section by collecting some consequences of
the preceding result. 

\begin{corollary}
   \label{valuesattrans}
For $\bfl \in \cD_c,$ 
\begin{eqnarray*}
  \exp  \overline{\varphi}_{2;\bfl}  & = & \theta_2(\bfl)=1, \\
     \eta(\bfl) & = & 1,
\end{eqnarray*} 
and for $\bfl \in \mbox{int}(\cK),$ 
\begin{eqnarray*}
  \exp  \overline{\varphi}_{2;\bfl}  & = & \theta_2(\bfl) < 1, \\
     \eta(\bfl) & < & 1.
\end{eqnarray*} 
For each $j \in \cA,$ the functions
$ \lambda_j \rightarrow \overline{\varphi}_{2;\bfl}$ and
$  \lambda_j \rightarrow \eta(\bfl)$ are continuous for 
$\bfl  \in \mbox{int}(\cK),$ are left-continuous for $\bfl \in \cK$  
and $d>1,$ and, are right-continuous for $\bfl \in \cK^c$ and
$d=1.$ The functions
$ \lambda_j \rightarrow \overline{\varphi}_{2;\bfl}$ and
$  \lambda_j \rightarrow \eta(\bfl)$ are not right-continuous 
for $\bfl \in \cD_c$  
and $d>1$ but are left-continuous for $\bfl \in \cD_c$ and
$d=1.$ In addition, 
all above statements on (dis)-continuity are valid for
the functions $ \beta_j(\bfl)$  and
$\overb_j(\bfl)$ for $j \in \cA,$ and for 
 $ \overline{\varphi}_{\rho;\bfl}$ for each $\rho >r_u(\bfl).$
Furthermore, 
\begin{eqnarray*}
         r_u(\bfl) & < & 2 \qquad \qquad \mbox{for }
                     \bfl \in \mbox{int}(\cK), \\ 
         r_u(\bfl) & \leq & 2 \qquad \qquad \mbox{for }
                  \bfl \in \cD_c. 
\end{eqnarray*}                   
\end{corollary}

From this result, it follows that Hypothesis I is {\em automatically}
satisfied, and thus, superfluous whenever stated.

\begin{proof}
First, the claims about $r_u(\bfl)$ are obvious from the definition
of $r_u$ along with Propositions \ref{rubound}, \ref{potmonotone}
and \ref{secondtrans}.
 
From Proposition \ref{secondtrans} together with equation
(\ref{mequation}), it follows that  
$ \exp  \overline{\varphi}_{2;\bfl}  = \theta_2(\bfl)=1$ 
for $\bfl \in \cD_c.$ Proposition \ref{potmonotone} implies
that, for $\bfl \in \mbox{int}(\cK),$ we obtain
$ \exp  \overline{\varphi}_{2;\bfl}  <1,$ and thus
by Proposition \ref{etalessthanone},
$\eta(\bfl) < 1.$

Also, Proposition \ref{contpotential} provides that, 
for each $j \in \cA,$ the functions
$ \lambda_j \rightarrow \overline{\varphi}_{2;\bfl}$ and
$  \lambda_j \rightarrow \eta(\bfl)$ are continuous for 
$\bfl  \in \mbox{int}(\cK).$ Hence for $d>1,$
because $\cD_c \subset \cK,$ 
the function $\lambda_j \rightarrow 
\eta(\bfl)$ is left-continuous for $\bfl \in \cK,$ 
and so is $\overline{\varphi}_{2;\bfl}.$ 
In other words, $\eta(\bfl) <1$ for $\bfl \in \mbox{int}(\cK)$ 
and  $\eta(\cdot)$ is continuous and strictly
increasing for $\bfl  \in \mbox{int}(\cK)$ 
with $\eta(\bfl) \uparrow 1$ as $\bfl \in \cK$ approaches
 $\bfl  \in \cD_c$ for $d>1.$
For $d=1,$ due to monotonicity,
$\eta(\bfl)=1$ for $\bfl \in \cK^c,$ thus, 
$\eta(\cdot)$ is constant and   
certainly right-continuous. The statement about the discontinuity
is an instance of Proposition \ref{secondtrans}. 
The statements pertaining to $ \overline{\varphi}_{2;\bfl}$
carry over to $ \beta_j(\bfl)$  and
$\overb_j(\bfl)$ for $j \in \cA.$ 
Similar reasoning leads to the results for
$ \overline{\varphi}_{\rho;\bfl}$ for each $\rho >r_u(\bfl).$
\end{proof}

Note that, due to the fact that $\theta_2(\bfl) = 1$ for
$\bfl \in \cD_c,$ it follows that
$r_u(\bfl) \geq 2$ for $\bfl \in \cD_c,$ thus, in view of
Corollary \ref{valuesattrans} and the left-continuity,
that $r_u(\bfl)=2$ for $\bfl \in \cD_c$ and $d>1.$


\subsection{Weak Survival Region} 
 
\begin{corollary}
  \label{eigenvaluelessthanone}
For $\bfl \in \overR_1 \cap \overline{\cR_1^c},$ we have
$\theta_1(\bfl) \leq 1.$
\end{corollary} 

\begin{proof}
For $d=1$ indeed, by Corollary \ref{valuesattrans},
$\theta_1(\bfl) = \theta_2(\bfl) =1$ for
$\bfl \in \cD_c,$ as desired.

Next consider the case $d>1.$
Note that $ \overR_1 \cap \overline{\cR_1^c} \subset \cK$
because, for $\bfl \in \cD_c,$ each $b_i(\bfl) < 1.$
Hence by Corollary \ref{valuesattrans},
 $ \theta_1(\bfl) $ is left-continuous
in each variable $\lambda_j$ for 
$\bfl \in \overR_1 \cap \overline{\cR_1^c}.$ Now Lemma
\ref{meanone} proves our claim for $d>1.$
\end{proof}

\begin{corollary} 
   \label{weaksurvivalexists}
The weak survival region $\cR_2$ has nonempty interior for $d>1.$
\end{corollary} 

\begin{proof}
This is an immediate consequence of the following facts.
Observe as in the previous proof that 
$ \overR_1 \cap \overline{\cR_1^c} \subset \cK.$ 
By Corollary \ref{eigenvaluelessthanone}, $\theta_1(\bfl) \leq 1$
for $\bfl \in \overR_1 \cap \overline{\cR_1^c}.$
By Corollary \ref{valuesattrans}, $\theta_2(\bfl)=1$  
for $\bfl \in \cD_c = \overline{\cD_c}.$
In view of Proposition \ref{potmonotone},
for $\bfl \in \overR_1 \cap \overline{\cR_1^c},$ 
$$ 
  \theta_2(\bfl) < \theta_1(\bfl) \leq 1.
$$ 
Consequently, $\cD_c \cap (\overR_1 \cap \overline{\cR_1^c}) =
\emptyset.$
Furthermore, because $\theta_2(\bfl)$ is continuous in each
variable for $\bfl \in \mbox{int}(\cK),$ left-continuous
for $\bfl \in \overline{\cK},$ and strictly increasing along directions
of increase on
$ \cK,$ we conclude that the complement of $\cR_1$ in 
$\cK$ has nonempty interior.
 
Combining these with Proposition \ref{secondtrans} gives
$$ 
 \cK^c = \cR_3
$$
and $\cD_c =\overline{\cR_3^c} \cap \overR_3.$ This implies
that $\cK \setminus \cR_1 = \cR_2,$ with $\cR_2$ having
nonempty interior.
\end{proof}
  
The discussion of the last proof has the following  
corollary.

\begin{corollary}
 \label{kcequalss}
For $d>1,$ we have
\begin{equation}
  \label{kidentified}
    \cK^c = \cR_3.
\end{equation}
\end{corollary} 

It remains to see what happens for $d=1.$
Combining Proposition \ref{secondtrans} with the left-continuity 
and strict monotonicity of $\beta(\lambda)$ implies that $\beta$
strictly increases in a continuous fashion all the way up to take 
the value $1$ at $\cD_c.$ Hence, the function $\beta(\lambda)$ is
continuous on $(0, \infty).$
Since the contact process does not survive
as long as $\beta <1,$ it follows that $\cR_3 \subset \cK^c.$

We return to more than one dimension.
For $d>1,$ the complement of $\cK$ coincides
with the strong survival region $\cR_3$ and $\cD_c$ coincides
with $\overR_2 \cap \overR_3.$ 
As a consequence, $\overR_2 \cap \overR_3 \subset \cR_2,$
which leads to the following observation.

\begin{corollary}
  \label{cpatsecondphasetrans}
The symmetric anisotropic contact process on $\cT_{2d}$ for $d > 1$ 
with $\bfl \in \overR_2 \cap \overR_3$ survives {\em weakly}.
\end{corollary} 

This proves Theorem \ref{survivalattransition}.
Replacing $\cK$ by $\cR_1 \cup \cR_2$ in all our previous results
completes the proof of Theorem \ref{properties} and combining
with Corollary \ref{valuesattrans} accomplishes Theorem 
\ref{etaatboundary}.
 
\begin{corollary}
   \label{propertiesinr}
The statements of Proposition \ref{potmonotone} hold
for each $\bfl \in \cR_1 \cup \cR_2$ with each $\lambda_k >0$
and the statements in Propositions \ref{contpotential} 
and \ref{continuity} hold for each 
$\bfl \in \mbox{int}(\cR_1 \cup \cR_2).$ Furthermore for $d>1,$ 
if continuity is replaced by left-continuity, then 
the statements in Propositions \ref{contpotential} 
and \ref{continuity} hold for each 
$\bfl \in \cR_1 \cup \cR_2$ with each $\lambda_k >0.$ 
\end{corollary}

\begin{proof}
It remains to see the last claim. But this is a consequence
of the fact that $\eta(\bfl) <1$ for $\bfl \not \in \overR_3,$ 
and thus, $\eta(\cdot)$ is continuous and strictly
increasing for $\bfl  \in \mbox{int}(\cR_1 \cup \cR_2)$ 
 with $\eta(\bfl) \uparrow 1$ as $\bfl \in \cR_2$ approaches
 $\overR_2 \cap \overR_3$ for $d>1.$ Hence, for each
$j \in \cA,$ the function $\lambda_j \rightarrow 
\eta(\bfl)$ is left-continuous for $\bfl \in \overR_2.$    
But this implies left-continuity for $\bfl \in \overR_2$  
for {\em all} functions that are continuous for 
$\bfl  \in \mbox{int}(\cR_1 \cup \cR_2).$  
\end{proof}

For $d=1,$ the functions described in
Corollary \ref{propertiesinr} are constant on $\overR_3,$
thus certainly right-continuous for $\bfl \in \overR_3.$
 
\begin{corollary}
  \label{criticalexponents}
For $d>1,$ we have 
\begin{eqnarray*}
     r_u(\bfl) & = & 1 \qquad \qquad \mbox{for }
            \bfl \in \overR_1 \cap \overR_2, \\
      r_u(\bfl) & = & 2 \qquad \qquad \mbox{for }
            \bfl \in \overR_2 \cap \overR_3. 
\end{eqnarray*}     
Moreover, for every $\bfl \in \overR_1 \cap \overR_2,$       
\begin{eqnarray}
    \label{thetaoneone}
    \theta_1(\bfl) & = &1 .
\end{eqnarray}   
\end{corollary}

\begin{proof}
Indeed, we have already noticed that $r_u(\bfl)=2$ for $\bfl \in 
\cD_c =\overR_2 \cap \overR_3$ for $d>1$ (remark following 
Corollary \ref{valuesattrans}).
As observed earlier, $r_u(\bfl) \geq 1$ for $\bfl \not \in \cR_1.$
By \mbox{Lemma \ref{meanone},} we have $r_u(\bfl) \leq 1 $ for 
            $\bfl \in \overR_1 \cap \overR_2.$ All claims now follow
from continuity and strict monotonicity considerations.
\end{proof}

Next we complete the proof of  Theorem \ref{weaksurvival}.

\begin{proposition}
 \label{weakregion}
The weak survival region $\cR_2$ 
enjoys the following properties:
\begin{enumerate}
\item[(a)]
\vspace*{-0.2cm}
The boundary $ \overline{\cR}_1 \cap \overline{\cR}_2$ consists
of all $\bfl$ such that the $\overb_i(\bfl)$ satisfy
$$
       \sum_{i \in \cA} \, \frac{ \overb_i(\bfl)}{
     1 + \overb_i(\bfl)} = 1
$$
and the boundary $ \overline{\cR}_2 \cap \overline{\cR}_3$
consists of all $\bfl$ such that the $\overb_i(\bfl)$ satisfy
$$
       \sum_{i \in \cA} \, \frac{ \overb_i(\bfl)^2}{
     1 + \overb_i(\bfl)^2} = 1.
$$ 
\item[(b)]
\vspace*{-0.2cm}
Every line in the interior of the first quadrant in 
${\bf R}^d$ that passes through the origin has an 
intersection with $\cR_2$ that is a line segment. 
\item[(c)]
\vspace*{-0.2cm}
The region $\cR_2$ has positive $d$-dimensional Lebesgue measure.
\item[(d)]
\vspace*{-0.2cm}
The region $\cR_2$ is connected and is a symmetric region 
in the $d$ parameters $\lambda_i.$
\end{enumerate} 
\end{proposition}

\begin{proof}
First we list some observations. Let $d >1.$
Define the functions $R(x) = x/(1+x)$ and  
$Q_{\rho}(\bfl) = \sum_{i \in \cA} R(\overb_i(\bfl)^{\rho}).$
The function $R(x)$ is continuously differentiable and has first
derivative $ 1/(1+x)^2 > 0$ for nonnegative $x,$
thus, is strictly increasing in $x.$
Recall from (\ref{thetaoneone}) that $\theta_1(\bfl)=1$
for $\bfl \in \overR_1 \cap \overR_2$ and from
Proposition \ref{secondtrans} that $\theta_2(\bfl)=1$
for $\bfl \in \overR_2 \cap \overR_3.$ 

{\bf (a).} 
Whence by (\ref{mequation}),
for every $\bfl \in \overR_1 \cap \overR_2,$
$$ 
   Q_1(\bfl)=1,
$$
and for every $\bfl \in \overR_2 \cap \overR_3,$
$$ 
  Q_2(\bfl)=1.
$$

{\bf (b)} and {\bf (c).} 
In fact, (c) will follow once (b) has been proved.
Fix $\bfl \in \overR_1 \cap \overR_2.$ Thus 
by Proposition \ref{potmonotone}, $ Q_2(\bfl) <1.$ For
each $\alpha = \{ \alpha_k \}_{k \in \cA }$ with at least one
$\alpha_k >0,$ due to the continuity properties 
of the functions $\overb_i(\cdot),$ there exists a smallest
real $t_* > 0$ such
that
$$ 
   Q_2( \{ \lambda_j + \alpha_j t_* \}_{j \in \cA} ) = 1 .
$$
The positivity of $t_*$ now guarantees claims (b) and (c).
 
{\bf (d).}
The first portion of (d) is an immediate consequence of the
statement in (b), whereas the second portion follows from
the symmetric role that is being played by each
function $R(\overb_i(\bfl)^{\rho}),$ $i \in \cA,$
in defining $Q_1$ and $Q_2$ together with the statement in (a).
\end{proof}

The question whether, for $d>1,$
the critical contact process at the
first phase transition behaves as in the lower or upper
phase is more subtle than the analogous question at the
second phase transition. At the second phase transition,
the discontinuity settles the issue since the phases are
well separated in a certain sense, whereas at the first
phase transition the contact process behaves continuously.
However, we exploit the fact that for $\bfl \not \in \cR_1,$
the expected value of $E \vert \cY_m \vert$ must grow without
bound as $m \rightarrow \infty,$
and thus, a Galton-Watson tree may
be embedded in the set of vertices ever to be infected that
is dominated by the contact process and whose attached 
Galton-Watson process is supercritical.

The next result will prove Theorem \ref{extinctionattransition}. 

\begin{proposition}
  \label{firsttransition}
$ \overR_1 \cap \overR_2 \subset \cR_1,$ that is,
for $\bfl \in \overR_1 \cap \overR_2,$ the contact process 
on $ \cT_{2d}$ for $d >1$ almost surely becomes extinct.
\end{proposition} 
 
\begin{proof}
Let $d>1.$ Assume that $\lambda_a \geq \lambda_b \geq \lambda_i$ 
for all $i \in \cA.$ Thus, $\lambda_b >0$ and $\beta_b(\bfl)>0.$
In view of (\ref{wumeans}), Lemma \ref{meanfromroot}
and Corollary \ref{almostmean}, we obtain
 $ \sum_{x \in \cL_m^*} w_{x} \rightarrow \infty$ as $m \rightarrow
 \infty$ if and only if 
 $ \sum_{x \in \cG_m} u_{x} \rightarrow \infty$ as $m \rightarrow
 \infty.$ 

It suffices to show that the Galton-Watson trees $\tilde{\tau}_r$
with mean offspring
numbers $ \sum_{x \in \cG_r} u_x,$  constructed
in parallel with the Galton-Watson trees $\tau_r,$ have  
corresponding {\em subcritical} 
Galton-Watson processes $ \{ \tilde{Z}_n(r) \}_{n>0}$
for $\bfl \in \overR_1 \cap \overR_2.$
Indeed, observe that the Galton-Watson process
associated with $\tilde{\tau}_r$ dominates 
the set of infected vertices,
in particular, on the event of survival,
the limit set of $\tilde{\tau}_r$ contains the limit
set $\Lambda$ of the contact process. Thus, if the Galton-Watson
process dies out, then the limit set $\Lambda$ will be empty and
the contact process cannot survive.

The infection is most likely to follow trails along vertices with
addresses composed of mostly letters $a$ and $b.$
For every integer $m,$ 
for each real $0<s<1,$ define $H_s \subset \cG_m$ to be the
set of all vertices at distance $m$ from the root vertex 
which contain $sm$ letters $b$ and $(1-s)m$ letters $a.$
We first determine the cardinality $\vert H_s \vert.$
Indeed, note that $  \vert H_s \vert =
{m \choose ms}$ and expand by means of the Stirling formula for factorials,
namely, $m ! = (m/e)^m \sqrt{2 \pi m} e^c$ with 
$ (12m + 2/5m)^{-1} \leq c \leq (12m)^{-1},$ to get
\begin{eqnarray*}
    \vert H_s \vert & = & {m \choose ms}  
                      =  \frac{e^{-c}}{ \sqrt{2 \pi s(1-s) m}}
                        \, [ \frac{1}{s^s (1-s)^{1-s}} ]^m.
\end{eqnarray*}
Thus, by the mean value theorem
along with the fact that each $u_x(\bfl)$ decays
exponentially in $\vert x \vert$ because $\bfl \not \in \cR_3,$
there is some $0 < u_*(\bfl,s) < 1$ such that
\begin{eqnarray} 
   \label{hscard} 
   \sum_{x \in  H_s}  u_x(\bfl)  
              & =  & \frac{e^{-c}}
               { \sqrt{2 \pi s(1-s) m}}
                \, \left[ \frac{u_*(\bfl,s) }
                {s^s (1-s)^{1-s}} \right]^m.
\end{eqnarray} 
Set $\alpha(\bfl,s) = u_*(\bfl,s)/ (s^s (1-s)^{1-s}).$ 

Now let $\bfl \in \overR_1 \cap \overR_2.$
If $ \sum_{x \in \cG_m} u_{x}(\bfl) \rightarrow \infty$ as
$m \rightarrow \infty$ (so that $m$ is running
through powers of $2$), then it must be the case that there
exists some $s >0$ such that $\sum_{x \in  H_s}  u_x(\bfl) 
 \rightarrow \infty$ as $m \rightarrow \infty.$ In that case,
we must have $\alpha(\bfl,s) > 1.$ By (\ref{ubarineq}),
the continuity and
strict monotonicity properties discussed earlier, there
would exist some $\bfl' \in \mbox{int}(\cR_1)$ 
in a direction of decrease for $\bfl$ such that
$$ 
   1 < \alpha(\bfl',s) <  \alpha(\bfl,s).  
$$
Hence, it would follow that
  $ \sum_{x \in \cG_m} u_{x}(\bfl') \geq  \sum_{x \in  H_s}  u_x(\bfl')  
  \rightarrow \infty$ as $m \rightarrow \infty,$
equivalently by Proposition \ref{wmeanlimit},
$ \sum_{x \in \cL_m^*} w_{x}(\bfl')  
  \rightarrow \infty$ as $m \rightarrow \infty.$ 
By Lemma \ref{infmeansupercrit}, for sufficiently large $m>0,$
the Galton-Watson process $Z_n(m)$ attached to the
 Galton-Watson tree $\tau_m$ would then be supercritical, and thus, 
the Galton-Watson tree $\tau_m$ would be infinite with 
positive probability. Since this Galton-Watson tree is dominated
by the set of vertices ever to be infected, this contradicts $\bfl'
\in \cR_1.$ Hence, we conclude that  $\alpha(\bfl,s) \leq 1$  
for all $0<s<1.$ 

Finally, since $\alpha(\bfl,s) \leq 1,$  
an elementary calculation together with (\ref{hscard}) yields
that, for all sufficiently large $m,$ the sum
 $\sum_{x \in H_s} u_x \leq  e^{-c} / \sqrt{2 \pi s(1-s) m},$
which $\rightarrow 0$ as $m \rightarrow \infty.$ Therefore, 
for each fixed $1> s>0$ and all sufficiently large $m,$
each sub-Galton-Watson tree, restricted to the vertices in $H_s,$
has attached a {\em subcritical} Galton-Watson process, thus,
survives with probability zero.
This implies that the full Galton-Watson tree $\tilde{\tau}_m,$ 
which is dominating the contact process, has an associated 
{\em subcritical} Galton-Watson process $ \tilde{Z}_n(m).$ 
Consequently, the contact process cannot survive
with positive probability. This completes the proof.
\end{proof}

Finally, recall that $\mu_n(V) = P \{ (\omega_n , \omega_{n+1}, \ldots) 
\in V \} $   for every $ \omega = \omega_1 \omega_2  \ldots 
 \in \Omega$ and any Borel set $V \subset \Omega$ and that the function
$\varphi : \Sigma \rightarrow \bR$ was defined by
 $ \varphi (\ldots x_1 x_2 \ldots ) = \log b_{x_1}. $ 
Then the following result is an immediate consequence of 
Proposition \ref{pstatprocess} and Corollary
\ref{kcequalss}.

\begin{theorem}
   \label{statprocess}
For $\bfl \in \cR_2$
and every $n \geq 1,$ the measure $\mu_n$ 
is absolutely continuous with respect to $\mu_{\varphi}$ and
$\mu_n \stackrel{\cD} {\rightarrow} \mu_{ \varphi}$
as $n \rightarrow \infty.$ Furthermore, the pressure $P (\varphi) =0.$ 
\end{theorem}


\section{Upper Bounds}
\setcounter{equation}{0}

Now we are ready to prove the easier bounds for the Hausdorff 
dimensions. Recall from the beginning paragraph in Section 4
that $\theta_1 = \theta$ is the lead eigenvalue of $M_1.$

\medskip
{\bf Upper bound for the Hausdorff dimension $\delta(\bfl)$ 
of the limit set.}

\begin{lemma}
  \label{firstupper}
With probability $1$ on the event of survival, $\delta(\bfl) \leq
   - \log \theta_1/\log \alpha .$
\end{lemma}
 
\begin{proof}
Recall that $\cY_m$ denotes the set of all vertices $x \in \cG_m$
that are ever infected. Pick $\varepsilon >0.$
Then by Lemma \ref{meanfromroot}, 
for all sufficiently large integers $m,$ 
\begin{eqnarray*}
    E \vert \cY_m \vert 
                  & = & \sum_{x \in \cG_m} u_x \\
                  &  \leq & \theta_1(\bfl)^{(1+ \varepsilon)m} .
\end{eqnarray*}
It follows from the Chebyshev-Markov inequality and
the Borel-Cantelli lemma that almost surely, eventually
\begin{equation}
  \label{uppercard}
   \vert \cY_m \vert \leq (\theta_1(\bfl)^{1 + \varepsilon} +
    \varepsilon)^m.
\end{equation}
         
Observe that the sets $\cY_m$ provide a sequence of open
covers of $\Lambda,$ in particular, if $\cE_x$ denotes the
set of all ends of $\cT$ that pass through $x,$ then
$$
    \Lambda \subset \cup_{x \in \cY_m} \cE_x.
$$
For each $x \in \cY_m,$ the diameter of $\cE_x$ 
(in the $d_{\alpha}$ metric) is $\alpha^m.$ Hence,
by (\ref{uppercard}) for all sufficiently large $m,$
$$
     \sum_{x \in \cY_m} 
   \mbox{diameter}_{\alpha}(\cE_x)^{ - \log(
   \theta_1(\bfl)^{1 + \varepsilon} +
    2 \varepsilon)/ \log \alpha} \leq 1
$$
(see e.g.\ \cite{falc}).
Because $\varepsilon > 0$ was arbitrary,
this implies that with probability $1$ on the event of survival,
the Hausdorff dimension of $\Lambda$ is
$$ 
   \delta(\bfl) \leq \frac{-\log \theta_1}{\log \alpha},
$$
as required.
\end{proof}

\smallskip
{\bf Upper bound for $\delta(\bfl;\mu).$}
Let $\mu$ be an ergodic, $\sigma$--invariant probability 
measure on the space $\Omega$ of semi-infinite reduced words
from $\cA.$ Recall that $\delta(\bfl; \mu)$ denotes the 
Hausdorff dimension of $\Lambda \cap \Omega_{\mu}$
(in the metric $d_{\alpha}$),  where $\Omega_{\mu}$ is the
subset of $\Omega$ consisting of all sequences $\omega$ that
are ``generic" for $\mu$ in the sense of definition 
(\ref{mugeneric}). Recall that $\varphi_{\bfl} : \Omega 
\rightarrow {\bf R}$ is the function defined by 
$\varphi_{\bfl}(x_1 x_2 \ldots) = \log \overb_{x_1}(\bfl).$
Since $\varphi_{\bfl}$ is continuous on $\Omega, $ relation
(\ref{mugeneric}) holds with $f = \varphi.$

\begin{lemma}
  \label{corrfreq}
For every $\varepsilon > 0,$ there exist sets $\Gamma_m =
\Gamma_m (\mu) \subset \cG_m$ of vertices at distance $m$ from
the root $1$ such that
\begin{equation}
    \label{cardgamma}
      \lim_{m \rightarrow \infty} \frac{1}{m} \log \vert 
      \Gamma_m \vert \leq h(\mu) + \varepsilon,
\end{equation}
\begin{equation}
  \label{closegamma}
    \lim_{m \rightarrow \infty} \sup_{x_1 x_2 \ldots x_m \in 
             \Gamma_m} \vert \frac{1}{m} \sum_{j=1}^m 
                 \log \overb_{x_j}(\bfl) - \int_{\Omega} 
             \varphi_{\bfl} d\mu \vert \leq \varepsilon,
\end{equation}
and such that for every sequence $x_1 x_2 \ldots \in \Omega,$
\begin{equation}
   \label{iogamma}
     x_1 x_2 \ldots \in \Omega_{\mu} \, \Rightarrow  \, x_1 x_2 \ldots
    x_m \in \Gamma_m(\mu) \, \mbox{ infinitely often}.
\end{equation}
\end{lemma}

\begin{proof}
The Shannon-McMillan theorem and the ergodic theorem along with
the ergodicity of $\mu$ guarantee that there are sets
$ \Gamma_m \subset \cG_m$ such that
\begin{eqnarray}
   \label{fullmeas}
   \lim_{m \rightarrow \infty} \mu( \bigcup_{x \in \Gamma_m} 
           \Omega(x)) =1, \\
      \label{smmtheor} 
    \lim_{m \rightarrow \infty} \max_{x \in \Gamma_m} 
             \vert - \log \mu(\Omega(x))/m - h(\mu) \vert = 0,
\end{eqnarray}
and such that (\ref{closegamma}) holds.                       
A consequence of (\ref{smmtheor}) is that for every $\varepsilon >0,$
for sufficiently large $m,$ and every $x \in \Gamma_m ,$
\begin{eqnarray}
  \exp \{ -m (1+ \varepsilon) h(\mu) \} & \leq &
        \mu(\Omega(x)) \leq \exp \{ -m(1 - \varepsilon) h(\mu) \}, \\
\mbox{thus,} 
  \quad \qquad \qquad \qquad \qquad \qquad & & \nonumber \\      
     \vert \Gamma_m \vert \exp \{ -m (1+ \varepsilon) h(\mu) \} 
          & \leq & \mu(\bigcup_{x \in \Gamma_m} \Omega(x)) 
            \leq \vert \Gamma_m \vert
           \exp \{ -m (1- \varepsilon) h(\mu) \},
\end{eqnarray}            
which together with (\ref{fullmeas}) provides (\ref{cardgamma}).
It remains to verify that (\ref{iogamma}) then holds.
Pick a subsequence $\{ m_n \}_{n>0}$ of integers such
that 
$$
     \sum_{n=1}^{\infty} [1 - \mu(\bigcup_{x \in \Gamma_{m_n}}
               \Omega(x))]  < \infty.
$$
In view of the Borel-Cantelli lemma, for $\mu$-almost every
sequence $x=x_1 x_2 \ldots \in \Omega,$ for sufficiently large
$n,$ the initial segment
$x_1 x_2 \ldots x_{m_n}$ is contained in  $ \Gamma_{m_n}.$
This completes the proof.               
\end{proof}

Our next result proves half of (\ref{dimintersect}). 

\begin{proposition}
 \label{secondupper}
If $h(\mu) + \int \varphi_{\bfl} d \mu < 0,$ then with probability
one, $\Lambda \cap \Omega_{\mu} = 
\emptyset.$ If \linebreak
$h(\mu) + \int \varphi_{\bfl} d \mu \geq 0, $
then with probability one on the event of survival, 
\begin{equation}
  \label{uppertwo}
      \delta(\bfl; \mu) \leq - \frac{h(\mu) + \int_{\Omega} 
        \varphi_{\bfl} d \mu}{ \log \alpha}.
\end{equation}
\end{proposition}

\begin{proof}
Let $\Gamma_m = \Gamma_m(\mu)$ be as stated in Lemma 
\ref{corrfreq}. If a reduced semi-infinite word $\omega=
x_1 x_2 \ldots $ is an element of $\Lambda \cap \Omega_{\mu},$
it must be the case that $x_1 x_2 \ldots x_m \in \cY_m \cap
\Gamma_m$ for infinitely many integers $m.$ Define
$$
  \Lambda_m(\mu) = \{ \omega = x_1 x_2 \ldots \in \Lambda \cap
                     \Omega_{\mu} : \,  x_1 x_2 \ldots x_m \in
                     \cY_m \cap \Gamma_m \}.
$$
Consequently for each $m \geq 1,$
\begin{equation}
   \label{partcovering}
   \Lambda \cap \Omega_{\mu} \subset \bigcup_{n \geq m} \Lambda_n(\mu).
\end{equation}

Hence, the set $\bigcup_{n \geq m} \Lambda_n(\mu)$ is a covering
of $\Lambda_n(\mu)$ by sets of
diameter $\alpha^m.$ Therefore, in order to find 
an upper bound for the Hausdorff dimension of 
$\Lambda_n(\mu),$ it suffices to find an upper bound for the
cardinality of $\cY_m \cap \Gamma_m.$  By the same reasoning as
in the proof of Lemma \ref{firstupper}, for every $\varepsilon >0$
and sufficiently large $m,$
$$
  E \vert \cY_m \cap \Gamma_m \vert \leq  (1+ \varepsilon)^m \,
         \sum_{x_1 x_2 \ldots x_m \in \Gamma_m} \prod_{j=1}^m 
      \overb_{x_j}(\bfl) = (1+ \varepsilon)^m  \,
       \sum_{x \in \Gamma_m} \exp \{ \sum_{j=1}^m 
      \varphi_{\bfl}(\sigma^j x) \}.
$$
From inequality (\ref{cardgamma}), it follows that
$\vert \Gamma_m \vert \leq \exp \{ m (h(\mu) + \varepsilon) \}.$
Furthermore, by (\ref{closegamma}), for every word $x \in
\Gamma_m,$ $\sum_{j=1}^m \varphi_{\bfl}(\sigma^j x) \leq
m ( \int \varphi_{\bfl} d \mu + \varepsilon) .$ With this in
mind, we conclude that, for all sufficiently large $m,$ the
expected cardinality of $\cY_m \cap \Gamma_m$ is no greater
than $ \exp \{ m ( h(\mu) + \int \varphi_{\bfl} d \mu + 
3 \varepsilon) \}.$ The Borel-Cantelli lemma thus implies
that, with probability one, eventually
\begin{equation}
   \label{cardintersect}
      \vert \cY_m \cap \Gamma_m \vert \leq
        \exp \{ m ( h(\mu) + \int \varphi_{\bfl}
       d \mu +  4 \varepsilon) \}.
\end{equation}
In case $ h(\mu) + \int \varphi_{\bfl} d \mu +  4 \varepsilon
< 0,$ eventually, $\cY_m \cap \Gamma_m $ is empty. Thus, 
$\Lambda_m (\mu)$ must be empty, and so by 
(\ref{partcovering}), $\Lambda \cap \Omega_{\mu} = \emptyset.$
In the other case when $ h(\mu) + \int \varphi_{\bfl} d \mu +  
4 \varepsilon \geq 0,$ then as a consequence of 
inequality (\ref{cardintersect}), for every $n \geq 1$ and 
all sufficiently large $m,$ the set $\Lambda_n(\mu)$ is covered
by $  \exp \{ m ( h(\mu) + \int \varphi_{\bfl}
       d \mu +  5 \varepsilon ) \} $ sets of diameter $\alpha^m.$
Since $\varepsilon > 0$ can be chosen arbitrarily small,
$$
   \delta_H(\Lambda_n(\mu)) \leq - (h(\mu) +
   \int \varphi_{\bfl} d \mu ) / \log \alpha.
$$
Since $\Lambda \cap \Omega_{\mu} \subset \bigcup_{n \geq m} 
\Lambda_n(\mu),$
the required inequality (\ref{uppertwo}) follows.
\end{proof} 


\section{Lower Bounds}
\setcounter{equation}{0}

To verify formulae (\ref{dimlimit}) and 
(\ref{dimintersect}) for the Hausdorff dimensions
of the random sets $\Lambda$ and $\Lambda \cap \Omega_{\mu},$
we need to establish the lower bounds for the Hausdorff dimensions.
For this purpose, we shall again consider
Galton-Watson trees, embedded
in the set of vertices of $\cT$ that are ever infected, whose
limit sets are subsets contained in $\Lambda$ and $\Lambda \cap
\Omega_{\mu}$ that have Hausdorff dimensions which approach the
required bounds. In turn, these Hausdorff dimensions are
calculated by invoking a theorem of Hawkes \cite{hawk,lyon}
and an extension given in 
\cite{lase3}.
First, we show that the Hausdorff dimensions of both limit
sets are almost surely constant.

\begin{lemma}
  \label{constdimensions}
The Hausdorff dimensions $\delta(\bfl)$ and $\delta(\bfl; \mu)$ 
are almost surely constant on the event of survival.
\end{lemma}

\begin{proof}
The reasoning is parallel for each of both Hausdorff
dimensions. We shall give an argument for $\delta(\bfl),$ largely 
borrowed from \cite{lase2}.
Let $\delta_*$ be the essential supremum of the random variable
$\delta(\bfl).$ Then for any $\delta < \delta_*,$ there is
positive probability $p$ that the limit set of a contact process
initiated at the root $1$ has Hausdorff dimension at least 
$\delta.$ Since the distribution of the Hausdorff dimension
of a subset of the boundary $\partial \cT$ does not change
by an isometry of $\cT,$ due to the geometry of the tree,
it follows that for any vertex $x \in \cT,$ there is positive
probability $p$ that the limit set of a contact process 
initiated at $x$ has Hausdorff dimension at least $\delta.$ 
Therefore, if $\cF_t$ is the $\sigma$-algebra generated by
the percolation structure up to time $t,$ then
$$
   P \{ \delta(\bfl) \geq \delta \vert \cF_t \}  \geq
    p I_{\{ \vert A_t \vert \geq 1 \}}.
$$
But the martingale convergence theorem implies that 
this conditional probability converges to the indicator
function of the event $\{ \delta(\bfl) \geq \delta \}$
almost surely as $t \rightarrow \infty.$ Since, obviously,
the indicator of the event $\{ \vert A_t \vert \geq 1 \}$
converges to that of the event that the contact process
survives, it follows that 
$$
  I_{\{ \delta(\bfl) \geq \delta \}} \geq p     
    I_{\{ \mbox{survival}\}}
$$
almost surely. Since the indicators take only the values 
$0$ and $1,$ it immediately follows that $\delta(\bfl) \geq
\delta$ almost surely on the event of survival. Hence,
$\delta(\bfl) = \delta_* $ almost surely on the event of
survival, as desired.
\end{proof}


\subsection{Hawkes' Theorem and an Extension}

Recall the Galton-Watson trees $\tau$ from Section 2.4

\begin{theorem}[Hawkes \cite{hawk}]
If the offspring distribution has mean $\mu > 1$ and finite
second moment, then almost surely on the event of nonextinction,
the limit set $\Lambda_{GW}$ of the Galton-Watson tree $\tau$
has Hausdorff dimension (in metric $d_{\alpha}$)
$$
   \delta_H(\Lambda_{GW}) = - \frac{\log \mu}{\log \alpha} .
$$
\end{theorem}

Hawkes examines only the case $ \alpha = 1/2$ but the
result and its proof hold for all $\alpha \in (0,1)$ \cite{lyon,lase2}.
Now recall from Section 2.3 that the labelled Galton-Watson processes 
have label set $\cB.$
For each label $i \in \cB,$ let $q_i = \sum _{F \subset \cB; i 
\in F} Q(F) $ denote the probability that label $i$ is included
in a random set with distribution $Q.$ Define a function 
$ \psi: \cB^{\bf N} \rightarrow {\bf R}$ by 
\begin{equation}
     \label{littleq}
   \psi(x_1 x_2 \ldots) = \log q_{x_1}.
\end{equation}
For any ergodic, shift-invariant probability measure $\mu$ on the 
sequence space $\cB^{\bf N},$
define $\cB^{\bf N}_{\mu}$ to be the set of $\mu$-generic
sequences $\omega$ such that for every continuous function
$f: \cB^{\bf N} \rightarrow {\bf R},$
\begin{equation}
   \label{expmutwo}
       \lim_{n \rightarrow \infty} \frac{1}{n} \sum_{i=1}^n 
                 f(\sigma^i \omega) = \int_{\cB^{\bf N}} f d\mu.
\end{equation}

\begin{theorem}[Lalley and Sellke \cite{lase3}] Let $\tau$
be the labelled Galton-Watson tree attached to a supercritical 
labelled Galton-Watson process with label set $\cB$ and 
offspring distribution $Q,$ and let $\mu$ be any ergodic,
$\sigma$--invariant probability measure on $\cB^{\bf N}.$ 
If $h(\mu) + \int \psi d\mu < 0$ then with probability one,
\begin{equation}
     \partial \tau \cap \cB^{\bf N}_{\mu} = \emptyset.
\end{equation}             
If $h(\mu) + \int \psi d\mu \geq 0,$ then almost surely on
the event of nonextinction, the Hausdorff \linebreak
dimension of $\partial \tau \cap \cB^{\bf N}_{\mu} $ 
in the metric $d_{\alpha}$ is
\begin{equation}
    \delta_H(\partial \tau \cap \cB^{\bf N}_{\mu}) 
     =  - \frac{h(\mu) + \int \psi d\mu}{ \log \alpha}.
\end{equation}
\end{theorem}


\subsection{Lower bounds for $\delta(\bfl)$ and $\delta(\bfl; \mu)$}

The following corollary together with Lemma \ref{firstupper}
completes the proof of equation (\ref{dimlimit}).

\begin{corollary}
  \label{firstlower}
   With probability $1$ on the event of survival, 
  $\delta(\bfl) \geq
   - \log \theta_1(\bfl)/ \log \alpha  .$
\end{corollary}

\begin{proof}
The Hausdorff dimension $\delta_H(\Lambda)$ of $\Lambda$
is almost surely constant on the event of survival by Lemma 
\ref{constdimensions}, thus, it is enough to 
show that for any real number $\delta^* <  
 - \log \theta_1(\bfl) / \alpha, $ there is positive probability
that $ \delta_H(\Lambda) \geq \delta^*. $ Consider the embedded
Galton-Watson trees $\tau_k.$ Since $\bfl \in \cR_2,$ we know
that $\theta_1(\bfl) > 1.$ Indeed, if $\bfl$ is off the region
$\cR_3$ and in $\cR_2,$ then $\bfl \in \mbox{int}(\cR_2),$
because $\overR_1 \cap \overR_2 \subset \cR_1.$ Then
the strict monotonicity of $\theta_1(\cdot)$ along any direction of 
increase and the fact that $\theta_1(\bfl)=1$ for $\bfl \in \overR_1
\cap \overR_2$ imply that $\theta_1(\bfl) >1$ for $\bfl \in \cR_2.$ 
The second moment condition in the Hawkes' theorem
is satisfied since the offspring random variable is a sum of indicator
variables
$$
 \sum_{x \in \cL_k^*} I_{ \{ \exists \mbox{ downward infection
trail } \{ \mbox{root} \} \rightarrow \mbox{ x beginning at t=0}\}}
$$
for every $k.$ This sum of indicators obviously has finite variance
for each $k.$
By the Hawkes' theorem, on the set
of nonextinction, the limit set $\Lambda_k$ of $\tau_k$ has
Hausdorff dimension (in the metric $d_{\alpha}$)
$$
    \delta_k = - \frac{\log \mu_k^{1/k}}{\log \alpha} .
$$
Since the probability of nonextinction is positive, it follows
that, with positive probability, the set $\Lambda$ has a subset
$\Lambda_k$ of dimension $\delta_k.$ By (\ref{gwmeans}), for
sufficiently large $k,$ $\delta_k \geq \delta^*,$ thus, with
positive probability, the Hausdorff dimension of $\Lambda$ is
at least $\delta^*,$ as desired.
\end{proof}

Proposition \ref{secondupper} together 
with the following result finishes the
proof of equation (\ref{dimintersect}).

\begin{corollary}
  \label{secondlower}
 For any ergodic, shift-invariant probability measure $\mu$ on $\Omega,$
with probability one on the event of survival,
\begin{equation}
   \label{ineqdiminter}
     \delta(\bfl;\mu) \geq - (h(\mu) +
 \int \varphi_{\bfl} d \mu) / \log \alpha.
\end{equation}
\end{corollary}

\begin{proof}
The proof is nearly parallel to the one of Corollary \ref{firstlower}.
Let $\psi = \varphi$ in (\ref{littleq}). 
In the case that the righthand side of (\ref{ineqdiminter}) 
is negative, the inequality apparently holds.
Suppose otherwise. Consider the limit set $\Lambda_k$ of the
Galton-Watson tree $\tau_k.$ 
As in the previous proof, the underlying Galton-Watson process is 
supercritical and 
the second moment condition in the Extended
Hawkes' theorem is satisfied.
By the shift-invariance and ergodicity of $\mu$ and
the definition of $\varphi_{\bfl}$ together with
the Extended Hawkes' theorem, almost surely on
the event of survival, the intersection of $\Lambda_k$
with $\Omega_{\mu}$ has Hausdorff dimension (in the metric $d_{\alpha}$)
$$
  - \frac{h(\mu) + k^{-1} \int_{\Omega} 
  \log \prod_{j=1}^k \overline{b}_{x_j}^w d\mu(x)}{
           \log \alpha},
$$
where the $\overline{b}_{j}^w$ denote the entries of the matrix
$B_1^w$ associated with the functions $w_x$ (as opposed to
$B_1$ associated with the $u_x,$ that we were using before).
But by Proposition \ref{wmeanlimit}, the matrices $B_1^w$ and $B_1$ have
the same lead eigenvalue, and thus, we conclude that
$\overb_j =  \overline{b}_{j}^w$ for every $j \in \cA.$ 
Since $\Lambda_k \subset \Lambda,$ the last display is
a lower bound for the Hausdorff dimension of $\Lambda \cap \Omega_{\mu}.$
\end{proof}


\section{Backscattering Inequalities}
\setcounter{equation}{0}
  
Throughout this section, we only address the case $d>1.$
The term ``backscattering inequality'' was coined in \cite{lase1}
since the inequality $\delta(\bfl) \leq \delta_H(\Omega)/2$ may be
proved by a backscattering argument in different context. The following
argument is short and does not involve the backscattering idea. 
Another line of reasoning will be explained below that
is based on the Gibbs Variational Principle. 

\begin{proposition}
 \label{schwartz}
For all $\bfl \not \in \cR_3,$
\begin{equation}
  \theta_1 (\bfl)^2 \leq \theta_2(\bfl)(2d-1) \leq 2d-1,
\end{equation}
with strict inequality $\theta_1^2(\bfl) < 2d-1$ holding
except possibly for $\bfl \in
\overline{\cR}_2 \cap  \overline{\cR}_3.$
\end{proposition}

\begin{proof}
Recall that $\theta_1(\bfl)$ and $\theta_2(\bfl),$ respectively,
are the lead eigenvalues of the matrices 
$M_1(\bfl)$ and $M_2(\bfl),$ respectively. Thus, for $i=1,2,$
$ \theta_i(\bfl) = \lim_{m \rightarrow \infty} ({\bf 1}^t
M_i(\bfl)^m {\bf 1})^{1/m}, $ where ${\bf 1}$ denotes the 
$2d$-vector with all entries $1.$ As an appeal 
to the Cauchy-Schwarz inequality
\begin{eqnarray*}
       ({\bf 1}^t M_1(\bfl)^m {\bf 1})^{1/m} & = &
                 \left ( \sum_{i_1 i_2 \ldots i_m} 
                        \prod_{j=1}^m \overb_{i_j}(\bfl) 
                     \right )^{1/m} \\
                 & \leq & \left ( \sum_{i_1 i_2 \ldots i_m} 
                        \prod_{j=1}^m \overb_{i_j}(\bfl)^2 
                     \right )^{1/2m} \left (
            \sum_{i_1 i_2 \ldots i_m} 1 \right )^{1/2m} .
\end{eqnarray*}
The last line of the display tends to $\sqrt{\theta_2(\bfl)}
\sqrt{2d-1}$ as $m \rightarrow \infty$ because the number of
reduced words $i_1 i_2 \ldots i_m$ of length $m$ equals
$ (2d-1)^m.$ Taking the limit as $m \rightarrow \infty$ also
implies that $\theta_1(\bfl)^2 \leq \theta_2(\bfl) (2d-1).$
It follows from  Proposition \ref{secondtrans}, Corollary
\ref{kcequalss}, and the strict monotonicity of $\theta_2(\bfl)$
along any direction of increase that
$\theta_2(\bfl) \leq 1$ with strict inequality except for
$\bfl \in \overline{\cR}_2 \cap  \overline{\cR}_3.$
\end{proof}

\begin{corollary}
 \label{backscatter}
For all $\bfl \not \in \cR_3,$ with probability one,
\begin{equation}
  \label{possbackscatter}
  \delta (\bfl) \leq \frac{1}{2} \delta_H(\Omega),
\end{equation}
where strict inequality holds except possibly for 
$\bfl \in \overline{\cR}_2 \cap  \overline{\cR}_3.$
\end{corollary}

\begin{proof}
It is an easy exercise that the Hausdorff dimension of
$\Omega$ in the metric $d_{\alpha}$ equals \linebreak $ - \log (2d-1)/
\log \alpha. $ Set $\theta = \theta_1.$
Hence, the result follows from equation
(\ref{dimlimit}) and Pro\-po\-si\-tion \ref{schwartz}.
\end{proof}

\smallskip
{\bf Gibbs Variational Principle.}
Consider any ergodic, shift-invariant probability measure
$\mu $ on $\Omega$ and the measure-preserving system 
$(\Omega, \sigma, \mu).$ Let $h(\mu)$ denote the 
measure-theoretic entropy, let $\varphi: \Omega \rightarrow
{\bf R}$ be any H\"{o}lder continuous function, and let
$P(\varphi)$ denote the thermodynamic 
pressure of the potential function $\varphi,$
as described in Section 3.
The {\em Gibbs Variational Principle} \cite{bowe,pesi}
states that 
\begin{equation}
  \label{gibbs}  
  h(\mu) + \int_{\Omega} \varphi d\mu \leq  
         P (\varphi),
\end{equation} 
where the inequality is strict unless $\mu= \mu_{\varphi}$
is the Gibbs state for $\varphi.$
If we let $\varphi = \rho \varphi_{\bfl},$ where 
$\varphi_{\bfl}(x_1 x_2 \ldots )   = \log \overb_{x_1}(\bfl),$
the pressure functional is given by
\begin{equation}
  \label{pressure}
                P(\rho \varphi_{\bfl}) 
   = \lim_{m \rightarrow \infty} \frac{1}{m} \log 
           \left ( \sum_{i_1 i_2 \ldots i_m} \prod_{j=1}^m
           \overb_{i_j}(\bfl)^{\rho}  \right ) 
                = \log \theta (\rho; \bfl).
\end{equation}
Since the potential functions depend only on the first 
coordinate, the corresponding Gibbs states
$\mu_{\rho \varphi_{\bfl}}$ are the probability distributions
of the stationary Markov chains with transition probabilities
$$
   {\bf P}(\rho; \bfl)_{ij} = \frac{\overb_j(\bfl)^{\rho} v_j}{
          \theta(\rho;\bfl) v_i} \, ( 1 - \delta_i(j^{-1})),
$$
where $v$ is the lead right eigenvector of the Perron-Frobenius
matrix $M_{\rho}(\bfl)$ and $\delta_{.}(\cdot)$ denotes 
the Kronecker delta function. 
Observe that these Gibbs states coincide with the ones
discussed in Section 3. In particular, 
when $\rho=2$ and $  \bfl \in \overline{\cR}_2 \cap 
\overline{\cR}_3,$ this Gibbs state coincides with  
the measure $\mu_*,$ announced in Theorem \ref{halfboundary}
(and discussed after that result).

The following result now completes the proof of Theorem
\ref{halfboundary}.

\begin{proposition}
For every ergodic, $\sigma$--invariant probability measure
$\mu$ on $\Omega,$
\begin{equation}
  \label{secbackscatter}
    \delta(\bfl;\mu) \leq \frac{1}{2} \delta_H(\Omega_{\mu})
\end{equation}
almost surely, where the inequality is strict unless 
$  \bfl \in \overline{\cR}_2 \cap 
\overline{\cR}_3$ and $\mu = \mu_*.$
\end{proposition}

\begin{proof}
Let $\rho=2$ and $\varphi = 2 \varphi_{\bfl}.$  We combine 
(\ref{gibbs}) and (\ref{pressure}), divide both
sides of the inequality by $2,$ and recall formula 
(\ref{dimintersect}) for the Hausdorff dimension 
$\delta(\bfl;\mu)$ of $\Lambda \cap \Omega_{\mu},$ 
\begin{eqnarray*}
  h(\mu) + 2 \int_{\Omega} \varphi_{\bfl} d\mu & \leq & 
   \log \theta_2(\bfl)  \\
   \Leftrightarrow \quad
     h(\mu) + \int_{\Omega} \varphi_{\bfl} d\mu & \leq & 
   \frac{1}{2} ( h(\mu) + \log \theta_2(\bfl))   \\
    \Rightarrow \qquad  \qquad \quad
      \delta(\bfl; \mu) & \leq &- \frac{1}{2} (h(\mu) +
        \log \theta_2(\bfl)) / \log \alpha.
\end{eqnarray*}
Unless $\mu $ coincides with the Gibbs state
$\mu_{2 \varphi_{\bfl}} = \mu_* ,$ strict inequality holds. 
Recall again that $\theta_2(\bfl) < 1$ for $\bfl$ in the interior
of $\cR_2$ and $\theta_2(\bfl) =1$ for $\bfl \in \overR_2 \cap \overR_3.$
This brings along
$$
  \delta(\bfl; \mu) \leq - \frac{1}{2} h(\mu)/ \log \alpha,
$$
and strict inequality holds unless 
$\mu = \mu_*$ and $  \bfl \in \overline{\cR}_2 \cap 
\overline{\cR}_3.$ The statement now follows from
the observation that $\Omega_{\mu} $ has Hausdorff dimension,
in the metric $d_{\alpha},$ equal $-h(\mu)/ \log \alpha,$
thanks to a theorem of Billingsley.
\end{proof}

Along with Corollary \ref{backscatter}, the following proposition
ends the proof of Theorem \ref{limitset}.

\begin{proposition}
Strict inequality holds in (\ref{halfdim1})
for every $\bfl \in \overline{\cR}_2 \cap \overline{\cR}_3$
except when the contact process is isotropic, that is,
when all infection rates $\lambda_i =  \lambda.$
\end{proposition}

\begin{proof}
By Proposition \ref{schwartz}, it suffices to verify that
$$
   \theta_1(\bfl)^2 < \theta_2(\bfl) (2d-1) 
$$ 
for every $\bfl \in \overline{\cR}_2 \cap \overline{\cR}_3$
except when the contact process is isotropic. 
Let $\bfl \in \overline{\cR}_2 \cap \overline{\cR}_3.$
By applying formula (\ref{pressure}) and the Gibbs Variational
Principle to the potential functions
$\varphi  = \varphi_{\bfl}$ and
$\varphi = 2 \varphi_{\bfl},$ respectively, with Gibbs state
$\mu=\mu_1,$ we find by (\ref{gibbs})
\begin{eqnarray}
   2 h(\mu_1) + 2 \int \varphi_{\bfl} d\mu_1 & = &
        2 \log \theta_1(\bfl)    \\
    \label{twopotential}
     h(\mu_1) + 2 \int \varphi_{\bfl} d\mu_1 & \leq &
         \log \theta_2(\bfl),
\end{eqnarray}
with strict inequality in (\ref{twopotential}) except
when the Gibbs states $\mu_1$ and $\mu_2$ coincide. 
An instance of \cite{bowe}, Theorem 1.28, allows this 
to happen when
the difference of the potential functions $\varphi_{\bfl}$
and $2 \varphi_{\bfl} $ is (co-)homologous to a constant
function, equivalently, if and only if there exists
constants $\epsilon > 0 $ and $c>0$ such that
$$
   \epsilon <  \frac{1}{c^m} 
    \frac{\prod_{j=1}^m \overb_{i_j}(\bfl)^2}{
            \prod_{j=1}^m \overb_{i_j}(\bfl)} < \frac{1}{
              \epsilon}
$$
for every $m \geq 1$ and all finite reduced words $i_1 i_2
\ldots i_m  \in \cG_m.$ This is possible if and only if
all of the values $\overb_i, $ $i \in \cA,$ are identical.
However, by Proposition \ref{equalbeta}, this happens if
and only if the contact process is isotropic or each $\overb_i =1.$
But the latter would contradict that $\bfl \in \overR_2 \cap \overR_3,$
thus, in $\cR_2.$ 
\end{proof}


\section{Isotropic Contact Process}
\setcounter{equation}{0}

We first outline the existence proof for an immediate phase 
in the isotropic case prior to summarizing additional
consequences in the isotropic setting.
Indeed, we recall that $\lambda \rightarrow \beta(\cdot)$ is a
strictly increasing and continuous function as long as 
$\beta(\bfl) <1,$ left-continuous everywhere,
that by Corollary \ref{eigenvaluelessthanone},
$\beta(\lambda) \leq  1/(2d-1)$ 
for $\lambda \in \overline{\cR_1^c}
\cap \overline{\cR}_1,$ and by Proposition
\ref{secondtrans},
$\beta(\lambda) = 1/\sqrt{2d-1}$ for $\lambda \in \cD_c.$
Thus, there must exist an intermediate interval between
extinction and strong survival that has positive length.

\begin{corollary}
For $d>1,$ consider the {\em isotropic} contact process, 
that is, with $ \lambda_i \equiv \lambda,$ in the 
weak survival phase
$\lambda_1 < \lambda \leq \lambda_2.$
Let $\beta(\lambda) = \beta_i(\lambda).$ 
Then there is some $0 <   D_{\lambda} < \infty$ such that
$$ 
   1/(2d D_{\lambda_1}) < \lambda \leq 1/(2 \sqrt{2d-1}D_{\lambda_2}) ,
$$
$$
     \beta(\lambda) = \theta_1(\lambda)/(2d-1) \qquad \mbox{and} \quad 
        \, \delta(\lambda) = - \log [(2d-1) \beta(\lambda)] / \log \alpha,
$$
where 
$$
    \theta_1(\lambda) = 1/(2 \lambda D_{\lambda}) - \{ 
     1/(4 \lambda^2 D_{\lambda}^2) -
       (2d-1)   \}^{1/2}.
$$
In addition, 
$$
  \theta_1(\lambda) = (2d-1) \beta(\lambda), \; \; 
      \theta_2(\lambda) = (2d-1) \beta(\lambda)^2, \; \;
   \theta_2(\lambda) = \beta(\lambda) \theta_1(\lambda).
$$
Especially at the phase transitions, 
\begin{eqnarray*}
    \theta_1(\lambda_1) & = & 1,  \quad \qquad  \quad
         \beta(\lambda_1) = 1/(2d-1),  \qquad
                \   \\*[0.1cm]
    \theta_1(\lambda_2) &= & \sqrt{2d-1}, 
         \quad \beta(\lambda_2) = 1/\sqrt{2d-1},
           \nonumber  \\*[0.1cm]
    \theta_2(\lambda_1) & = & 1/\sqrt{2d-1}, \qquad \qquad \qquad \qquad  
                \nonumber   \\*[0.1cm]
    \theta_2(\lambda_2) &= & 1. \qquad \qquad   \qquad \qquad  
                   \nonumber
\end{eqnarray*}
\end{corollary}

\begin{proof}
Suppose that $d>1.$ (For $d=1,$ while equation (\ref{singidentity})
holds, relation (\ref{recurrence}) that will be relied on 
below is not valid since the set $\cD_c$ is not in $\cK.$) 
Note that continuity allows to include the first critical value
$\lambda= \lambda_1$ in our calculations below. 
Also, observe that, in the isotropic case $ \lambda_i \equiv \lambda,$
by Proposition \ref{equalbeta}, $\overb_i(\lambda) = \overb(\lambda)
= \beta(\lambda).$ 

First, it is an easy exercise to solve equation (\ref{mequation}) 
for $\beta,$ thus, $\beta = \theta_1 /(2d-1) .$ With this in mind,
relation (\ref{recurrence}) 
\begin{equation}
   \label{recurrisotropic}
   \beta = \lambda D_{\bfl} [ 1 + \beta^2 (2d-1) ]
\end{equation}
may be restated as
$$
  \theta_1/ (2d-1) = \lambda D_{\lambda} [ 1 + \theta_1^2/(2d-1)] ,
$$
where there are some positive finite constants $C_1$ and $C_2$
such that for each $\lambda_1 \leq \lambda \leq \lambda_2,$
$ C_1 \leq D_{\lambda} = D_{\lambda}(x) \leq C_2 (n+1)$
by Proposition \ref{betarecursion}. 
This quadratic equation in $\theta_1$ has two solutions 
$$
   \theta_{1(1,2)} =
   \frac{1}{2 \lambda D_{\lambda}} \pm 
         \sqrt{ \frac{1}{4 \lambda^2 D_{\lambda}^2}
                -  (2d-1)}, 
$$ 
the relevant root being the one with the negative sign in front of the
radical because $\beta$ is a nondecreasing function in $\lambda,$
so is $\theta_1.$ 

Next we find the endpoints of the interval 
$(\lambda_1, \lambda_2],$ the interval of the weak survival phase.
Since $\theta_1(\lambda_1) = 1,$ by Corollary \ref{criticalexponents},
solving (\ref{mequation})
yields $\beta = 1/(2d-1),$ and, solving (\ref{recurrisotropic})
gives $\lambda = \lambda_1 = 1/(2d D_{\lambda}) = 1/(2d D_{\lambda_1}) .$ 
Moreover, since $\theta_2(\lambda_2) = 1$ by Corollaries
\ref{valuesattrans} and \ref{kcequalss},
solving equation (\ref{mequation})
with $\rho =2 $ yields $\beta = 1/\sqrt{2d-1},$ and, again solving
(\ref{recurrisotropic}) provides $\lambda= \lambda_2 
= 1/(2 \sqrt{2d-1} D_{\lambda}) = 1/(2 \sqrt{2d-1} D_{\lambda_2}) .$
Additionally, from equation (\ref{mequation}) it follows that
$\theta_2 = \beta \theta_1,$ which equals $ \beta^2 (2d-1).$
This finishes our proof.
\end{proof}

Lower bounds for each $\lambda_k,$ $k \in \cAp,$ for $\bfl 
\in \overR_1 \cap \overR_2$ or $\bfl \in \overR_2 \cap \overR_3$ 
may be found by determining the critical values
for the anisotropic branching random walk (for a description
of an algorithm, see \cite{lahu}, Section 3.2), whose population
dominates the one of the contact process.


\section{Growth Profile}
\setcounter{equation}{0}

This section will be devoted to the proof of Theorem
\ref{tgrowth}. Recall that $r_t$ and $R_t$ denote
the {\em smallest} and {\em largest} distances among the infected
sites $x \in A_t$ and $N_n(ns)$ denotes the number of vertices $x \in
A_{ns}$ at distance $n$ from the root that are infected at time
$ns.$

\medskip
{\bf Proof of Theorem \ref{tgrowth}.}
Recall from (\ref{twoaverages}) that
$ n_{\zeta,s} \overline{\Phi}_{2,s} = n_{\zeta} \overline{\varphi}_2
= \log \zeta$ for some $s>0.$
For each $\rho > r_u(\bfl),$ the analogous identity is
$ n_{\zeta,s} \overline{\Phi}_{\rho,s} = n_{\zeta}
 \overline{\varphi}_{\rho}.$ 
Therefore, it follows from the continuity and monotonicity 
properties of $\overline{\Phi}_{1,s;\cdot}$ and 
$ \overline{\varphi}_{1; \cdot}$ in $\bfl$ and
the fact that $ \overline{\varphi}_{1;\bfl} =0$ for each 
$\bfl \in \overR_1 \cap \overR_2,$ by Theorem \ref{weaksurvival},
along with (\ref{linearscale}) that there exists at least 
one solution $s$ of $\overps = 0.$ 

\smallskip
{\bf Proof of (\ref{tnumberinfected}).}
Suppose that $\exp \{\overps \} $ is the lead eigenvalue of the 
time-dependent equivalent of the matrix $M_1$ defined
via the entries of the matrix associated with 
$\overline{\Phi}_{2,s} .$
Fix $s>0$ so that $\overps >0.$
Then by the same arguments that led to
the result in Lemma
\ref{meanfromroot}, we find, by relying on the time-dependent
functionals, that
\begin{equation}
   \label{limtimemeanfromroot}
   \lim_{n \rightarrow \infty} \, 
      ( \, \sum_{x \in \cG_n} u_{x,ns}  \, )^{1/n}  
     =  \lim_{n \rightarrow \infty} \, ( E N_n(ns) )^{1/n} 
      = \exp \overps
\end{equation} 
for every $s>0.$
Hence, for any $\varepsilon >0$ and all sufficiently large $n>0,$
$$
   E N_n(ns)   \leq   \exp \{ n(\overps + \varepsilon) \}.
$$
Therefore, by the Borel-Cantelli lemma and the Markov inequality,
almost surely,
\begin{eqnarray*}
  \limsup_{n \rightarrow \infty}
   \, \frac{1}{n} \, \log N_n(ns) & \leq & \log 
                           ( \exp \{ \overps + \varepsilon \} +
                             \varepsilon ), \\
\mbox{thus, since $\varepsilon >0$ is arbitrary, } \qquad
\qquad \qquad \quad \qquad & & \\                              
\limsup_{n \rightarrow \infty}
   \, \frac{1}{n} \, \log N_n(ns) & \leq &  \overps. 
\end{eqnarray*} 

To prove the reverse direction, we grow some time-dependent
labelled Galton-Watson trees embedded in the set of vertices
that are infected at time $ns.$ The construction is parallel
to the one described in Proposition \ref{wmeanlimit}, from which we
conclude that 
\begin{equation}
     \overps^w =  \overps
 \end{equation}
for every $s>0.$ 
Along the same lines as before (Corollary \ref{almostmean}),
we also conclude that 
the mean offspring numbers $\mu_m$ for the
Galton-Watson trees at time $ns$
 satisfy
\begin{equation}
   \label{gwtimemeans}
         \lim_{m \rightarrow \infty} \mu_m^{1/m} = 
         \exp  \overps  . 
\end{equation}
Since $N_n(ns)$ dominates the corresponding Galton-Watson chain,
it follows that, on the event of survival of the Galton-Watson 
process,
$$
   \lim \inf_{n \rightarrow \infty} \, \frac{1}{n} \,
             \log N_n(ns) \geq \log  \exp  \overps = \overps
$$              
because by picking $m$ suitably large,
(a) the probability of the 
event of nonextinction of the Galton-Watson process
is as close as desired to the one of the event of survival of the
contact process, and (b) the expected number of offspring of
the Galton-Watson process in $\cG_n$ is as close as desired to
the expected number of vertices in $\cG_n$ that are infected at 
time $ns.$ This finishes the proof of (\ref{tnumberinfected}).

\medskip
{\bf Proof of (\ref{tsmalldist}) and  (\ref{tlargedist}).}  
We will mimic the sketch of proof in \cite{lal2}.
Note that we have seen above that 
the equation $\overps = 0$ has at least
one solution. Let $s_1$ be the smallest solution and $s_2$ the
largest solution. Moreover, for any interval $(a,b)$,
let $N_t(a,b)$ denote the number of vertices $x \in A_t$
with $ at < \vert x \vert < bt.$ 

It suffices to show that, almost surely on the event of survival, 
for any $\delta>0,$ eventually $N_t(t/s_1 + \delta t, \infty) =0$
and  $N_t(0, t/s_2 - \delta t) =0.$ Once this is verified, it will
then follow that, almost surely on the event of survival,
$\limsup_{t \rightarrow \infty} R_t/t \leq s_1$ and
$\liminf_{t \rightarrow \infty} r_t/t \geq s_2,$ which, in view
of the part of the proof we have already seen, will imply 
(\ref{tsmalldist}) and  (\ref{tlargedist}). 

Fix $\delta >0.$ 
A moment's thought shows that it is enough to consider {\em integer}
times $t.$ Observe that for each $s < s_1,$ we have 
$\overps < 0.$ Pick $\epsilon >0 $ small enough so that
$\overps + \epsilon <0 $ with $s = s_1/(1+ \delta s_1) <s_1.$ 
For sufficiently large fixed $t,$ 
the probability that $N_t(t/s_1 + \delta t, \infty) >0,$
by (\ref{limtimemeanfromroot}), 
is no larger than  
$$
  \sum_{n \geq t/s_1 + \delta t} \, \sum_{x \in \cG_n} u_{x,t}
     \leq   \sum_{n \geq t/s_1 + \delta t} \exp 
          \{ n(\overps + \epsilon) \} 
$$
with $s = s_1/(1+ \delta s_1) <s_1.$ Since $\overps + \epsilon$
is negative, the sum is bounded above by 
$  c \exp \{ - \gamma t \}$ for some positive constants $\gamma$
and $c$ (depending only on $\delta$). Since $\sum_{t=0}^{\infty}
\exp \{ - \gamma t \} < \infty,$ the Borel-Cantelli lemma 
provides that, almost surely, eventually
 $N_t(t/s_1 + \delta t, \infty)=0 .$ Since $\delta>0$ was arbitrary,
this proves one half. But in fact, the proof that, almost surely,
eventually $N_t(0, t/s_2 - \delta t) =0$ runs in parallel. 
This completes our proof. 
$\hfill \Box$

\smallskip
\bigskip
{\bf Acknowledgment.}
I would like to thank the Statistics Department of the University
of Chicago for their warm hospitality during my visit in Fall 1997,
during which a major portion of this work was completed.

\end{document}